\documentclass[12pt,a4paper]{article}

\usepackage{amsmath}    
\usepackage{amssymb}
\usepackage{graphicx}   
\usepackage{verbatim}   
\usepackage{color}      
\usepackage{hyperref}   
\usepackage{enumerate}
\usepackage{listings}
\usepackage{xcolor}
\usepackage{slashed}
\usepackage{amsthm}
\usepackage{amsmath}    
\usepackage{amssymb}
\usepackage{mathrsfs}
\usepackage{bbold}
\usepackage{lscape}

\usepackage[british]{babel}
\usepackage[font={small,it}]{caption}

\usepackage{cleveref}
\crefformat{footnote}{#2\footnotemark[#1]#3}

\newtheorem{definition}{Definition}[section]

\newtheorem{theorem}{Theorem}[section]
\newtheorem{example}{Example}[section]

\newcommand{\del}{\partial}
\renewcommand{\theta}{\vartheta}
\renewcommand{\phi}{\varphi}

\newcommand{\veccc}[3]{\left ( \begin{array}{c}#1\\#2\\#3\\ \end{array}\right )}

\newcommand{\dd}{\mathrm{d}}

\newcommand{\id}{\mathbb{1}}

\renewcommand{\vec}{\mathbf}

\newcommand{\ii}{\mathbb{i}}

\renewcommand{\title}{Stationarity preservation and the low Mach number behaviour of the Discontinuous Galerkin method on Cartesian grids}

\newcommand{\authorOne}{Wasilij Barsukow\footnote{CNRS, Bordeaux Institute of Mathematics, UMR5251, Talence, 33405 France, and Imperial College London, IRL2004, Huxley Building, South Kensington Campus, London, United Kindgom, wasilij.barsukow@math.u-bordeaux.fr}}

\begin{document}

\begin{center} \Large
\title

\vspace{1cm}

\date{}
\normalsize

\authorOne
\end{center}

\begin{abstract}
Due to added numerical stabilization (diffusion), the stationary states of numerical methods for hyperbolic problems need not be consistent discretizations of those of the PDEs. A closely related phenomenon is the lack of consistency of common finite volume methods for the Euler equations in the limit of low Mach number. In this work, the stationary states of the Discontinuous Galerkin (DG) method for linear acoustics on Cartesian grids are explored theoretically and experimentally, thus extending previous studies in the context of first-order finite difference methods. It is found that for a polynomial degree above some threshold, DG is stationarity preserving, but depending on the choice of numerical flux can suffer from a reduction of the order of accuracy at stationary state. This allows to explain the behaviour of the method for the Euler equations at low Mach number.



\end{abstract}

\newcommand{\ndof}{N_\mathrm{dof}}
\newcommand{\fluxlowmach}{Low Mach}

\section{Introduction}

\subsection{Low Mach number problems} 

Time-explicit Finite Volume methods for the compressible Euler equations employ Riemann solvers, which introduce numerical diffusion that is excessive for subsonic flows and prevents numerical simulations from resolving accurately the limit of low Mach number. 

The low Mach number limit of the Euler equations can be expressed by considering the rescaled Euler equations
\begin{subequations}
\begin{align}
 \del_t \rho + \nabla \cdot (\rho \vec v ) &= 0 \\
 \del_t (\rho \vec v) + \nabla \cdot (\rho \vec v \otimes \vec v) + \frac{\nabla p}{\epsilon^2} &= 0  \\
 \del_t e + \nabla \cdot (\vec v(e + p)) &= 0 \\
 e = \frac{p}{\gamma-1} + \frac12 \rho \epsilon^2 |\vec v|^2 &
\end{align} \label{eq:euler}
\end{subequations}
in the limit $\epsilon \to 0$. The first low-Mach-problem is that the wave speeds are $\mathcal O(\frac{1}{\epsilon})$, which results in very small time-steps when explicit methods are used. This is commented on later. However, even if one is willing to accept a long computation time, the numerical diffusion brought about by the Riemann solvers will ``dominate'' in the limit and prevent the simulation from resolving a low Mach number flow accurately. This second low-Mach-problem is not one of convergence: upon increasing the grid resolution one does obtain satisfactory results. While usually, however, grid resolution is chosen according to the smallest spatial scales that one wishes to resolve (e.g. vortices), here, the necessary resolution would be dictated by the Mach number of the flow. Even if the features of the flow remain of same size, an increase in background pressure and thus in sound speed, or a longer runtime of the simulation would necessitate a finer grid. 

The second low-Mach-problem thus is one of excessive numerical diffusion. There is a large amount of literature devoted to studies and modifications of finite volume methods to address it (e.g. \cite{turkel87,weiss95,klein95,guillard04,li08,thornber08,dellacherie10,rieper11,li13,chalons13,birken16,barsukow16,dellacherie16} and many more).

\subsection{Relation to stationary states of linear acoustics}

Linearization of the Euler equations \eqref{eq:euler} around the state of constant density and pressure and zero velocity yields (up to a linear transformation of the variables) the (linear) acoustic equations
\begin{subequations}
\begin{align}
 \del_t \vec v + \frac{\nabla p}{\epsilon^2} &= 0  \\
 \del_t p + \nabla \cdot \vec v &= 0 
\end{align} \label{eq:acoustics}
\end{subequations}
The notation $\vec v = (u,v)$ in 2-d is used below frequently. Here, $\epsilon \to 0$ is the limit of long time.
The similitude in the behaviours of numerical methods for \eqref{eq:euler} and \eqref{eq:acoustics} in the respective limits $\epsilon \to 0$ has been first considered in \cite{barsukow17a}. Therein, the behaviour of a numerical method for the Euler equations in the low Mach number limit was connected with the question whether its linearization (which is a method for linear acoustics) possesses discrete stationary states that correctly discretize those of \eqref{eq:acoustics}, coining the term \emph{stationarity preserving} for such methods. In \cite{jung22} an asymptotic two-time-scale expansion of the isentropic Euler equations has been used to show that the correct incompressible limit is only achieved discretely if the corresponding numerical method for \eqref{eq:acoustics} possesses the correct long-time limit. At PDE level, it is clear that the limit of long time for \eqref{eq:acoustics} may not even exist (in particular on periodic domains). However von-Neumann-stable linear numerical methods for \eqref{eq:acoustics} dissipate all but the stationary Fourier modes, such that the long-time limit of a stable numerical method for \eqref{eq:acoustics} will be one of its stationary states.

Discrete preservation of stationary states for the acoustic equations is by no means trivial, which is due to the stabilization added. The stationary states of a numerical method often are discretizations of only a small subset of all the stationary states of the PDE (e.g. only constants). They then poorly represent the stationary states of the PDE. This is described next in more detail.

Generally speaking, a consistent numerical method can have stationary states not consistent with the stationary states of the PDE. 
Here is an extreme example: $$\del_t q = 0 \qquad\qquad \qquad q \colon \mathbb R^+_0 \times \mathbb R \to \mathbb R$$ is a PDE which goes as far as keeping \emph{all} initial data stationary, and a numerical method with the modified equation $$\del_t q = \Delta x^p \del_x^2 q,$$ where $\Delta x$ is the spatial cell size, is consistent for $p \geq 1$. However, its steady states are only those affine in $x$, independently of $p$. 

The choice of an upwind numerical flux in a first-order method for \eqref{eq:acoustics} will result in a method, whose modified equation reads
\begin{subequations}
\begin{align}
 \del_t  u + \frac{\del_x p}{\epsilon^2} &= \Delta x \del_x^2 u\\
 \del_t  v + \frac{\del_y p}{\epsilon^2} &= \Delta y \del_y^2 v  \\
 \del_t p + \del_x u + \del_y v &= \Delta x \del_x^2 p+ \Delta y \del_y^2 p
\end{align} \label{eq:acousticsmodified}
\end{subequations}
Its stationary states are no longer $\nabla p = 0$, $\del_x u + \del_y v = 0$, but $\nabla p = 0$, $\del_x u = 0$, $\del_y v = 0$ individually, i.e. only shear flows. An initial datum of a stationary vortex, say, will be diffused away because it is not a shear flow. In \cite{barsukow17a}, numerical methods are called stationarity preserving if the relations governing their stationary states are discretizations of \emph{all} the stationary states of the PDE, without any further restriction. A more precise and testable definition is given herein based on the discrete Fourier transform of the method. This is also the tool used in this paper.

In \cite{barsukow17a} it has also been shown that a stationarity preserving method is vorticity preserving, i.e. it possesses a discretization of $\nabla \times \vec v$ kept stationary during the discrete evolution. This allowed (e.g. in \cite{barsukow20cgk,barsukow21yee}) to use involution preserving methods for linear acoustics (which by themselves have been an object of active study in e.g. \cite{morton01,sidilkover02,jeltsch06,mishra09preprint,lung14,barsukow23nodal}) to develop new low Mach number compliant methods for the Euler equations. 

Discontinuous Galerkin (DG) methods (\cite{cockburn1989a,cockburn12}) definitely make use of Riemann solvers, and in particular of the Rusanov flux, which in the context of finite volume methods would equivocally be considered the worst possible choice for low Mach number flow.
The above-mentioned large number of low Mach fixes for finite volume methods is in strong contrast to the sparsity of general results or studies for Discontinuous Galerkin methods. The modification based on ideas from \cite{turkel87} (in e.g. \cite{bassi09,bassi13}) remains the only low-Mach-fix that seems to have ever been tested for DG.  

\subsection{Time integration and grid geometry}

A first comment is due concerning the nature of the integration in time.
Time-implicit (or semi-implicit) methods allow to stabilize central derivatives in space; this approach can be seen as the most radical low-Mach modification -- to get rid of numerical diffusion altogether (e.g. \cite{viallet11,miczek15,abbate19} using fully implicit and e.g. \cite{degond07,cordier12,haack12,dimarco17,boscarino19,boscheri20,thomann20,boscheri21,boscheri21a} using semi-implicit approaches). These methods also have the advantage of alleviating the time-step restriction, since the CFL condition no longer depends on the fast waves. However, due to the increased computational effort, implicit (or semi-implicit) methods are faster than explicit ones only below Mach numbers of about $10^{-3}$. Examples of (semi-)implicit DG in the context of low Mach number include \cite{feistauer07,hennink21} among others. This paper focuses on methods that add sufficient numerical stabilization to be amenable to explicit time integration. 

A second comment concerns the grid geometry: Triangular/simplicial grids are a special case. For both finite volume (\cite{guillard09,dellacherierieper10,perrier24}) and DG methods (\cite{ern23,perrier24dg}) it is sufficient to use a contact-wave-preserving numerical flux to obtain low Mach compliance in the context of the Euler equations, or stationarity/involution preservation in the context of linear acoustics\footnote{In \cite{gjonaj07}, an erroneous argument was given to support the contrary claim.}. This explains why in the triangular-grid DG methods in \cite{bassi97,persson08} no modifications related to the accuracy of spatial discretization in the low Mach number limit are reported.

Whether a low Mach number fix is suitable or not seems also to depend on the boundary conditions, as pointed out in \cite{lannabi24}.

\subsection{Outline of the paper}

This paper aims at analyzing the low Mach number behaviour for the Euler equations/stationarity preservation for the acoustic equations of DG in the context of Cartesian grids and periodic boundaries.
First, the experimental results are presented on a test case for the Euler equations that is well-understood and often used for the analysis of finite volume methods. The finding is that DG methods with a polynomial degree of at least 1 for the Roe flux and at least 2 for the Rusanov flux are able to resolve low Mach number flow without the need for any kind of fix (Section \ref{sec:numerical1}). 

This finding is then analyzed in the context of the acoustic equations, and it is confirmed that DG methods are, indeed, stationarity preserving in these situations (Sections \ref{sec:fourier} and \ref{sec:statpres}). However, depending on the choice of the numerical flux, one might observe a decrease in the order of accuracy on the stationary states. This is the case even for the central flux; no loss of accuracy is seen when using a particular low Mach compliant numerical flux. It is important to mention, though, that the decrease in the order of accuracy is a much less dramatic problem than the inconsistency at steady state / inability to resolve the low Mach number limit altogether. The stationarity preservation analysis uses the discrete Fourier transform, in a way similar to how it has already been used for the analysis of the Active Flux method in \cite{barsukow18activeflux,barsukow24affourier}. Here, for the first time a way to analyze the order of accuracy at steady state is presented.

Objects whose natural dimension is that of space (called $d$) are set boldface; indices never denote derivatives. Matrices with $n$ rows and $m$ columns are denoted by $\mathscr M^{n \times m}(X)$ with entries in $X$. 
In two spatial dimensions, the curl and the cross product are understood as the rotated gradient and rotated scalar product, respectively. 
$\id_m$ denotes the identity matrix operating on $\mathbb R^m$ or $\mathbb C^m$ and $P^K(X)$ denotes the space of univariate polynomials of degree $\leq K$ with coefficients in $X$. $Q^K := P^{K,K}$ denotes the space of bivariate polynomials, of degree $\leq K$ in each variable.

For $d=2$, a Cartesian computational grid consists of cells
\begin{align}
 C_{ij} = \left [ x_{i-\frac12}, x_{i+\frac12} \right ] \times \left [ y_{j-\frac12}, y_{j+\frac12} \right ]
\end{align}
with $x_{i+\frac12} - x_{i-\frac12} = \Delta x$ $\forall i$, $y_{j+\frac12} - y_{j-\frac12} = \Delta y$ $\forall j$. Denote by $(x_i, y_j)$ the centroid of $C_{ij}$.

\section{Behaviour of DG in the low Mach number regime} \label{sec:numerical1}

\begin{figure}
\centering
 \includegraphics[width=0.49\textwidth]{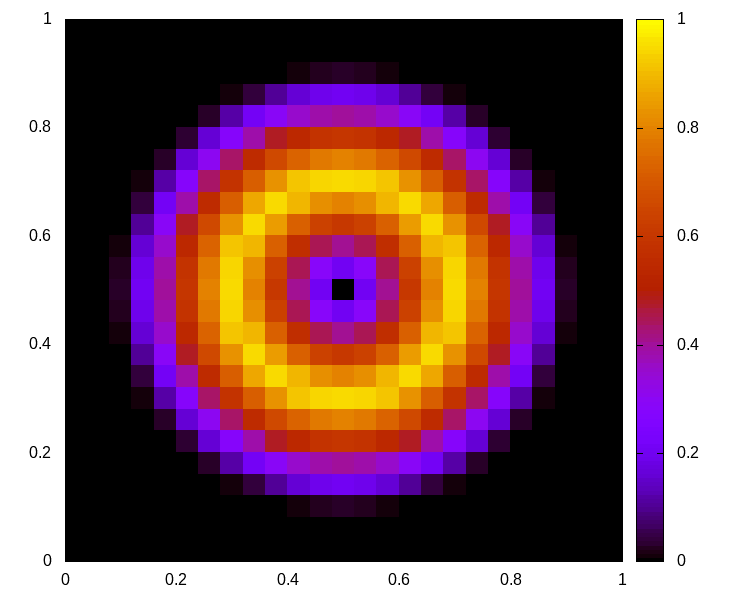} \hfill \includegraphics[width=0.49\textwidth]{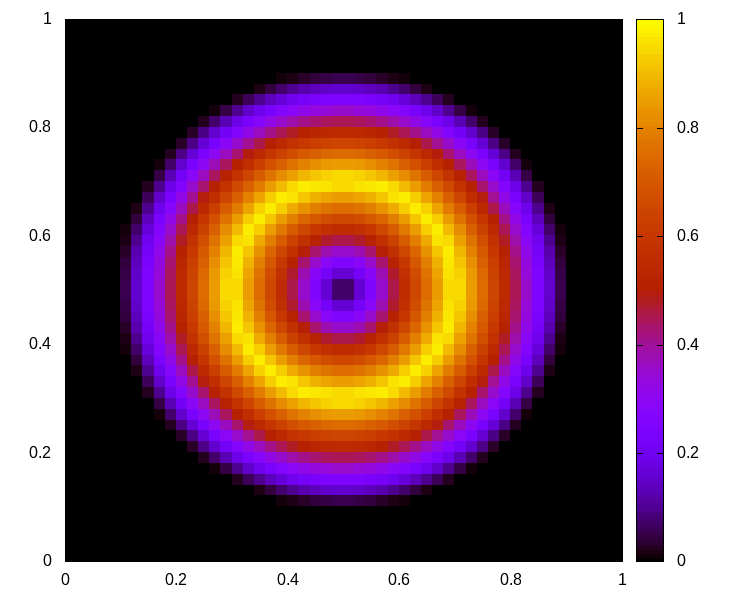}
 \caption{The low Mach number setup for the Euler equations, colour coded is the magnitude $|\vec v|$ of the velocity. \emph{Left}: Grid of $25 \times 25$. \emph{Right}: Grid of $50\times 50$.}
 \label{fig:eulersetup}
\end{figure}

Below, for the Euler equations \eqref{eq:euler} with $\gamma = 1.4$ the following setup of a stationary vortex is considered (\cite{gresho90,barsukow16}):
\begin{align}
  \rho &= 1  \qquad  v_\phi = \begin{cases}
                         	5 r & r < 0.2  \\ 2-5 r & 0.2 \leq  r < 0.4 \\ 0 &\text{else} \end{cases}
  \\p &=  \begin{cases} p_0 + 12.5 r^2 & r < 0.2 \\ p_0 + 4 \ln(5r) + 4 - 20r + 12.5 r^2 & 0.2 \leq  r < 0.4 \\ p_0 + 4 \ln 2 - 2 & \text{else}
                        \end{cases}
\end{align}
with $p_0 := \frac{1}{\gamma \epsilon^2} - \frac12$. This setup is shown in Figure \ref{fig:eulersetup}.
Contrary to questions of stationarity preservation for linear acoustics, here it is not important that the setup is a stationary state, but only that the flow is of low Mach number. In the test case at hand the Mach number is chosen via the background pressure, i.e. by modifying the sound speed: $p_0$ is chosen such that the maximum local Mach number is $\epsilon$.

This test case can be used to investigate experimentally the low Mach number compliance of a numerical method. Upon decreasing the Mach number of the flow (while keeping its time and length scales as well as the grid spacing constant) one should observe (asymptotically) no difference in the numerical solution. Note that the grid is not refined here, i.e. one is not studying the question whether the numerical solutions are in any way already close to the exact ones.

\begin{figure}
\centering
 \includegraphics[width=0.49\textwidth]{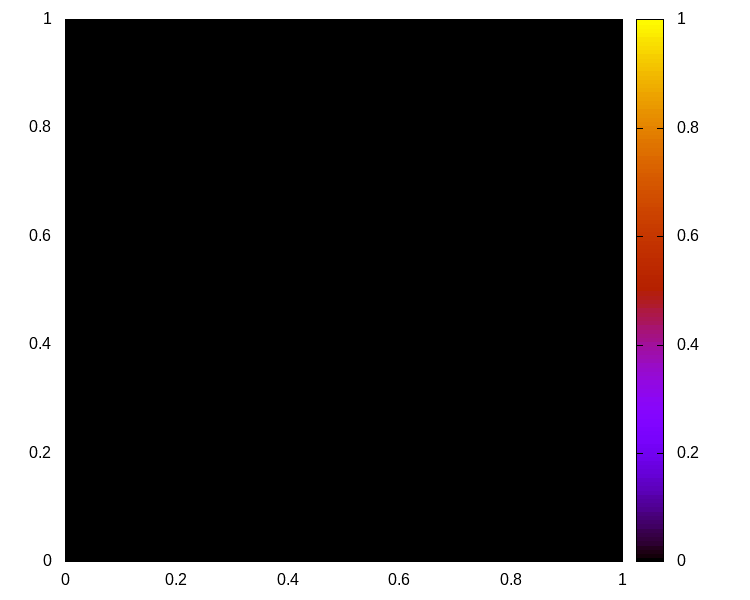} \hfill \includegraphics[width=0.49\textwidth]{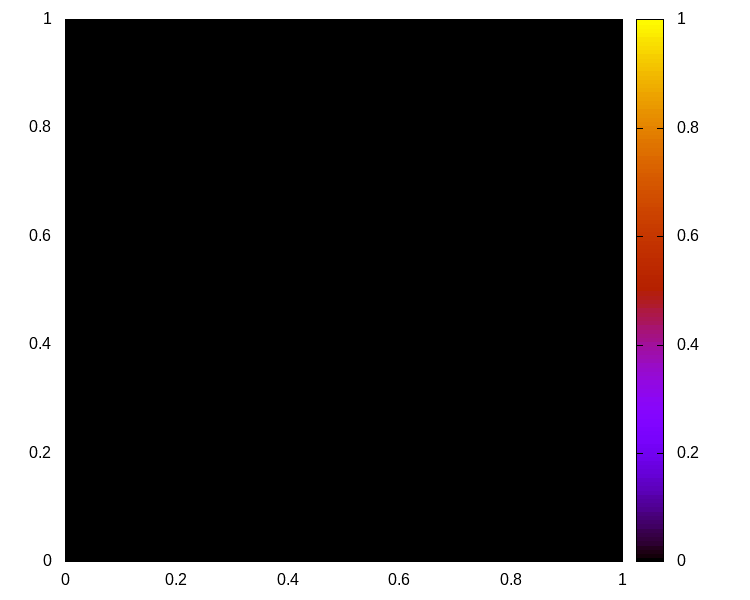}
 \includegraphics[width=0.49\textwidth]{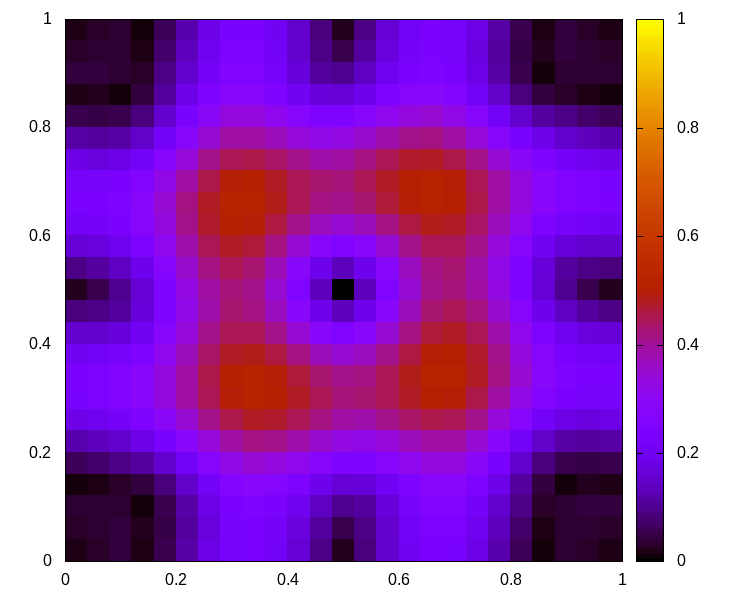} \hfill \includegraphics[width=0.49\textwidth]{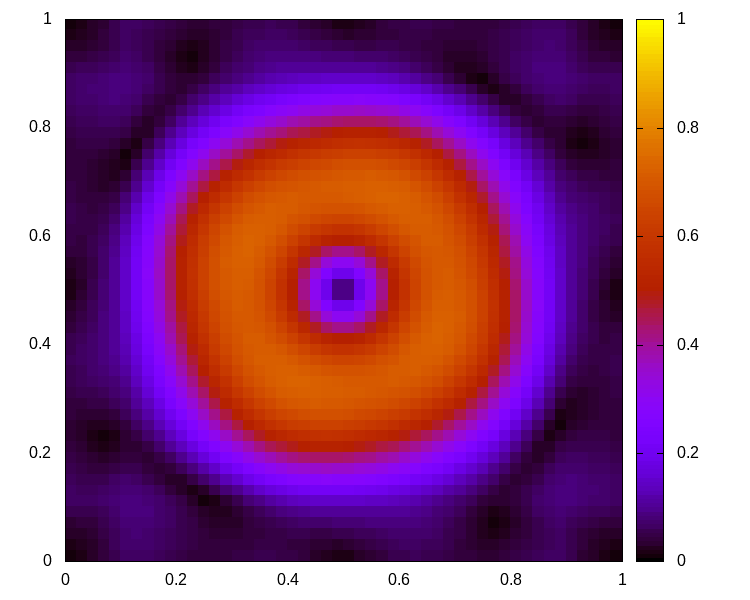} 
 \includegraphics[width=0.49\textwidth]{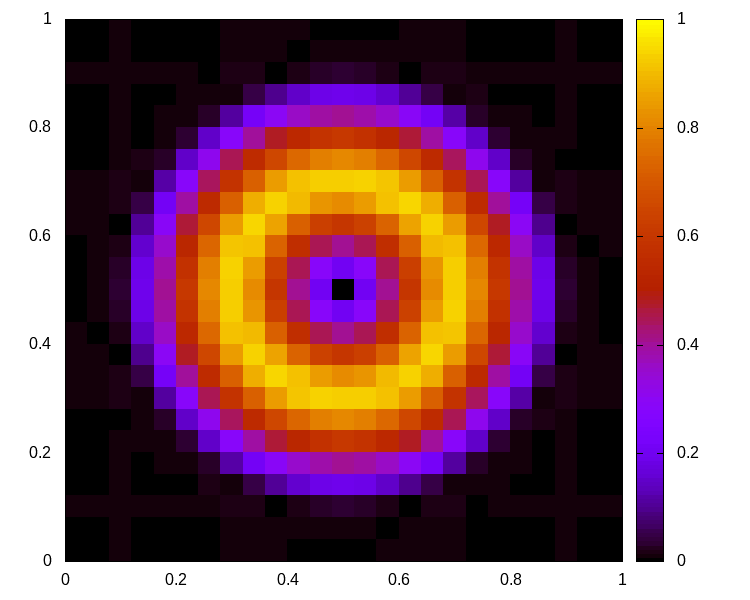} \hfill \includegraphics[width=0.49\textwidth]{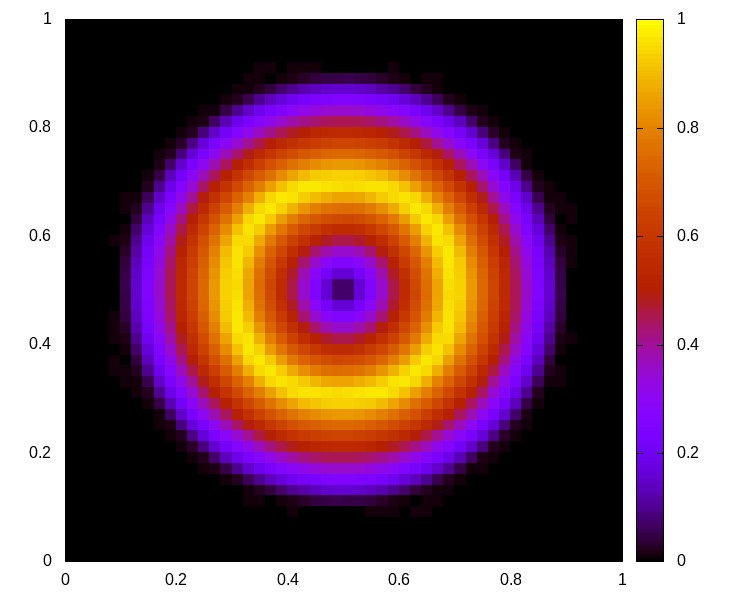} 
 \caption{Rusanov DG for the Euler equations with $\epsilon = 10^{-3}$, colour coded is the magnitude $|\vec v|$ of the velocity. Numerical solution at time $t=1$ at cell center is shown. \emph{Top}: $K=0$, the solution is diffused away down to machine zero. \emph{Middle}:  $K=1$. \emph{Bottom}: $K=2$. \emph{Left}: Grid of $25 \times 25$. \emph{Right}: Grid of $50\times 50$.}
 \label{fig:eulerM3rusanov}
\end{figure}

\begin{figure}
\centering
 \includegraphics[width=0.49\textwidth]{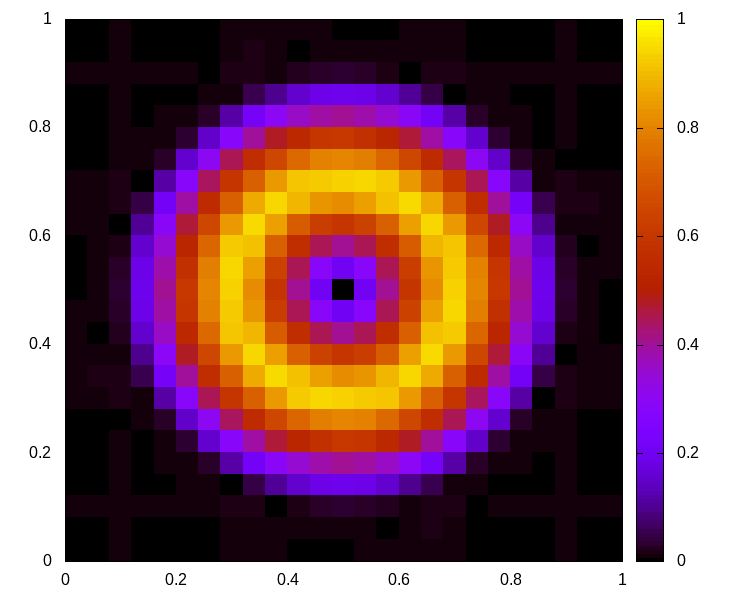} \hfill \includegraphics[width=0.49\textwidth]{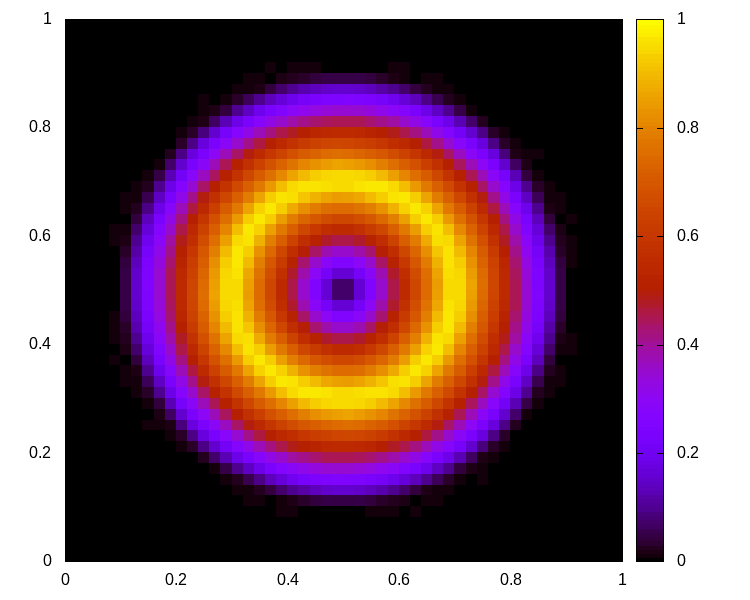} 
 \caption{Same setup as Figure \ref{fig:eulerM3rusanov}: Rusanov DG for the Euler equations for $K=2$, but with $\epsilon = 10^{-4}$.}
 \label{fig:eulerM4rusanov}
\end{figure}

Figures \ref{fig:eulerM3rusanov}--\ref{fig:eulerM4rusanov} show simulation results of DG with the Rusanov numerical flux for the Euler equations for $\epsilon = 10^{-3}$ and $10^{-4}$. One observes that for polynomial degree $K=1$, DG displays the typical artefacts associated with the lack of low Mach number compliance. For $K=2$, no such artefacts are visible. It is important to note here that this is \emph{not} just due to the higher order of accuracy, but a qualitative change of behaviour -- otherwise these artefacts would be visible again for $\epsilon = 10^{-4}$. In view of the inability of first-order finite volume methods endowed with the Rusanov flux to resolve the low Mach number limit, the good performance of DG for $K=2$ is a (pleasant) surprise.

\begin{figure}
\centering
 \includegraphics[width=0.49\textwidth]{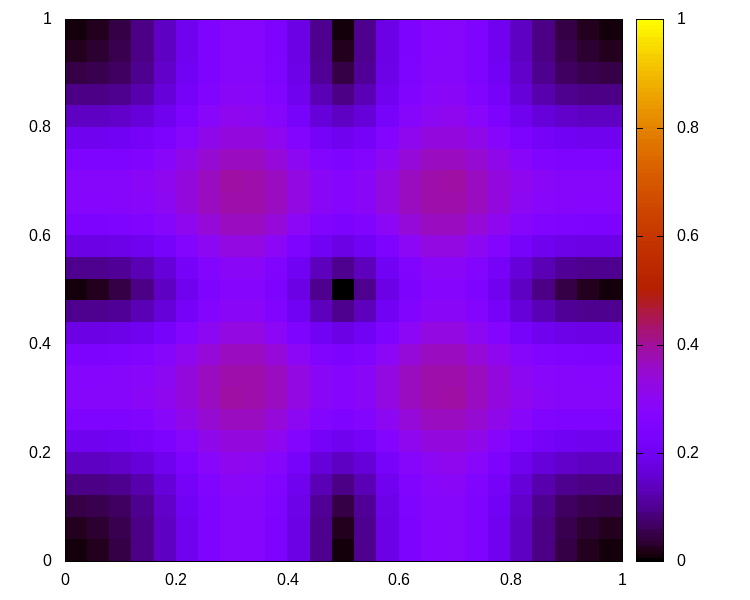} \hfill \includegraphics[width=0.49\textwidth]{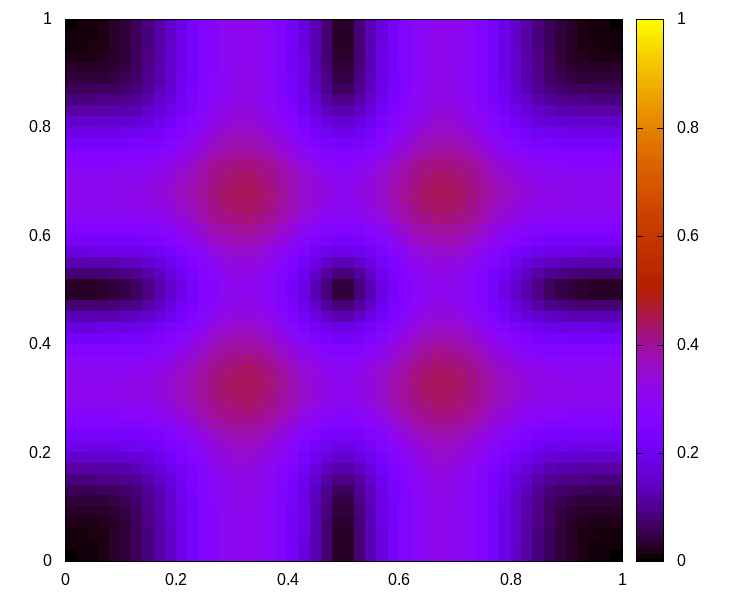}
 \includegraphics[width=0.49\textwidth]{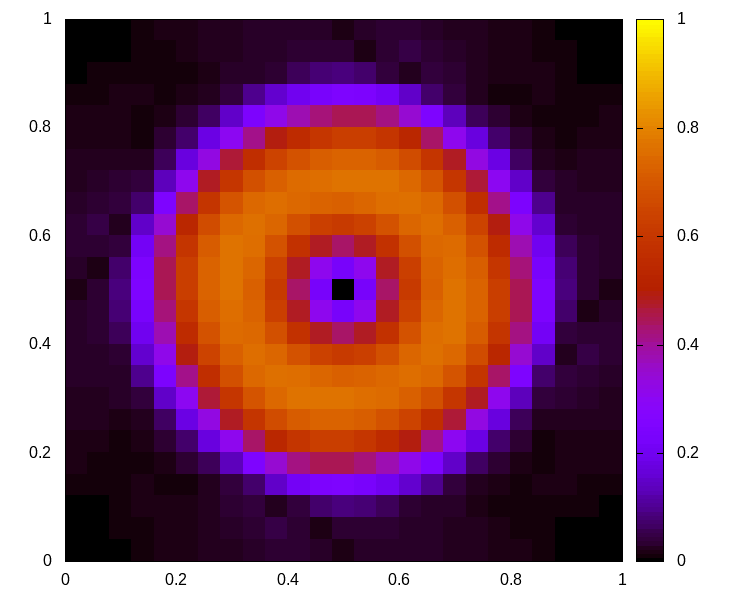} \hfill \includegraphics[width=0.49\textwidth]{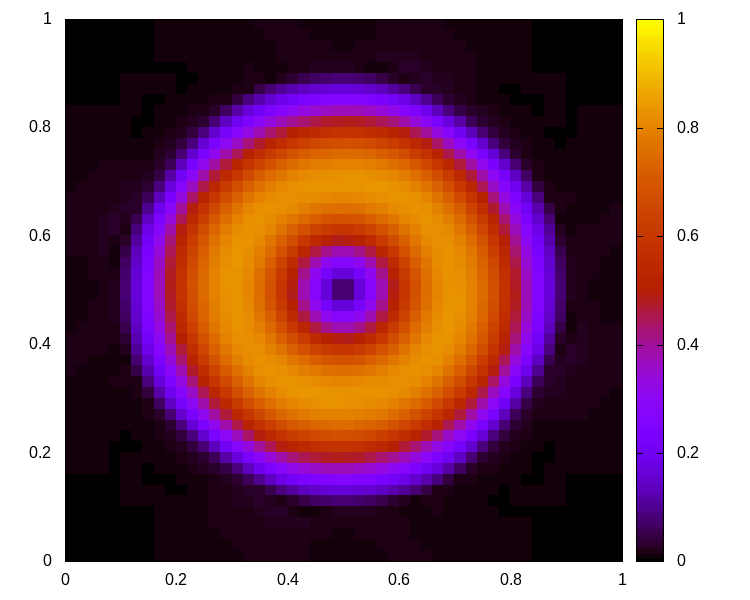} 
 \includegraphics[width=0.49\textwidth]{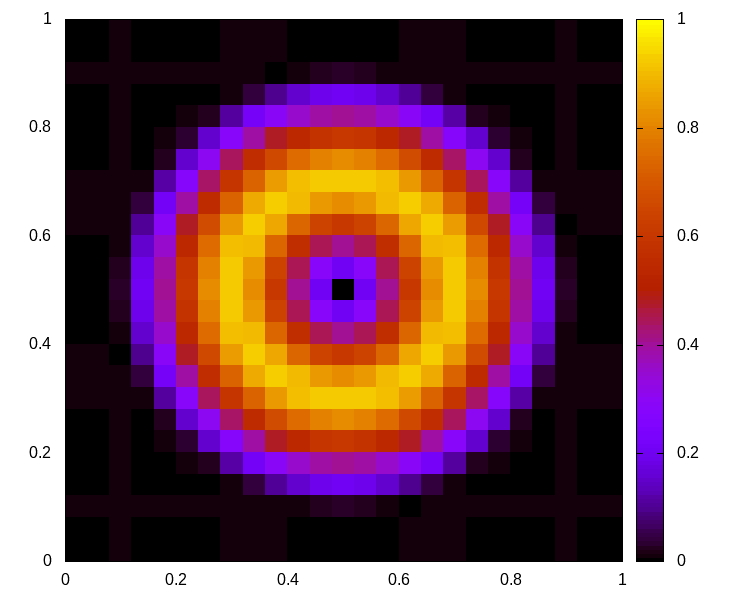} \hfill \includegraphics[width=0.49\textwidth]{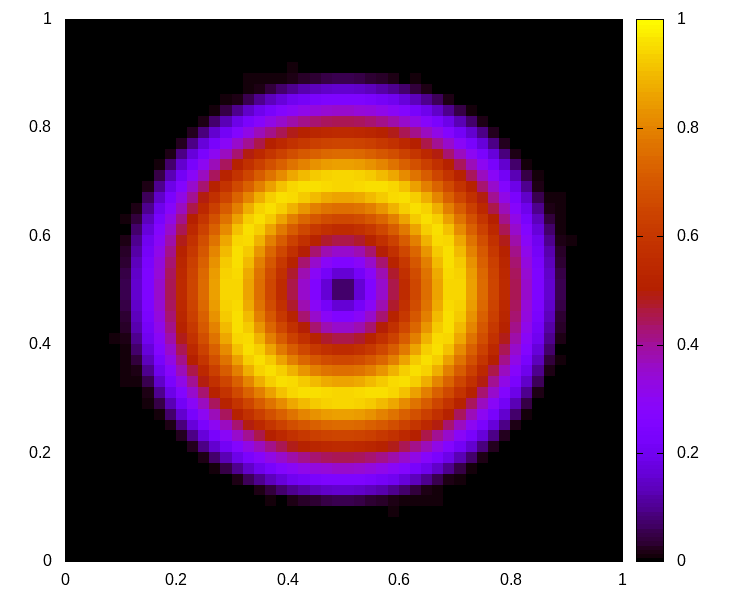} 
 \caption{Same as Figure \ref{fig:eulerM3rusanov}, but using the Roe flux.}
 \label{fig:eulerM3roe}
\end{figure}

Figure \ref{fig:eulerM3roe} shows results of the DG method with the Roe flux for $\epsilon = 10^{-3}$. One clearly observes its ability to resolve the low Mach number limit for $K=1$ and $K=2$. This conclusion follows from the fact that the solution is not diffused away leaving behind the typical cross-pattern artefacts which are a shear-flow steady state. While the solutions for both $K=1$ and $K=2$ can be considered equally ``good'' from the point of view of asymptotic consistency, their difference is due to the higher order of accuracy associated with a higher polynomial degree.

\begin{figure}
\centering
 \includegraphics[width=0.49\textwidth]{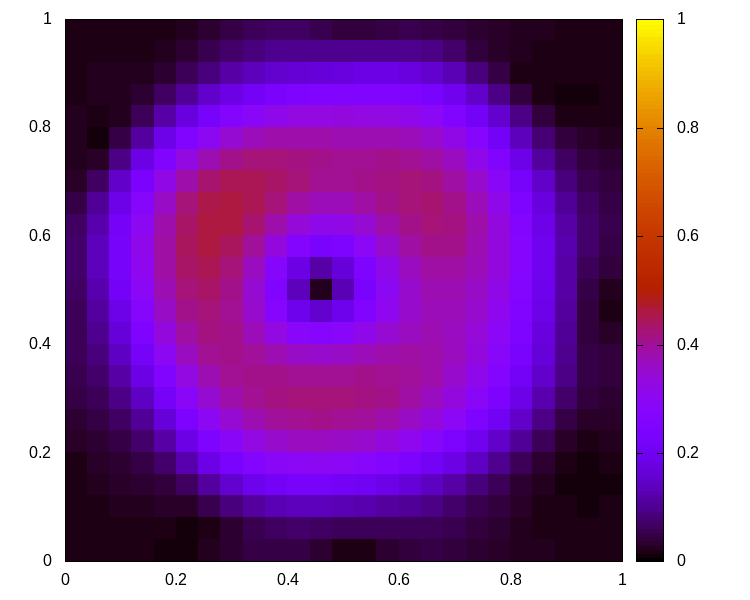} \hfill \includegraphics[width=0.49\textwidth]{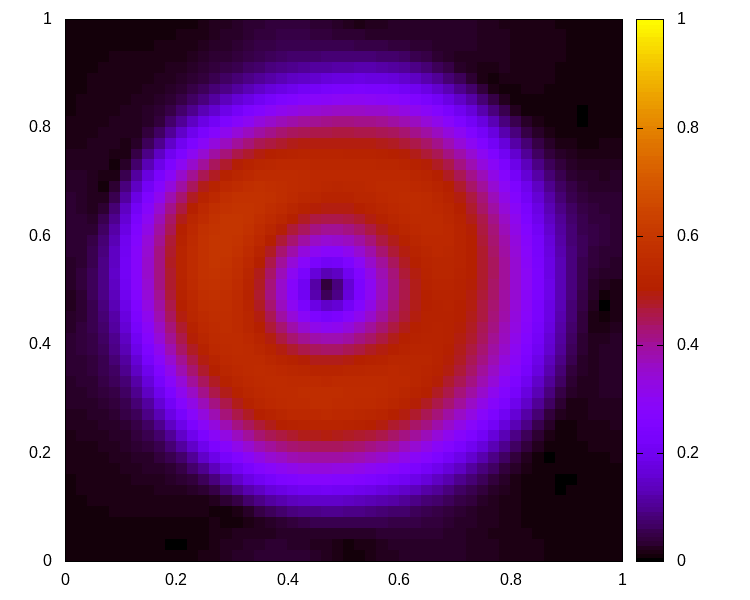}
 \includegraphics[width=0.49\textwidth]{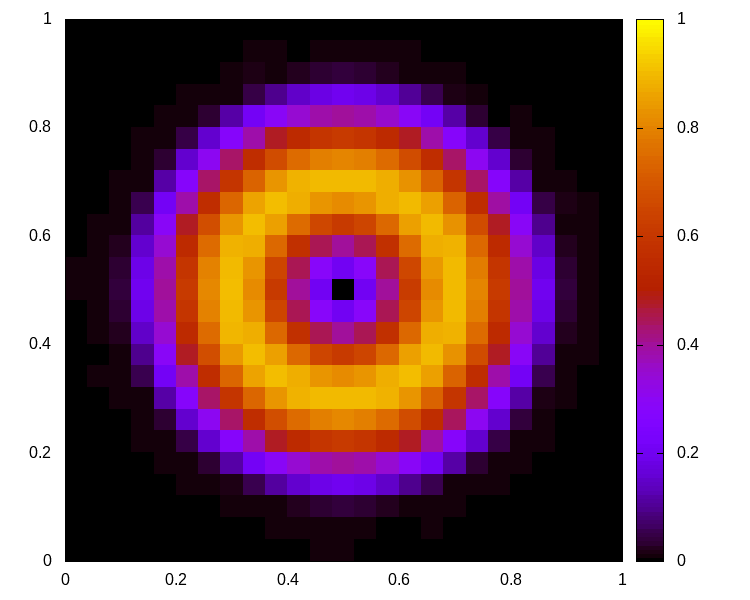} \hfill \includegraphics[width=0.49\textwidth]{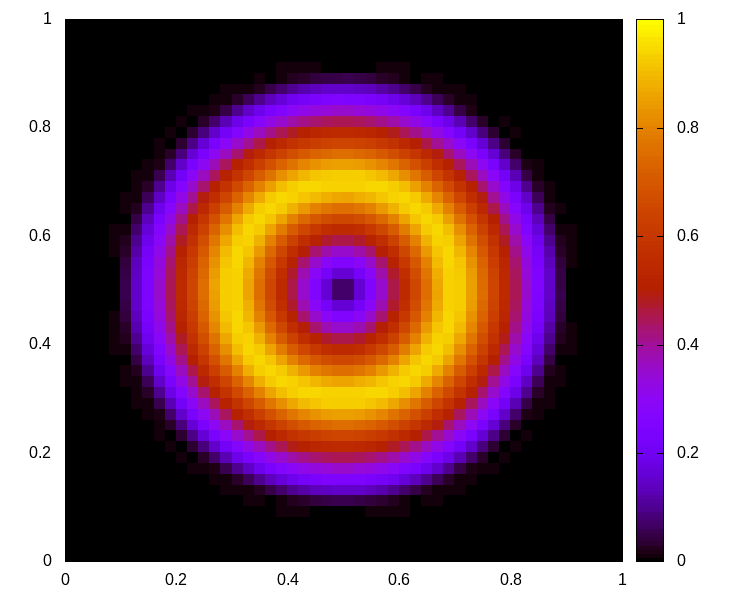} 
 \caption{Same as Figure \ref{fig:eulerM3rusanov}, but using the \fluxlowmach~flux \eqref{eq:euler0201}. \emph{Top}: $K=0$. \emph{Bottom}: $K=1$.}
 \label{fig:eulerM30201}
\end{figure}

Finally, Figure \ref{fig:eulerM30201} shows the behaviour of DG upon usage of a ``low Mach flux''. Since modifications for low Mach flows change the numerical diffusion, they can affect the stability of the method (see e.g. \cite{birken05}). In particular, certain low Mach fixes are only stable with a CFL condition that depends on $\Delta x^2$, rather than $\Delta x$. In \cite{barsukow18thesis} a careful study of the stability of low Mach fixes was performed, and the following numerical diffusion matrix (instead of the absolute value of the Jacobian in the Roe method) has been proposed:

\begin{subequations}
\begin{align}
 D_x &= \left( \begin{array}{cccc} 
               |u| + f & -1 & 0 & 0 \\
			   0 & \sqrt{u^2 + v^2} + f & 0 & 1\\
               0 & 0 & |u| + f & 0 \\
               0 & -c^2 0 & 2c
              \end{array} \right )
\\
D_y &= \left( \begin{array}{cccc}
             |v| + f & 0 & -1 & 0\\ 	
             0& |v| + f&0&0\\
			 0&0& \sqrt{u^2 + v^2} + f& 1\\
                0 & 0&-c^2  &2c 
             \end{array}\right )
\end{align} \label{eq:euler0201}
\end{subequations}

The matrix is understood in primitive variables $(\rho, u, v, p)$; for details of the derivation see \cite{barsukow18thesis}. Here, $f = 0.1$ and $c = \sqrt{\frac{\gamma p}{\rho}}$ is the speed of sound. The matrix is evaluated in the simple arithmetic mean of the two states involved. The resulting finite volume method has been found in \cite{barsukow18thesis} to have a CFL condition proportional to $\Delta x$, and to comply with the low Mach number limit. As Figure \ref{fig:eulerM30201} shows, it retains this property when applied in a DG method, with obvious improvements for higher polynomial degree due to the increased order of accuracy.

In order to explain these results, the remainder of the paper is concerned with an analysis of stationarity preservation of DG for linear acoustics. This requires to introduce the discrete Fourier transform (Section \ref{sec:fourier}) before the results of the analysis can be presented (Section \ref{sec:statpres}).

\section{The analysis of stationarity preservation on Cartesian grids} \label{sec:fourier}

\subsection{The Fourier transform} \label{sec:fouriercontinuous}

Consider a general linear hyperbolic $m \times m$ system of PDEs
\begin{align}
 \del_t q + \vec J \cdot \nabla q &= 0 & q \colon \mathbb R^+_0 \times \mathbb R^d \to \mathbb R^m \label{eq:linsystem}
\end{align}
with the directional Jacobians $\vec J = (J_1, J_2, \ldots, J_d)$, $J_k \in \mathscr M^{m \times m}(\mathbb R)$. Occasionally, the Jacobians will be denoted by $J_x, J_y, \ldots$. With the ansatz 
\begin{align}
 q(t, \vec x) = \hat q(t) \exp(\ii \vec k \cdot \vec  x)
\end{align}
of a Fourier mode \eqref{eq:linsystem} becomes 
\begin{align}
 \frac{\dd}{\dd t} \hat q + \ii (\vec J \cdot \vec k) \hat q &= 0
\end{align}
The wave vector $\vec k \in \mathbb R^d$ determines the spatial frequency of the mode; its components will be denoted by $k_1, k_2 \ldots$ or $k_x, k_y, \ldots$.
The real diagonalizability of the matrix $\vec J \cdot \vec k \in \mathscr M^{m \times m}(\mathbb R)$ is the definition of hyperbolicity of \eqref{eq:linsystem}. 

Discrete stationary states are given by the (right) kernel of $\vec J \cdot \vec k$ and discrete involutions by its left kernel (i.e. of $\ker \, (\vec J \cdot \vec k)^\text{T}$). The dimension of the kernel (i.e. the number of components of $\hat q$ such that $(\vec J \cdot \vec k) \hat q = 0$ that can be chosen freely), in general, depends on $\vec k$. In particular the choice $\vec k = 0$ is always stationary, which makes sense as then the Fourier mode is spatially constant. One can, however, only speak of an involution-preserving method if for \emph{every} $\vec k$ there exists a (row) vector $\omega$ with 
\begin{align}
 \omega (\vec J \cdot \vec k) = 0 \qquad \forall \vec k
\end{align}
From here one concludes $\frac{\dd }{\dd t} (\omega \hat q) = 0$. Note that $\omega$, in general, will depend on $k_x, k_y$.

It is helpful to illustrate this important distinction with an example. Consider linear acoustics \eqref{eq:acoustics} (with $\epsilon =1$ for simplicity) in $d=2$ spatial dimensions:
\begin{align}
 \vec J \cdot \vec k = J_x k_x + J_y k_y = \left( \begin{array}{ccc} 0 & 0 & k_x \\ 0 & 0 & k_y \\ k_x & k_y & 0 \end{array} \right)
\end{align}
Here, $\min_{\vec k} \dim \ker (\vec J \cdot \vec k) = 1$ since
\begin{align}
 (-\ii k_y,  \ii k_x, 0) \left( \begin{array}{ccc} 0 & 0 & k_x \\ 0 & 0 & k_y \\ k_x & k_y & 0 \end{array} \right) = 0 \qquad \forall \vec k
\end{align}
One thus concludes that 
\begin{align}
 (-\ii k_y,  \ii k_x, 0) \veccc{\hat u}{\hat v}{\hat p} = - \ii k_y \hat u + \ii k_x \hat v =  \ii \vec k \times \vec v
\end{align}
is stationary, and indeed this is the Fourier transform of the vorticity $\nabla \times \vec v$.

For $k_x = k_y = 0$ one obviously finds $\dim \ker (\vec J \cdot \vec k) = 3$, but such ``involutions'' would only remain stationary for particular solutions, i.e. they would lack generality. The involution $(-\ii k_y,  \ii k_x, 0) \hat q$, however, remains stationary for all data $\hat q$, even if these are time-dependent.

The dimensions of the left and right kernels being equal, it is helpful to be able to speak of a similar property when referring to the stationary states:

\begin{definition} \label{def:nontrivstatstates}
 If $\min_{\vec k} \dim \ker (\vec J \cdot \vec k) > 0$, then the system \eqref{eq:linsystem} is said to possess \emph{non-trivial stationary states}.
\end{definition}

\subsection{Cartesian grids and the discrete Fourier transform}

Linear numerical methods on Cartesian grids can be studied in great depth using the discrete Fourier transform. A usual understanding of numerical methods involves huge (but sparse) matrices whose properties are difficult to analyze theoretically in general. On Cartesian grids, the discrete Fourier transform allows to ``compress'' these matrices into a (complex-valued) matrix that has the same size whatever the number of grid cells. As shown in \cite{barsukow2025sbp}, the algebraic operations upon the discrete Fourier transform are the same as the classical approach of handling block-circulant matrices.

Consider a grid function $q \colon \mathbb Z^d \to \mathbb R^m$. The ansatz of a discrete Fourier mode $q_{ij}(t) = \hat q(t) \exp(\ii k_x \Delta x i + \ii k_y \Delta y j)$ is very useful on Cartesian grids, since shifts correspond to algebraic factors:
\begin{align}
 q_{i+I,j+J} &= q_{ij} t_x^I t_y^J
\end{align}
with $t_x := \exp(\ii k_x \Delta x)$ and $t_y := \exp(\ii k_y \Delta y)$. In this analysis, suitable (e.g. periodic) boundary conditions are tacitly assumed. On a one-dimensional grid with $N$ cells, a finite difference such as $q_{i+1} - q_i \in \mathbb R^m$, written in simple matrix notation would correspond to the $Nm \times Nm$ matrix
\begin{align}
 \left( \begin{array}{cccccc} 
         \id_{m} &  &&& -\id_{m}\\
         -\id_{m} &\id_{m}   \\
         & -\id_{m} & \id_{m} \\
         &&\ddots & \ddots & \\
         &&& -\id_{m} & \id_{m} \\
         \id_{m}&&&& -\id_{m}   \\
        \end{array}
\right ) \left( \begin{array}{c}
          q_1 \\ \vdots \\ q_i \\ q_{i+1} \\ \vdots \\ q_N
         \end{array} \right) \label{eq:matrixreal}\end{align} 
while upon the Fourier transform it merely becomes $ q_{ij} (t_x - 1) \in \mathbb C^m$, independently of $N$.

A linear finite difference method
\begin{align}
 \frac{\dd}{\dd t} q_{ij} \,\, +\!\!\! \sum_{(I,J) \in S \subset \mathbb Z^2} \!\! \alpha_{IJ} q_{i+I,j+J} &= 0 & \alpha_{IJ} &\in \mathscr M^{m \times m}(\mathbb R)
\end{align}
for \eqref{eq:linsystem} upon the discrete Fourier transform becomes
\begin{align}
 \frac{\dd}{\dd t} \hat q + \mathcal E \hat q &= 0 & \mathcal E &:= \sum_{(I,J) \in S \subset \mathbb Z^2} \alpha_{IJ} t_x^I t_y^J \in \mathscr M^{m \times m}(\mathbb C) \label{eq:fourierfindiff}
\end{align}

The method is entirely characterized by the \emph{evolution matrix} $\mathcal E$, which plays the role of $\ii (\vec J \cdot \vec k)$ (in particular observe that it depends on $t_x,t_y$ and thus on $\vec k$). However, observe also that $\mathcal E$ does not increase in size when the number of grid cells, or even when the number of spatial dimensions $d$ changes, while the analogue of the matrix in \eqref{eq:matrixreal} grows in size exponentially with $d$.

Everything that has been said about involutions and (non-trivial) stationary states in Section \ref{sec:fouriercontinuous} can now be immediately transferred to the discrete setting by analyzing the left and right kernels of $\mathcal E$. For finite difference methods, this kind of analysis has been used in \cite{barsukow17a} in order to understand the behaviour of involution preserving methods and in order to study their stationary states, with the following

\begin{definition} \label{def:statpres}
 If $\min_{\vec k} \dim \ker \mathcal E = \min_{\vec k} \dim \ker (\vec J \cdot \vec k) > 0$, then the numerical method is said to be \emph{stationarity preserving}.
\end{definition}

While it is not very difficult to construct numerical methods that preserve some given involution, it is much more difficult to find the discrete involution kept stationary by a given numerical method, or even to prove that none exist. The same is true abot the discrete stationary states. Using the discrete Fourier transform this is straightforward, since one only needs to investigate the kernel of $\mathcal E$.

\subsection{Discrete Fourier transform of the DG method}

The difference between finite difference methods and DG, or finite element methods in general, is the way how higher order of accuracy is achieved: finite differences increase the stencil, and finite element methods include additional degrees of freedom. This has an immediate consequence on the Fourier transform, because additional degrees of freedom correspond to additional Fourier modes: the discrete Fourier transform can only be applied to a grid function, and thus every degree of freedom has to be considered an independent such grid function. 

This enlarges the matrix $\mathcal E$, which now has size $(m\ndof) \times (m\ndof)$ if $\ndof$ degrees of freedom per cell are used. For example, when stationarity preservation of the Active Flux method was studied in \cite{barsukow18activeflux,barsukow24affourier}, the matrix to be analyzed was $12 \times 12$ because linear acoustics in 2-d has three variables, and because the Active Flux method under consideration had 4 degrees of freedom per cell. This, however, is still much smaller than the $30000 \times 30000$ matrix that would be involved for the very modest $50 \times 50$ grid, if written down in the style of \eqref{eq:matrixreal}, i.e. without resorting to the Fourier transform. While the latter matrix is sparse, the matrices involved in a numerical method's Fourier transform are dense, and also become increasingly difficult to analyze theoretically as their size grows with the order of accuracy of the method. This is why the analysis in the next Section is only performed for $K \leq 3$, which involves analyzing kernels of matrices up to $48 \times 48$.

The Fourier transform of the DG method shall be derived next in detail. Even though the reader might be familiar with the DG method, the derivation is useful for introducing the notation and, most of all, to show how the Fourier transform acts on the expressions. Consider basis functions $\{ b_k \}_{k}$, $b_k \colon [-\frac12,\frac12] \to \mathbb R$ that form a basis of $P^K$, $K \geq 0$ and a conservation law 
\begin{align}
\del_t q + \del_x f^x (q) + \del_y f^y(q) &= 0 & q &\colon \mathbb R^+_0 \times \mathbb R^2 \to \mathbb R^m \label{eq:conslaw}
\end{align}
In a tensor-basis framework, the multi-dimensional test and basis functions are products of the one-dimensional ones in different directions. Consider $q$ thus to be approximated by $q_h$ in the space
\begin{align}
 q_h &\in (P^{K,K}_\text{br})^m & P^{K,K}_\text{br} &= \big\{ \pi \in L^\infty, \pi \big|_{C_{ij}} \in P^{K,K}(\mathbb R) \big\}
\end{align}
Define the restriction $q_h \big|_{C_{ij}}$ on one cell (including a coordinate shift)
\begin{align}
 q_{ij} &\colon  \left[-\frac{\Delta x}{2}, \frac{\Delta x}{2}\right] \times \left[-\frac{\Delta y}{2}, \frac{\Delta y}{2}\right] \to \mathbb R^m \\ 
q_{ij}(t, x, y) &:= q\big|_{C_{ij}}(t, x_i + x, y_j + y) 
\end{align}
and expand
\begin{align}
q_{ij}(t, x, y) = \sum_{\alpha,\beta = 0}^{K} q^{(\alpha \beta)}_{ij}(t) b_\alpha\left(\frac{x}{\Delta x}\right) b_\beta\left(\frac{y}{\Delta y}\right) \label{eq:expansionreal}
\end{align}
where $q^{(\alpha \beta)}_{ij}(t) \in \mathbb R^m$ are the degrees of freedom in cell $(i,j)$.

Multiply \eqref{eq:conslaw} with a test function, integrate over one cell and by parts:
\begin{align}
 \label{eq:dgfirstform}
 \sum_{\alpha\beta} \frac{\dd}{\dd t} q^{(\alpha \beta)}_{ij}(t) \int_{-\frac{\Delta x}2}^{\frac{\Delta x}2} \dd x\, b_k\left(\frac{x}{\Delta x}\right) b_\alpha\left(\frac{x}{\Delta x}\right) \int_{-\frac{\Delta y}2}^{\frac{\Delta y}2} \dd y\,  b_\ell\left(\frac{y}{\Delta y}\right) b_\beta\left(\frac{y}{\Delta y}\right)  
 \\\nonumber - \int_{-\frac{\Delta x}2}^{\frac{\Delta x}2} \dd x\, \int_{-\frac{\Delta y}2}^{\frac{\Delta y}2} \dd y\, \Big (b_\ell\left(\frac{y}{\Delta y}\right) \del_x b_k\left(\frac{x}{\Delta x}\right) f_x(q) + b_k\left(\frac{x}{\Delta x}\right) \del_y b_\ell\left(\frac{y}{\Delta y}\right)  f_y(q)  \Big) 
 \\\nonumber+ b_k\left(\frac12\right) \int_{-\frac{\Delta y}2}^{\frac{\Delta y}2} \dd y\,  b_\ell\left(\frac{y}{\Delta y}\right) \hat f_x\left(q_{ij}\left(t, \frac{\Delta x}{2}, y\right), q_{i+1,j}\left(t, - \frac{\Delta x}{2}, y\right)\right)
 \\\nonumber- b_k\left(-\frac12\right) \int_{-\frac{\Delta y}2}^{\frac{\Delta y}2} \dd y\,  b_\ell\left(\frac{y}{\Delta y}\right) \hat f_x\left(q_{i-1,j}\left(t, \frac{\Delta x}{2}, y\right), q_{i,j}\left(t, - \frac{\Delta x}{2}, y\right)\right)
 \\\nonumber + b_\ell\left(\frac12\right) \int_{-\frac{\Delta x}2}^{\frac{\Delta x}2} \dd x\, b_k\left(\frac{x}{\Delta x}\right)\hat f_y\left(q_{ij}\left(t, x,\frac{\Delta y}{2}\right), q_{i,j+1}\left(t,x, - \frac{\Delta y}{2}\right)\right)
 \\\nonumber - b_\ell\left(-\frac12 \right) \int_{-\frac{\Delta x}2}^{\frac{\Delta x}2} \dd x\, b_k\left(\frac{x}{\Delta x}\right) \hat f_y\left(q_{i,j-1}\left(t, x,\frac{\Delta y}{2}\right), q_{ij}\left(t,x, - \frac{\Delta y}{2}\right)\right)
 &= 0 
\end{align}
Here, numerical fluxes $\hat f_x$, $\hat f_y$ have been introduced, assumed linear in the following. This is the DG method.
As usual, define the mass matrix $M \in \mathscr M^{K+1,K+1}(\mathbb R)$ as
\begin{align}
 M^k_{\alpha} &:= \frac{1}{\Delta x}\int_{-\frac{\Delta x}2}^{\frac{\Delta x}2} \dd x\, b_k\left(\frac{x}{\Delta x}\right) b_\alpha\left(\frac{x}{\Delta x}\right) = \int_{-\frac{1}2}^{\frac{1}2} \dd \xi\, b_k\left(\xi\right) b_\alpha\left(\xi\right)\\
 &= \frac{1}{\Delta y}\int_{-\frac{\Delta y}2}^{\frac{\Delta y}2} \dd y\,  b_\ell\left(\frac{y}{\Delta y}\right) b_\beta\left(\frac{y}{\Delta y}\right) 
\end{align}

For the numerical examples below, an orthonormal basis $\{ b_k \}_k$
\begin{align}
 \{ b_k \}_k = \left\{ 1, 2 \sqrt{3} x , \frac{\sqrt{5}}{2} (12 x^2 -1 ), \sqrt{7} x (20 x^2 - 3), \dots \right \}
\end{align}
is used, giving rise to a diagonal mass matrix
\begin{align}
 M^k_{\alpha}  = \delta_{k\alpha} 
\end{align}
Any other choice of basis amounts merely to a linear transformation that will not change any of the statements of interest.
Define also
\begin{align}
 D_\alpha^k &= \int_{-\frac{\Delta x}2}^{\frac{\Delta x}2} \dd x \, \del_x b_k\left(\frac{x}{\Delta x} \right) b_\alpha\left(\frac{x}{\Delta x}\right)  =  \int_{-\frac12}^{\frac12} \dd \xi \, b_k'(\xi) b_\alpha(\xi) \in \mathbb R
\end{align}

Apply now the discrete Fourier transform, e.g. write, for any $\alpha,\beta$ 
\begin{align}
 q_{ij}^{(\alpha \beta)}(t) = \hat q^{(\alpha\beta)} \exp(\ii k_x \Delta x i + \ii k_y \Delta y j) = \hat q^{(\alpha\beta)} t_x^i t_y^j
\end{align}
which via \eqref{eq:expansionreal} also implies
\begin{align}
q_{ij}(t, x, y) = t_x^i t_y^j\sum_{\alpha,\beta = 0}^{K} \hat q^{(\alpha \beta)}(t)  b_\alpha\left(\frac{x}{\Delta x}\right) b_\beta\left(\frac{y}{\Delta y}\right)
\end{align}

Use that, for linear hyperbolic systems, $f^x(q) = J^x q$ and $f^y(q) = J^y q$: {\footnotesize
\begin{align}
  \Delta x \Delta y\sum_{\alpha, \beta} M^k_{\alpha} M^\ell_{\beta} \frac{\dd}{\dd t}  q_{ij}^{(\alpha \beta)}(t)   
  -J^x \Delta y\sum_{\alpha, \beta} M^\ell_{\beta}   D^k_\alpha   q_{ij}^{(\alpha \beta)}(t)
  -J^y \Delta x\sum_{\alpha, \beta} M^k_{\alpha}   D^\ell_\beta  q_{ij}^{(\alpha \beta)}(t)
 \\\nonumber+ \left( b_k\left(\frac12\right) - b_k\left(-\frac12\right)\frac{1}{t_x} \right) \int_{-\frac{\Delta y}2}^{\frac{\Delta y}2} \dd y\,  b_\ell\left(\frac{y}{\Delta y}\right) 
 \hat f^x\left(q_{ij}\left(t, \frac{\Delta x}{2}, y\right), t_x q_{ij}\left(t, - \frac{\Delta x}{2}, y\right)\right)
 \\\nonumber + \left( b_\ell\left(\frac12\right) - b_\ell\left(-\frac12 \right) \frac{1}{t_y} \right) \int_{-\frac{\Delta x}2}^{\frac{\Delta x}2} \dd x\, b_k\left(\frac{x}{\Delta x}\right) \hat f^y\left(q_{ij}\left(t, x,\frac{\Delta y}{2}\right), t_y q_{ij}\left(t,x, - \frac{\Delta y}{2}\right)\right)
 &= 0 
\end{align}}

Write the numerical flux as
\begin{align}
 \hat f^x(q_\text{L}, q_\text{R}) = J_x \frac{q_\text{L} + q_\text{R}}{2} - D_x \frac{q_\text{R} - q_\text{L}}{2} \label{eq:numfluxx}\\
 \hat f^y(q_\text{L}, q_\text{R}) = J_y \frac{q_\text{L} + q_\text{R}}{2} - D_y \frac{q_\text{R} - q_\text{L}}{2} \label{eq:numfluxy}
\end{align}
with diffusion matrices $D_x,D_y \in \mathscr M^{m \times m}(\mathbb R)$. They will be specified later; recall that for linear acoustics, with the ordering $\vec v = (u, v), p$ of variables,
\begin{align}
 J_x &= \left( \begin{array}{ccc}  0 & 0 & 1 \\ 0 & 0 & 0\\ 1 & 0 &0  \end{array} \right )
 & J_y &= \left( \begin{array}{ccc}  0 & 0 & 0 \\ 0 & 0 &1 \\0 & 1 &0  \end{array} \right )
\end{align}
Observe that now the common factor $t_x^i t_y^j$ cancels out:

 {\footnotesize
\begin{align}
 \sum_{\alpha, \beta} M^k_{\alpha} M^\ell_{\beta}  \frac{\dd}{\dd t} \hat q^{(\alpha \beta)}(t) 
  -J^x \frac{1}{\Delta x} \sum_{\alpha, \beta} M^\ell_{\beta}  D^k_\alpha \hat q^{(\alpha \beta)}(t)
-J^y \frac{1}{\Delta y} \sum_{\alpha, \beta} M^k_{\alpha}  D^\ell_\beta\hat q^{(\alpha \beta)}(t)
 \\\nonumber+ \sum_{\alpha, \beta}  \frac{ b_k\left(\frac12\right)t_x - b_k\left(-\frac12\right)}{t_x \Delta x}   
\left( J^x \frac{ b_\alpha\left(\frac{1}{2}\right) + t_x b_\alpha\left( - \frac{1}{2}\right) }{2} - D^x \frac{t_x b_\alpha\left( - \frac{1}{2}\right) - b_\alpha\left(\frac{1}{2}\right) }{2} \right) M^\ell_{\beta} \hat q^{(\alpha \beta)}(t)
 \\\nonumber + \sum_{\alpha, \beta} \frac{ b_\ell\left(\frac12\right)t_y - b_\ell\left(-\frac12 \right)}{t_y \Delta y}  \left ( J^y \frac{b_\beta\left(\frac{1}{2}\right) + t_y b_\beta\left(- \frac{1}{2}\right)}{2} - D^y \frac{ t_y b_\beta\left(- \frac{1}{2}\right) - b_\beta\left(\frac{1}{2}\right)}{2} \right ) M^k_{\alpha} \hat q^{(\alpha \beta)}(t) \label{eq:DGfourier}
 &= 0 
\end{align}}

This can be written as
\begin{align}
 \sum_{\alpha, \beta} M^k_{\alpha} M^\ell_{\beta}   \frac{\dd}{\dd t} \hat q^{(\alpha \beta)}(t) + \sum_{\alpha, \beta} \mathcal E_{\alpha \beta}^{k \ell} \hat q^{(\alpha \beta)} = 0
\end{align}
Observe that $M^k_{\alpha} \in \mathbb R$, but $\mathcal E_{\alpha \beta}^{k \ell} \in \mathscr M^{m \times m}(\mathbb C)$: 
\begin{align}
  \mathcal E_{\alpha \beta}^{k \ell} &:= -J^x \frac{1}{\Delta x}  M^\ell_{\beta}  D^k_\alpha
 -J^y \frac{1}{\Delta y} M^k_{\alpha}  D^\ell_\beta 
 \\\nonumber&+  \frac{ b_k\left(\frac12\right)t_x - b_k\left(-\frac12\right)}{t_x \Delta x}   
\left( J^x \frac{ b_\alpha\left(\frac{1}{2}\right) + t_x b_\alpha\left( - \frac{1}{2}\right) }{2} - D^x \frac{t_x b_\alpha\left( - \frac{1}{2}\right) - b_\alpha\left(\frac{1}{2}\right) }{2} \right)M^\ell_{\beta}
 \\\nonumber &+ \frac{ b_\ell\left(\frac12\right)t_y - b_\ell\left(-\frac12 \right)}{t_y \Delta y}  \left ( J^y \frac{b_\beta\left(\frac{1}{2}\right) + t_y b_\beta\left(- \frac{1}{2}\right)}{2} - D^y \frac{ t_y b_\beta\left(- \frac{1}{2}\right) - b_\beta\left(\frac{1}{2}\right)}{2} \right ) M^k_\alpha
\end{align}

If $\{ b_k \}_k$ are taken as a basis for the space $P^K$ of polynomials with degree at most $K$, then the set of pairs $\{ (k,\ell) \}_{k,\ell = 0, \ldots, K}$ has $(K+1)^2$ elements. Each $\hat q^{(k\ell)}$ is a vector in $\mathbb C^m$. One can thus rewrite \eqref{eq:DGfourier} as
\begin{align}
 \frac{\dd}{\dd t} Q + \mathcal E Q &= 0
\end{align}
where $Q$ combines all $m(K+1)^2$ entries of the vectors $\{ \hat q^{(k\ell)} \}_{k,\ell = 0, \ldots, K} $ and $\mathcal E$ consists of all the entries of all the matrices $\mathcal E_{\alpha \beta}^{k \ell}$. The precise form of $\mathcal E \in \mathscr M^{m(K+1)^2 \times m(K+1)^2}(\mathbb C)$ of course depends on the ordering that one chooses. In the end, the method can thus be brought in the same form \eqref{eq:fourierfindiff} as previously the finite difference method.

\section{Stationarity preservation of DG and the order of accuracy of stationary solutions} \label{sec:statpres}

\subsection{General procedure}

In general, the kernel of $\mathcal E$ is multi-dimensional, and it is difficult to make sense of individual basis elements, since their choice and their normalizations are arbitrary. For finite difference methods, a natural definition of stationarity preservation (see Definition \ref{def:statpres}) is a one-to-one correspondence between elements in $\ker \vec J \cdot \vec k$ and those in $\ker \mathcal E$, in the sense that the latter converge to the former as $\Delta x, \Delta y \to 0$. DG achieves higher order of accuracy through incorporation of additional degrees of freedom inside the cell: while finite difference methods represent every scalar function by one value per cell, DG represents it through $\ndof =(K+1)^2$ values per cell (in 2-d). Thus, while elements of $\ker \vec J \cdot \vec k$ live in $\mathbb C^m$, those of $\ker \mathcal E$ for DG are in $\mathbb C^{m(K+1)^2}$ and cannot converge to those of $\ker \vec J \cdot \vec  k$. The number of elements in $\ker \mathcal E$ also can be larger than $\dim \ker \vec J \cdot \vec k$ for the same $\vec k$: e.g. for central DG for linear acoustics $\min_{\vec k} \dim \ker \mathcal E = (K+1)^2$ while $\min_{\vec k} \dim \ker \vec J \cdot \vec k = 1$. One observes that the dimension is larger by a factor which equals the number of degrees of freedom per cell, which itself is directly linked to the order of accuracy of the method. For other choices of the flux the kernel can be smaller, e.g. for upwind DG, $\min_{\vec k} \dim \ker \mathcal E$ is only $K^2$. One is thus led to think that this reduction might manifest itself in a reduction of the order of accuracy at stationary state. However, it turns out that this is not exactly the case, and to make the details of the connection precise is the aim of this Section. To begin with, the following definition formulates the least requirement:

\begin{definition}
Consider a linear numerical method for the system \eqref{eq:linsystem} with $\min_{\vec k} \dim \ker \vec J \cdot \vec k = 1$. Assume that it has $\ndof$ degrees of freedom per cell per variable. If 
\begin{align}
 1 \leq \min_{\vec k} \dim \ker \mathcal E \leq  \ndof 
\end{align}
then the numerical method is said to be \emph{stationarity preserving}.
\end{definition}
For problems with $\min_{\vec k} \dim \ker \vec J \cdot \vec k \geq 2$, outside of the scope of this paper, the definition will require refinement.

\newcommand{\dof}{\mathrm{DOF}}

First, let us define operators that allow to go back and forth between $\mathbb C^m$ and $\mathbb C^{m(K+1)^2}$.

\begin{definition}
 \begin{itemize}
  \item $\dof_{ij} \colon  (L^\infty(\mathbb C))^m \to \mathbb C^{m(K+1)^2}$ is the operator that evaluates the degrees of freedom\footnote{or, if the meaning of a degree of freedom is that of 1-forms on the dual space of $\mathbb R[x,y]^m$, as is customary in finite element literature, then $\dof$ \emph{are} the degrees of freedom, naturally extended to act on $(L^\infty)^m$.} in cell $i,j$ of a function $q$. For DG,
\begin{align}
 \sum_{\alpha\beta} M^k_\alpha M^\ell_\beta \dof_{ij}(q)_{\alpha\beta} = \frac{1}{\Delta x} \frac{1}{\Delta y} \int_{-\frac{\Delta x}{2}}^{\frac{\Delta x}{2}} \dd x \int_{-\frac{\Delta y}{2}}^{\frac{\Delta y}{2}} \dd y \, b_k\left( \frac{x}{\Delta x} \right ) b_\ell\left( \frac{y}{\Delta y} \right ) q(x_i + x, y_j+y)
\end{align}
  \item $\mathbb M^m_{\vec k} := \{ \hat Q \exp(\ii \vec k \cdot \vec x) \text{ with } \hat Q \in \mathbb C^m \}$ is the set of Fourier modes with $\vec k \in \mathbb R^d$. These are the natural objects with which to compare the discrete Fourier modes $\hat q \exp(\ii k_x \Delta x i + \ii k_y \Delta y j)$ with $\hat q \in \mathbb C^{m(K+1)^2}$ introduced in the previous Sections.
  \item We denote by $R \colon \mathbb C^{m(K+1)^2} \to \mathbb C[x,y]^m$ the \emph{reconstruction operator},
  defined by
  \begin{align}
   q \mapsto R q = \sum_{\alpha,\beta} q^{(\alpha\beta)} b_\alpha \left( \frac{x}{\Delta x} \right ) b_\beta \left( \frac{y}{\Delta y} \right )
  \end{align}
 \end{itemize}
\end{definition}

Recall that for the orthonormal basis considered here, $M^k_\alpha = \delta_{k\alpha}$.
Consider now what happens if one evaluates the degrees of freedom of a Fourier mode $q \in \mathbb M^m_{\vec k}$:
\begin{align*}
  \sum_{\alpha\beta} M^k_\alpha M^\ell_\beta \dof_{ij}(q)_{\alpha\beta} &= \\
  &\nonumber \!\!\!\!\!\!\!\!\!\!\!\!\!\!\!\!\!\!\!\!\!\!\!\! \!\!\!\!\!\!\!\!\!\!\!\!\!\!\!\!\!\!\!\!\!\!\!\!\!\!\!\!\!\!\!\!\!\!\frac{1}{\Delta x} \frac{1}{\Delta y} \int_{-\frac{\Delta x}{2}}^{\frac{\Delta x}{2}} \dd x \int_{-\frac{\Delta y}{2}}^{\frac{\Delta y}{2}} \dd y \,  b_k\left( \frac{x}{\Delta x} \right ) b_\ell\left( \frac{y}{\Delta y} \right ) \hat Q \exp(\ii k_x (x_i + x)) \exp(\ii k_y (y_j + y)) \\
 & \!\!\!\!\!\!\!\!\!\!\!\!\!\!\!\!\!\!\!\!\!\!\!\!\!\!\!\!\!\!\!\!\!\!\!\!\!\!\!\!\!\!\!\!\!\!\!\! = t_x^i t_y^j \hat Q  \frac{1}{\Delta x} \frac{1}{\Delta y} \int_{-\frac{\Delta x}{2}}^{\frac{\Delta x}{2}} \dd x\int_{-\frac{\Delta y}{2}}^{\frac{\Delta y}{2}} \dd y \,  b_k\left( \frac{x}{\Delta x} \right ) b_\ell\left( \frac{y}{\Delta y} \right ) \exp(\ii k_x  x) \exp(\ii k_y  y)
\end{align*}
It hence makes sense to make the following
\begin{definition}
$\widehat{\dof} \colon \mathbb C^m \to \mathbb C^{m(K+1)^2}$ is the operator that evaluates the discrete Fourier transform of a Fourier mode $q = \hat Q \exp(\ii \vec k \cdot \vec x)$ as {\footnotesize
\begin{align}
  \sum_{\alpha\beta} M^k_\alpha M^\ell_\beta \widehat{\dof}(\hat Q)_{\alpha\beta} := \hat Q  \frac{1}{\Delta x} \frac{1}{\Delta y} \int_{-\frac{\Delta x}{2}}^{\frac{\Delta x}{2}} \dd x \int_{-\frac{\Delta y}{2}}^{\frac{\Delta y}{2}} \dd y \, b_k\left( \frac{x}{\Delta x} \right ) b_\ell\left( \frac{y}{\Delta y} \right ) \exp(\ii k_x  x) \exp(\ii k_y  y)
\end{align}}
\end{definition}

To simplify notation, in the following $\Delta y = \Delta x$ is mostly assumed. We have the standard approximation result:

\begin{theorem} \label{thm:approximationfourier}
 Consider $q = \hat Q \exp(\ii \vec k \cdot \vec x) \in \mathbb M^m_{\vec k}$. If $v = \widehat{\dof}(\hat Q) + \mathcal O(\Delta x^{K+1})$, then $\| q -  R v \|_{L^2} \in \mathcal O(\Delta x^{K+1})$.
\end{theorem}
\begin{proof}
 The approximation $Rv$ matches the Taylor series of $q$ up to its degree $K$. Thus the remaining terms are a polynomial of at least degree $K+1$.
\end{proof}

Next, the object of study are linear subspaces $S_{\vec k}$ of $\mathbb C^{m(K+1)^2}$ and their approximation properties (they in general will depend on $\vec k$, hence the notation), in particular those arising as kernels of $\mathcal E \in \mathscr M^{m(K+1)^2 \times m(K+1)^2}$. Recall that the aim is to compare $\ker \vec J \cdot \vec k$ and $\ker \mathcal E$, but the former contains vectors in $\mathbb C^m$ and the latter in $\mathbb C^{m(K+1)^2}$. There are two natural ways of comparing them since there are two natural operators that map between the two: given an element $v \in S_{\vec k}$, either
\begin{itemize}
\item the reconstruction operator $Rv$ maps it to polynomials with coefficients in $\mathbb C^m$, which can be compared to Fourier modes $\mathbb M^m_{\vec k} \ni q$ as $\| q - R v \|_{L^2}$, or
\item it can be compared to the degrees of freedom $\widehat{\dof}(\hat Q) \in \mathbb C^{m(K+1)^2}$ of a Fourier mode $q = \hat Q \exp(\ii \vec k \cdot \vec x) \in \mathbb M^m_{\vec k}$ as $\| \widehat{\dof}(\hat Q) - v\|_{\ell^2}$. 
\end{itemize}
The latter version is considered easier and is the one used here. No notational distinction is made, but the latter norm is simply the Euclidean norm on $\mathbb C^{m(K+1)^2}$, while the former is an actual $L^2$ norm over the unit square.

The last conceptual difficulty is that we are not given elements to compare but subspaces. In particular, the kernel of $\mathcal E$ is given as a span of some rather arbitrarily chosen basis elements. Given a stationary ($\hat Q \in \ker \,\vec J \cdot \vec k$) Fourier mode $q := \hat Q \exp(\ii \vec k \cdot \vec x) \in \mathbb M^m_{\vec k}$, in general $\widehat{\dof}(\hat Q)$ will not be in $\ker \mathcal E$. One thus needs to find (in some sense) the closest or most relevant element in the kernel. The most natural way to proceed would be to decompose $\widehat{\dof}(\hat Q)$ in the eigenbasis of $\mathcal E$ and to keep only the elements in $\ker \mathcal E$. This would yield a map $\pi_{\vec k,\text{evo}} \circ \widehat{\dof}$ from $\mathbb C^m \ni \hat Q$ to $\ker \mathcal E \subset \mathbb C^{m(K+1)^2}$. Let us call $\pi_{\vec k,\text{evo}} \colon \mathbb C^{m(K+1)^2} \to \mathbb C^{m(K+1)^2}$ the \emph{evolutionary projector}, since upon time evolution of $\dof(q)$ the long-time limit of a von-Neumann-stable method will be just $\pi_{\vec k,\text{evo}} \widehat{\dof}(\hat Q)$. Unfortunately, it is very difficult to compute $\pi_{\vec k,\text{evo}}$ in practice. Instead, the following projectors will be used:

 \begin{definition}
 \begin{itemize}
  \item A projector $\pi_{\vec k} \colon \mathbb C^{m(K+1)^2} \to S_{\vec k} \subset \mathbb C^{m(K+1)^2}$ is called \emph{selectively exact}, if there exists an index set $I$ such that
 \begin{enumerate}
  \item $\# I = \dim S_{\vec k}$ and
  \item $(\pi_{\vec k} \hat q)_i =\hat q_i$ $\forall i \in I$, $\forall \hat q \in \mathbb C^{m(K+1)^2}$ \label{it:selectivelyexact}
 \end{enumerate}
 \item  It is called \emph{selectively optimal}, if instead of \ref{it:selectivelyexact}.
 \begin{enumerate}[1{a.}] \setcounter{enumi}{1}
  \item $(\pi_{\vec k} \hat q)_i = \hat q_i + \mathcal O(\Delta x^{K+1})$ $\forall i \in I$, $\forall q \in \mathbb C^{m(K+1)^2}$
 \end{enumerate}
 holds.
 \end{itemize}

\end{definition}

Such a projector selects certain degrees of freedom that are then approximated to the optimal order of accuracy, hence the name. The values of the other components then follow uniquely from the condition that the resulting vector needs to be in $S_{\vec k}$, and they might be anything between optimally accurate as well, or not even consistent with $\hat q$ (e.g. they could be identically zero). Clearly, a selectively optimal projector always exists: $S_{\vec k}$ has $\dim S_{\vec k}$ basis elements, and written as a matrix they can be brought to reduced row echelon form, which explicitly yields the set $I$.

It is not a priori clear that $\pi_{\vec k,\text{evo}}$ is selectively optimal for any choice of $I$, because the method might lose its order of accuracy at stationary state. However, if the method does retain its order of accuracy, $\pi_{\vec k,\text{evo}}$ is selectively optimal for any choice of $I$. To demonstrate the loss of accuracy for certain choices of numerical flux for DG, the strategy below therefore is to show that no selectively optimal projector $\pi_{\vec k}$ leads to an optimal order of accuracy in all components of $\pi_{\vec k}\widehat\dof(\hat Q)$. On the other hand, if such a selectively optimal projector can be found, it remains only a necessary condition for the claim that the stationary states are of optimal order of accuracy.

\begin{example} \label{example:upwindker}
The kernel of $\mathcal E$ for $K=1$ and upwind flux is one-dimensional and given by multiples of
 \begin{align}
  w := \left( - \frac{\sqrt{3}\Delta x(t_x+1)}{\Delta y(t_x-1)}, 0, -\frac{\Delta x}{\Delta y}, 0, \Big | \frac{\sqrt{3} (t_y+1)}{t_y-1}, 1, 0, 0,\Big | 0, 0, 0, 0 \right)^\text{T} \label{eq:kerupwind}
 \end{align}
 where the variables are ordered as 
 \begin{align*}
 \hat u^{(00)},\hat u^{(01)},\hat u^{(10)},\hat u^{(11)},\hat v^{(00)},\hat v^{(01)},\hat v^{(10)},\hat v^{(11)}, \hat p^{(00)},\hat p^{(01)},\hat p^{(10)},\hat p^{(11)}
 \end{align*}
 Let us construct some selectively optimal projectors; choose $I$ to contain the index associated to $v^{(00)}$. Then a selectively exact projector onto $\mathrm{span}(w)$ is
 \begin{align}
  \pi_{\vec k} \hat q := q_{v^{(00)}} \frac{t_y-1}{\sqrt{3} (t_y+1)} w
 \end{align}
where $q_{v^{(00)}}$ is the 5-th component of $q$.

Consider now $\hat Q = (-k_y, k_x, 0)^\text{T}$, such that $\hat Q \in \ker \, \vec J \cdot \vec k$. One easily computes 
\begin{align}
 \widehat{\dof}(\hat Q)_{00} = \hat Q \frac{4 \sin \left( \frac{\Delta x k_x}{2} \right ) \sin \left( \frac{\Delta y k_y}{2} \right )}{\Delta x \Delta y k_x k_y} =: \hat Q \alpha
\end{align}
Thus, {\footnotesize
\begin{align}
\pi_{\vec k} \widehat{\dof}(\hat Q) &= k_x \alpha \frac{t_y-1}{\sqrt{3} (t_y+1)} \left( - \frac{\sqrt{3}\Delta x(t_x+1)}{\Delta y(t_x-1)}, 0, -\frac{\Delta x}{\Delta y}, 0, \Big | \frac{\sqrt{3} (t_y+1)}{t_y-1}, 1, 0, 0,\Big | 0, 0, 0, 0 \right)^\text{T}\\
&= \left( -k_y \alpha + \mathcal O(\Delta x^2), \dots \Big | k_x \alpha, \dots,\Big | 0, 0, 0, 0 \right)^\text{T}
\end{align}}
which shows that with this projector, the $u^{(00)}$ would end up having the correct order of accuracy. One might now proceed to study the other terms. For the sake of easier computation, though, consider a different selectively optimal projector defined as
\begin{align}
  \pi_{\vec k} \hat q := \left(q_{v^{(01)}} + a_0 \Delta x^2\right ) w
 \end{align}
One then finds, for any $a_0$ and with $\Delta y = \Delta x$ 
\begin{align}
 \widehat{\dof}(\hat Q) - \pi_{\vec k}\widehat{\dof}(\hat Q)  = \left( 0, -\frac{\ii \Delta x k_y^2}{2\sqrt{3}}, 0, 0   , \Big |  0, 0, \frac{\ii \Delta x k_x^2}{2\sqrt{3}} , 0   , \Big | 0, 0, 0, 0     \right )^\text{T} + \mathcal O(\Delta x^2)
\end{align}
One thus has to conclude that the upwind flux for $K=1$ is at most first-order accurate, instead of second. As can be seen in Figure \ref{fig:upwindrusanov-order} below, this is what is observed experimentally.
\end{example}

It was mentioned above that the dimension of the kernel of $\mathcal E$ might be taken as indication of its approximation properties: since the kernel of $\vec J \cdot \vec k$ is one-dimensional, and since every scalar is represented by $(K+1)^2$ degrees of freedom per cell, one might expect that a $(K+1)^2$-dimensional kernel of $\mathcal E$ (as it happens for central flux) will correspond to a $(K+1)$-accurate stationary state, while a $K^2$-dimensional kernel (as that of upwind flux) will display one order less at stationary state. However, this is generally false: it is neither necessary to have a $(K+1)^2$-dimensional kernel for optimal accuracy (see Example \ref{example:singlekernel}), nor is it sufficient (Example \ref{example:centralker}). The approximation properties depend on the nature of the elements, not just on their number.

\begin{example} \label{example:singlekernel}
 Consider a subspace spanned by
 \begin{align}
  &\left(1, \frac{t_y-1}{\sqrt{3}(t_y+1)}, \frac{t_x-1}{\sqrt{3}(t_x+1)}, \frac{(t_x-1)(t_y-1)}{3(t_x+1)(t_y+1)}, \right . \label{eq:hypoker1}\\&\phantom{mmmmm} \left. \nonumber - \frac{k_x}{k_y}, -\frac{k_x(t_y-1)}{k_y\sqrt{3}(t_y+1)}, -\frac{k_x(t_x-1)}{k_y\sqrt{3}(t_x+1)}, -\frac{k_x(t_x-1)(t_y-1)}{3k_y(t_x+1)(t_y+1)}    , \Big | 0, 0, 0, 0 \right)^\text{T}
 \end{align}
 (with the same ordering of variables as in Example \ref{example:upwindker}). Define $\pi_{\vec k}$ such that $u^{(00)}$ is exact. Even though this is a one-dimensional subspace, one finds
\begin{align}
 \widehat{\dof}(\hat Q) - \pi_{\vec k}\widehat{\dof}(\hat Q) \in \mathcal O(\Delta x^2)
\end{align}
 A hypothetical numerical method with $\ker \mathcal E$ spanned by \eqref{eq:hypoker1} would be high-order on the stationary state, despite the one-dimensionality of the kernel.
\end{example}

\begin{example} \label{example:centralker}
 The kernel of $\mathcal E$ for central DG with $K=1$ for linear acoustics is spanned by {\tiny
 \begin{align*}
  w_1 = \left(-\frac{ \Delta x  (t_y -1)^2}{4  \Delta y   t_y },-\frac{\sqrt{3}  \Delta x  \left(t_y^2-1\right)}{4  \Delta y 
 t_y },-\frac{ \Delta x  (1+ t_x ) (t_y -1)^2}{4 \sqrt{3}  \Delta y  (t_x -1)  t_y },-\frac{ \Delta x  (1+ t_x ) \left(-1+ t_y^2\right)}{4
 \Delta y  (t_x -1)  t_y }, \Big | 0,0,0,1,\Big | 0,0,0,0\right)^\text{T}\\
w_2 = \left(\frac{ \Delta x  \left(t_y^2-1\right)}{4 \sqrt{3}  \Delta y   t_y },\frac{ \Delta x 
(t_y -1)^2}{4  \Delta y   t_y },\frac{ \Delta x  (1+ t_x ) \left(t_y^2-1\right)}{12  \Delta y  (t_x -1)  t_y },\frac{ \Delta x 
(1+ t_x ) (t_y -1)^2}{4 \sqrt{3}  \Delta y  (t_x -1)  t_y },\Big | 0,0,1,0,\Big | 0,0,0,0\right)^\text{T}\\
w_3 = \left(\frac{\sqrt{3}  \Delta x  (1+ t_x )
(t_y -1)^2}{4  \Delta y  (t_x -1)  t_y },\frac{3  \Delta x  (1+ t_x ) \left(t_y^2-1\right)}{4  \Delta y  (t_x -1)  t_y },\frac{ \Delta x 
( t_y -1)^2}{4  \Delta y   t_y },\frac{\sqrt{3}  \Delta x  \left(t_y^2-1\right)}{4  \Delta y   t_y },\Big | 0,1,0,0,\Big | 0,0,0,0\right)^\text{T}\\
w_4 = \left(-\frac{ \Delta x 
(1+ t_x ) \left(t_y^2-1\right)}{4  \Delta y  (t_x -1)  t_y },-\frac{\sqrt{3}  \Delta x  (1+ t_x ) (t_y -1)^2}{4  \Delta y 
(t_x -1)  t_y },-\frac{ \Delta x  \left(t_y^2-1\right)}{4 \sqrt{3}  \Delta y   t_y },-\frac{ \Delta x  (t_y -1)^2}{4  \Delta y 
 t_y },\Big | 1,0,0,0,\Big | 0,0,0,0\right)^\text{T}
 \end{align*}}
 The ordering of variables is the same as in Example \ref{example:upwindker}.
 These are $(K+1)^2 = 4$ vectors. However, any selectively optimal projector $\pi_{\vec k}$ based on the entries $4-8$ (related to degrees of freedom of $v$) yields
 \begin{align}
 \widehat{\dof}(\hat Q) - \pi_{\vec k}\widehat{\dof}(\hat Q)  = \left( 0, -\frac{\ii \Delta x k_y^2}{2\sqrt{3}}, -\frac{\ii \Delta x k_x k_y}{6\sqrt{3}}, 0   , \Big |  0, 0,0 , 0   , \Big | 0, 0, 0, 0     \right )^\text{T} + \mathcal O(\Delta x^2)
\end{align}
Surprisingly, therefore, central DG for $K=1$ is also only first, and not second order accurate. This is is surprising not only because there are four elements in the kernel of $\mathcal E$, but also because usually a loss of accuracy is associated with the numerical stabilization, absent in this case. This behaviour is confirmed experimentally and further comments are provided in Section \ref{sec:centralflux}.
\end{example}

Generally speaking, DG achieves its high order of accuracy by incorporating additional degrees of freedom. However, the kernel of $\mathcal E$ is spanned by vectors that contain polynomials in $t_x$ and $t_y$, i.e. their entries are finite difference formulas with a possibly large stencil. A derivative can be approximated locally by using the fact that the solution inside the cell is a polynomial and by differentiating it, or by using the cell averages in some stencil and applying a finite-difference-type formula. The kernel of $\mathcal E$ and thus the stationary states of DG make use of a mix of these two ways of achieving high order, and this is why deviations from the optimal order of accuracy can happen, and actually do happen, as shown above.

Next, a systematic presentation of the accuracy of DG at stationary state is shown for 4 choices of numerical flux.

\subsection{Order of accuracy at steady state}

As has been discussed in the Introduction, Finite Volume methods for linear acoustics (and hence DG with $K=0$) lose consistency at stationary state if an upwind or Rusanov flux is used, while with the central flux and fluxes originating from ``low Mach fixes'', the discrete stationary states remain consistent discretizations of the steady states of the PDE. The aim of this Section is to demonstrate how the picture changes for higher degrees of the polynomial reconstruction.

The basis elements are stated according to a lexicographic ordering of the degrees of freedom for each variable; the one for $K=1$ is stated in Example \ref{example:upwindker}; the one for $K=2$ is {\footnotesize
\begin{align*}
 u^{(00)},u^{(01)},u^{(02)},u^{(10)},u^{(11)},u^{(12)},u^{(20)},u^{(21)},u^{(22)},\Big |v^{(00)},v^{(01)},v^{(02)},v^{(10)},\dots,\Big |p^{(00)},\ldots,p^{(22)}
\end{align*}}
and the one for $K=3$ is {\footnotesize
\begin{align*}
 u^{(00)},u^{(01)},u^{(02)},u^{(03)},u^{(10)},u^{(11)},u^{(12)},\dots,\Big |v^{(00)},v^{(01)},v^{(02)},v^{(03)},v^{(10)},\dots,\Big |p^{(00)},\dots,p^{(33)}
\end{align*}}

Vertical lines help to separate the blocks associated to different variables $u$, $v$, $p$.
The notation $\boxed{N \times 0}$ abbreviates $N$ vanishing entries of the vectors in components associated to the pressure. The kernels of the matrices have been determined using \textsc{mathematica} due to the great length of the expressions.

\begin{figure}
 \centering
  \includegraphics[width=0.45\textwidth]{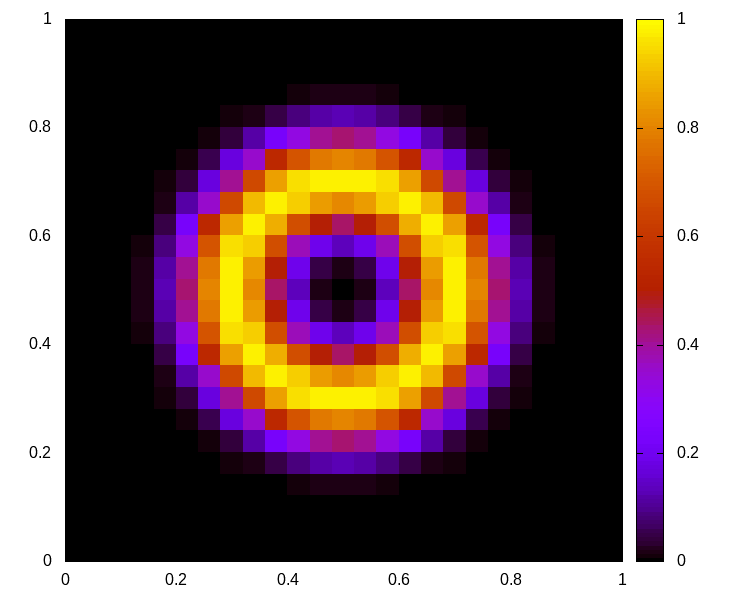} \includegraphics[width=0.45\textwidth]{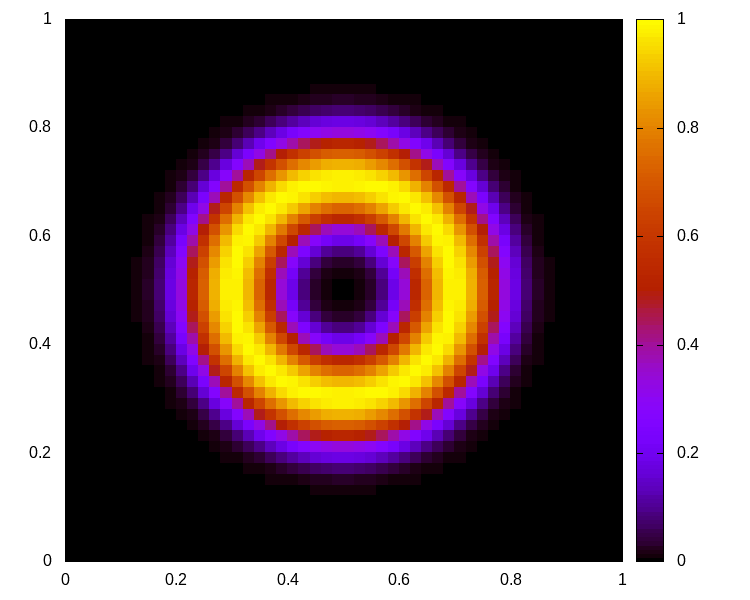} \\
  \caption{Initial data and exact solution used for the convergence studies for linear acoustics. The value of $|\vec v|$ at cell center is shown. \emph{Left}: Grid of $25 \times 25$. \emph{Right}: Grid of $50 \times 50$.}
 \label{fig:acoustcvortex}
\end{figure}

The experimental measurement of the order of accuracy is performed using the following (stationary) initial data $q_0 = (u_0, v_0, p_0 = 0)^\text{T}$ for linear acoustics:
\begin{align}
r &:= \sqrt{\left(x-\frac12\right)^2 + \left(y - \frac12\right)^2} &
V(r) &:= \begin{cases} 0 & r > 2w \\ \frac14 (1 - \cos(\pi r/w))^2 & \text{else} \end{cases} \label{eq:acousticvortex1}\\
u_0 &= -V(r) \frac{y-\frac12}{r} &
 v_0 &= V(r) \frac{x-\frac12}{r} \label{eq:acousticvortex2}
\end{align}
A value $w = 0.2$ is used. The data are shown in Figure \ref{fig:acoustcvortex}.
The setup is considered on a periodic domain $[0,1]^2$. The integration in time uses classical Runge-Kutta methods of the same order as the spatial discretization, unless stated otherwise. Note that the accuracy of the time integrator is irrelevant for the accuracy of the discrete stationary state. A CFL number of $0.03$ is used for all setups.

The order of accuracy at any time $t$ is obtained as $\log_2 \frac{\|q_0 - q_h^{25 \times 25}(t) \|_{L^2}}{\| q_0 - q_h^{50 \times 50}(t)\|_{L^2}}$. 
Occasionally, to obtain a more precise estimate of the order, the errors on grids of $50 \times 50$ and $100 \times 100$ cells are used instead; this is indicated in the Figures below as ``fine grids''. All integrals are systematically replaced by 5-point Gauss quadratures.

\subsubsection{Upwind flux}

\begin{figure}
 \centering
 \includegraphics[width=0.49\textwidth]{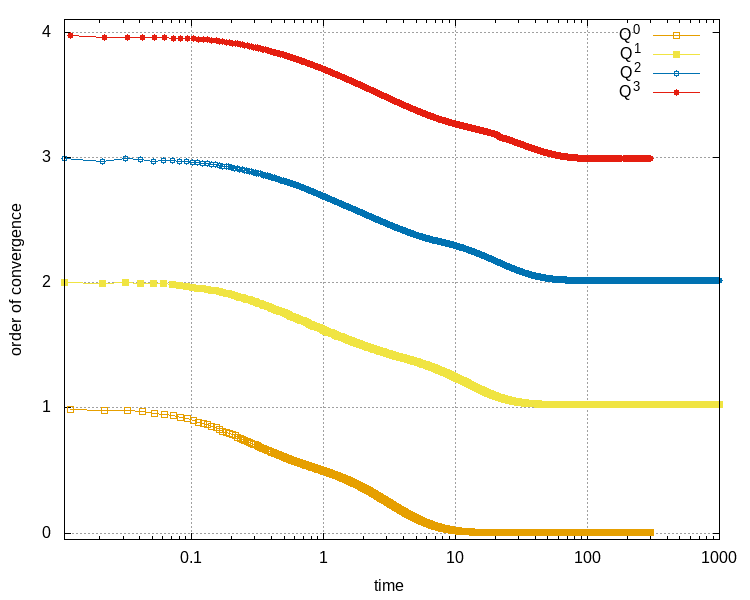} \hfill \includegraphics[width=0.49\textwidth]{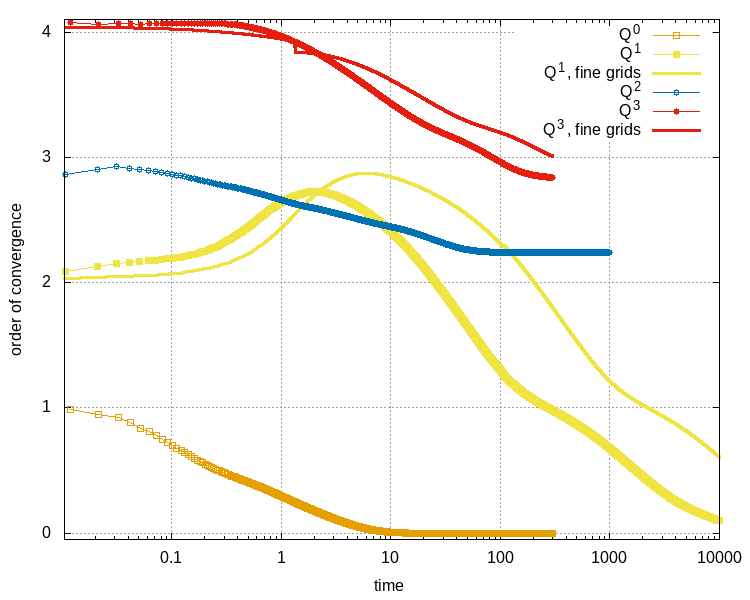}
 \caption{\emph{Left}: DG with an upwind flux: Order of convergence as a function of time, obtained from simulations on grids with $25 \times 25$ and $50\times 50$ cells. One observes one order less than expected. \emph{Right}: DG with a Rusanov flux: Order of convergence as a function of time, obtained from simulations on grids with $25 \times 25$ and $50\times 50$ cells (and using also $100 \times 100$ for certain cases). One observes that with a $Q^0$ or $Q^1$ approximation the stationary states are not consistently discretized; second order instead of third is observed for $Q^2$, third instead of fourth for $Q^3$. The long-time simulations are computationally very intensive and have not been performed for all setups.}
 \label{fig:upwindrusanov-order}
\end{figure}

Figure \ref{fig:upwindrusanov-order} (left) shows that when using an upwind numerical flux, i.e. diffusion matrices
\begin{align}
 D_x &= \left( \begin{array}{ccc}  1 & 0 & 0 \\ 0 & 0 & 0\\ 0 & 0 &1  \end{array} \right )
 & D_y &= \left( \begin{array}{ccc}  0 & 0 & 0 \\ 0 & 1 &0 \\0 & 0 &1  \end{array} \right )
\end{align}
in \eqref{eq:numfluxx}--\eqref{eq:numfluxy},
one experimentally observes a drop by one order when the scheme reaches its stationary state. For $Q^1$ approximations, for example, the discrete steady state is first- instead of second-order accurate. Even if sub-optimal, this, clearly, is better than losing consistency.

The analysis of the kernel of the evolution matrix for the cases $Q^K$, $K=1,2,3$ shows that its dimension is only $K^2$, while there are $(K+1)^2$ degrees of freedom per scalar. As has been discussed above, it is easy to conjecture that DG with the upwind flux will suffer a loss of one order. Even if a correct conclusion in this case, the argument is flawed, since the accuracy depends on the elements of the kernel, not just on its dimension. 

The accuracy can have two sources: One can be accurate because there are enough degrees of freedom per cell, or one can be accurate because of a large enough stencil. This latter aspect cannot be assessed by merely counting the basis elements of the kernel, it depends on what these basis elements are exactly. One can imagine hypothetical situations of basis elements spanning only the subspace of the $u$ variable, which would mean that they would not correspond to a consistent discretization of the steady states of the PDE, even if there are $(K+1)^2$ of them. Conversely, it is possible for just a one-dimensional kernel to correspond to a high-order accurate discretization of steady states, as has been shown in Example \ref{example:singlekernel}. In the case of the upwind flux, however, following the procedure outlined above one does confirm theoretically (for $K=1,2,3$), that the order of accuracy associated to the kernel of $\mathcal E$ is at most $K$ instead of $K+1$. 

In particular one finds the following. For $K=1$ the kernel is spanned by
\begin{align}
\left (-\frac{\sqrt{3}  \Delta x (1+ t_x)}{ \Delta y (t_x-1)},0,-\frac{ \Delta x}{ \Delta y},0,\Big |\frac{\sqrt{3} (1+ t_y)}{t_y-1},1,0,0,\boxed{4 \times 0}\right )^\text{T}
\end{align}
Its approximation properties were subject of Example \ref{example:upwindker}.
For $K=2$ the kernel of $\mathcal E$ is spanned by the four vectors {\footnotesize
\begin{align*}
w_1 &:= \left (0,\frac{\sqrt{5}  \Delta x}{ \Delta y},0,0,0,0,0,-\frac{ \Delta x}{ \Delta y},0,\Big |   0,0,0,-\sqrt{5},0,1,0,0,0,\boxed{9 \times 0}\right )^\text{T}\\ 
w_2 &:= \left (\frac{ \Delta x}{ \Delta y},0,0,0,0,0,-\frac{ \Delta x}{\sqrt{5}  \Delta y},0,0,\Big |0,0,0,\frac{\sqrt{3} (1+ t_y)}{t_y-1},1,0,0,0,0,\boxed{9 \times 0}\right )^\text{T} \\ 
w_3 &:= \left (0,-\frac{\sqrt{15}  \Delta x (1+ t_x)}{ \Delta y (t_x-1)},0,0,-\frac{\sqrt{5}  \Delta x}{ \Delta y},0,0,0,0,\Big |-\sqrt{5},0,1,0,0,0,0,0,0,\boxed{9 \times 0}\right )^\text{T} \\ 
w_4 &:= \left (-\frac{\sqrt{3}  \Delta x (1+ t_x)}{ \Delta y (t_x-1)},0,0,-\frac{ \Delta x}{ \Delta y},0,0,0,0,0,\Big |\frac{\sqrt{3} (1+ t_y)}{t_y-1},1,0,0,0,0,0,0,0,\boxed{9 \times 0}\right )^\text{T} 
\end{align*}}

A selectively optimal projector is given by
\begin{align*}
  \pi_{\vec k} \hat q &:= \left(q_{v^{(12)}} + a_0 \Delta x^3\right ) w_1 + \left(q_{v^{(11)}} + a_1 \Delta x^3\right ) w_2 \\&\qquad+ \left(q_{v^{(02)}} + a_2 \Delta x^3\right ) w_3 + \left(q_{v^{(01)}} + a_3 \Delta x^3\right ) w_4
\end{align*}
However, independently of the choice of $a_0, \dots, a_3$ one finds {\footnotesize
\begin{align}
 \widehat{\dof}(\hat Q) - \pi_{\vec k}\widehat{\dof}(\hat Q)  &= \\\nonumber & \!\!\!\!\!\!\!\!\!\!\!\!\!\!\!\!\left( *,*, \frac{\Delta x^2 k_y^3}{12\sqrt{5}},0,0,0,0,0,0,\Big| *,0,0,*,0,0, -\frac{ \Delta x^2 k_x^3}{12\sqrt{5}}, 0 ,0,  \boxed{9 \times 0} \right )^\text{T} + \mathcal O(\Delta x^3)
\end{align}}
where $*$ denotes some non-zero entry no worse than $\mathcal O(\Delta x^2)$ and thus of no further relevance. Thus, DG with upwind flux for $K=2$ is at most second-order accurate.

For $K=3$, the kernel is spanned by the elements given in Appendix \ref{app:ssec:upwindK3}; the corresponding calculation involving a selectively optimal projector is omitted here. The result is that it is unavoidable to have error terms $\mathcal O(\Delta x^3)$, while optimal order of accuracy would entail $\mathcal O(\Delta x^4)$. The drop of one order is, of course, most noticeable when the optimal order of accuracy is 1 itself ($K=0$), i.e. when upon stationarity the scheme loses consistency.

\begin{figure}
 \centering
  \includegraphics[width=0.45\textwidth]{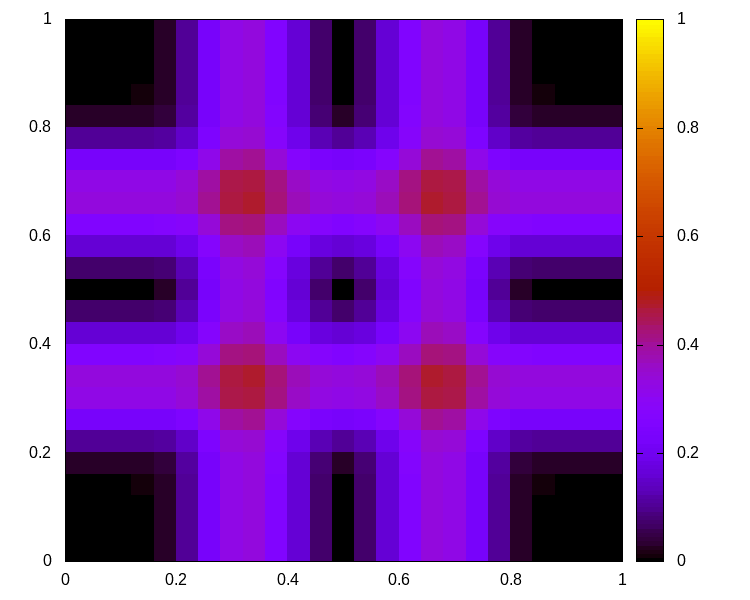} \includegraphics[width=0.45\textwidth]{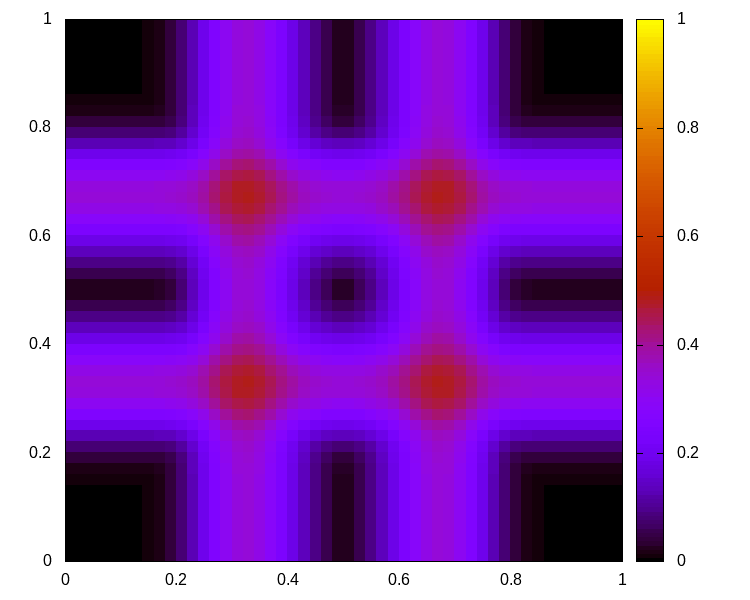} \\
 \includegraphics[width=0.45\textwidth]{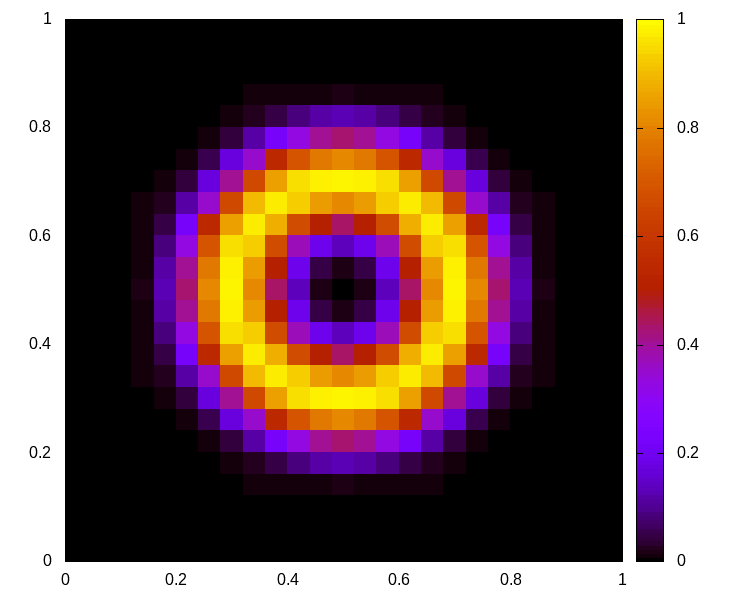} \includegraphics[width=0.45\textwidth]{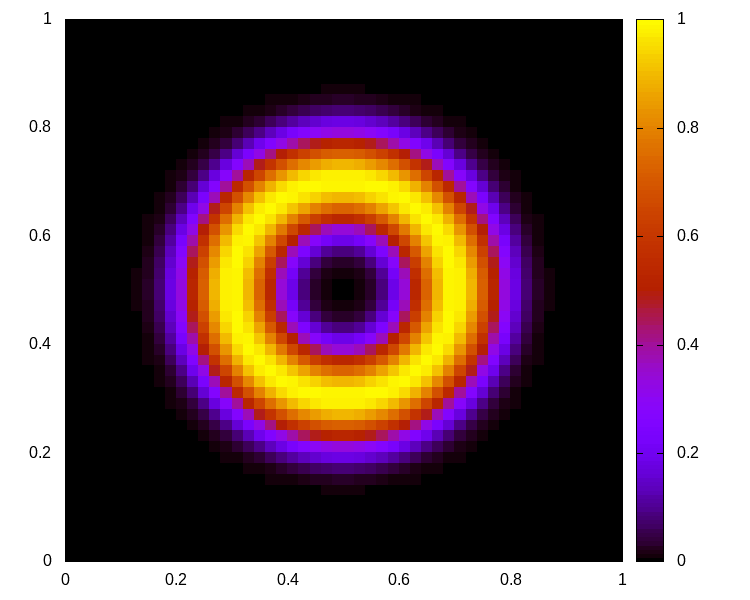} \\
 \includegraphics[width=0.45\textwidth]{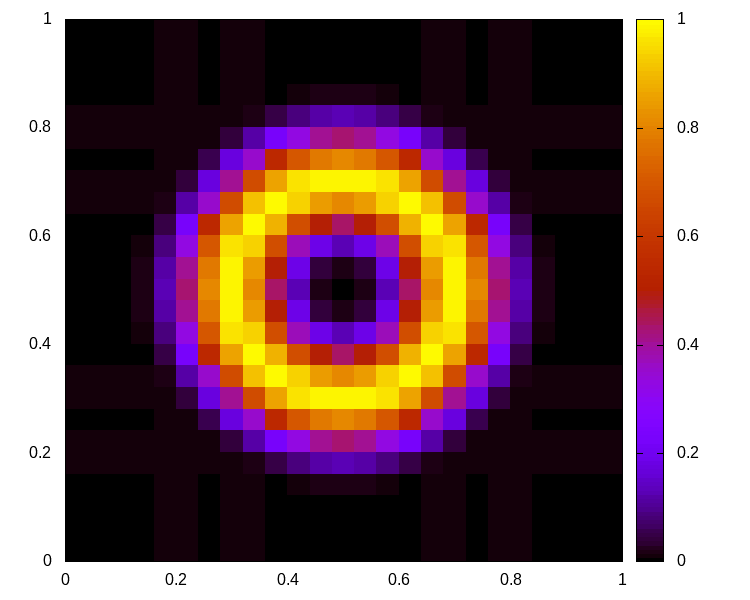} \includegraphics[width=0.45\textwidth]{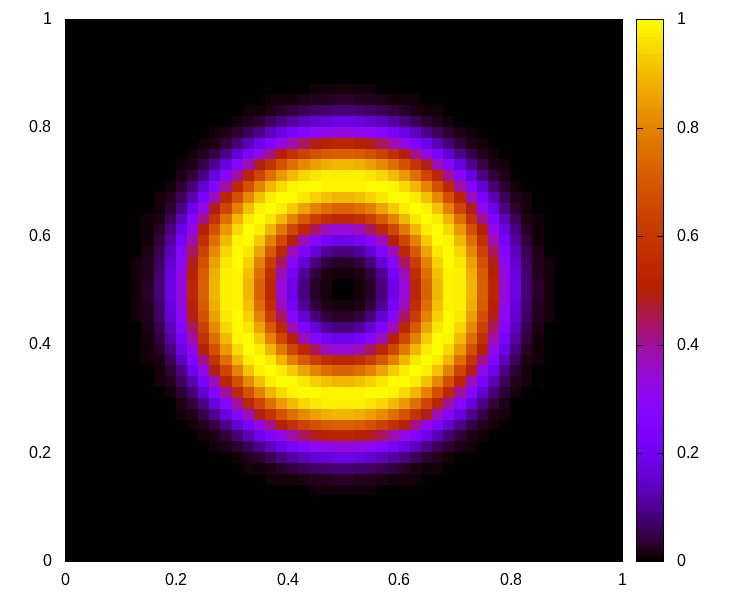} \\
  \caption{Upwind DG for linear acoustics; setup \eqref{eq:acousticvortex1}--\eqref{eq:acousticvortex2}. \emph{Top}: Numerical solution at time $t=300$ with $K=0$. \emph{Middle}: Numerical solution at time $t=10^3$ with $K=1$. \emph{Bottom}: Same for $K=2$. The solution $|\vec v|$ at cell center is shown. \emph{Left}: Grid of $25 \times 25$. \emph{Right}: Grid of $50 \times 50$.}
 \label{fig:upwind}
\end{figure}

\begin{figure}
 \centering
 \includegraphics[width=0.32\textwidth]{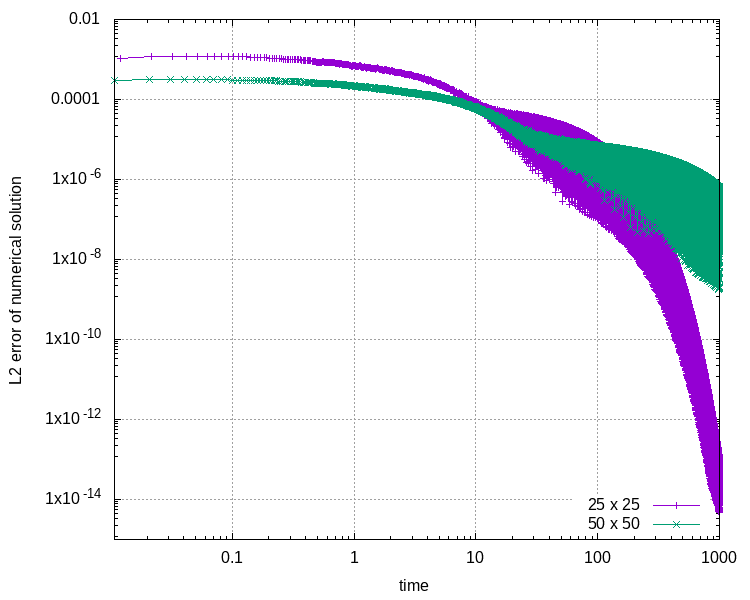} \hfill \includegraphics[width=0.32\textwidth]{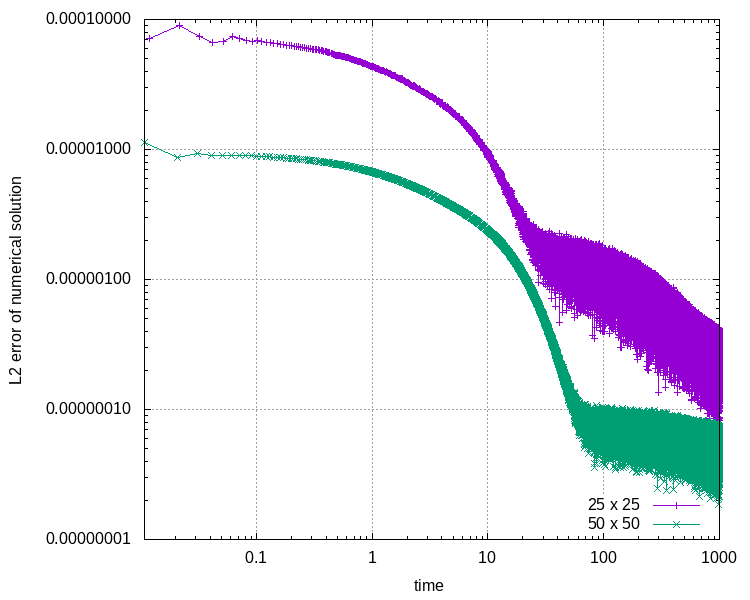} \hfill    \includegraphics[width=0.32\textwidth]{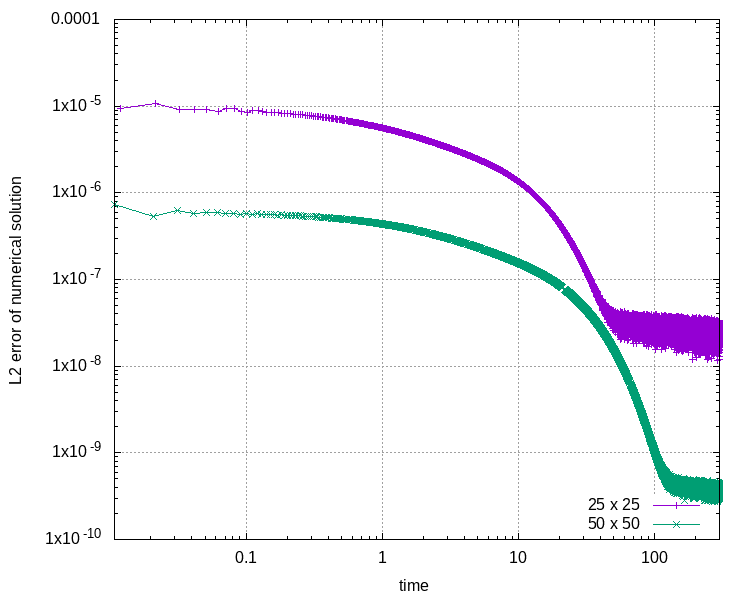}
 \caption{Error of the pressure for upwind DG as a function of time. \emph{Left}: $K=1$. \emph{Center}: $K=2$. \emph{Right}: $K=3$.}
 \label{fig:upwindconvergencepressure}
\end{figure}

Figure \ref{fig:upwind} shows the numerical solutions for various values of $K$; this Figure should be compared to \ref{fig:eulerM3roe} which shows the solutions of the Euler equations for low Mach number, computed with the Roe flux. Figure \ref{fig:upwindconvergencepressure} shows the error of the pressure as a function of time. Both the stationary states of the PDE and those of the numerical method are characterized by a constant pressure; by conservation one can expect the constant to be the same and thus the error should decay to machine precision. The plots thus allow to judge how close the simulation has got to the discrete stationary state.

\subsubsection{Rusanov flux}

The Rusanov flux for linear acoustics corresponds to the choice
\begin{align}
 D_x = D_y = \id_3
\end{align}
of the stabilization matrices in \eqref{eq:numfluxx}--\eqref{eq:numfluxy}.
Its sub-optimality is even more prominent, with the experimental orders of accuracy shown in Figure \ref{fig:upwindrusanov-order} (right). One observes that for a $Q^1$ approximation, the discrete stationary states are not discretizing those of the PDE consistently, the method is not stationarity preserving. For a $Q^2$ approximation one experimentally finds the discrete steady state to be second-order accurate instead of third, for $Q^3$ --- third-order accurate instead of fourth.

This matches the theoretical results: in fact, $\mathcal E$ does not have a non-trivial kernel neither for $Q^0$ nor for $Q^1$ approximations. The discrete steady states are essentially only those where $\del_x u = 0$ and $\del_y v = 0$ individually, or even just uniform constants. For a $Q^2$ approximation, the kernel of $\mathcal E$ is one-dimensional (i.e. of size $(K-1)^2$, compare this to $(K+1)^2$ degrees of freedom per scalar):

\begin{align}
w := \left (-\frac{2 \sqrt{15}  \Delta x (1+ t_x (4+ t_x)) (1+ t_y)}{ \Delta y (t_x-1)^2 (t_y-1)},-\frac{2 \sqrt{5}  \Delta x (1+ t_x (4+ t_x))}{ \Delta y (t_x-1)^2},0, \right . \\ \nonumber \left. -\frac{3 \sqrt{5}  \Delta x (1+ t_x) (1+ t_y)}{ \Delta y (t_x-1) (t_y-1)},-\frac{\sqrt{15}  \Delta x (1+ t_x)}{ \Delta y (t_x-1)},0,\right . \\ \nonumber \left.-\frac{\sqrt{3}  \Delta x (1+ t_y)}{ \Delta y (t_y-1)},-\frac{ \Delta x}{ \Delta y},0,\right . \\ \nonumber \left.\frac{2 \sqrt{15} (1+ t_x) (1+ t_y (4+ t_y))}{(t_x-1) (t_y-1)^2},\frac{3 \sqrt{5} (1+ t_x) (1+ t_y)}{(t_x-1) (t_y-1)},\frac{\sqrt{3} (1+ t_x)}{t_x-1},\right . \\ \nonumber \left.\frac{2 \sqrt{5} (1+ t_y (4+ t_y))}{(t_y-1)^2},\frac{\sqrt{15} (1+ t_y)}{t_y-1},1,0,0,0,\boxed{9 \times 0}\right )^\text{T}
\end{align}

Through merely a counting argument, one might expect only a first-order accurate approximation of steady states as a consequence -- recall that the kernel of $\mathcal E$ for the upwind flux in the $Q^1$ case is one-dimensional, and leads to a first-order approximation. However, as was emphasized earlier, the counting argument is not correct, and the actual values in the basis element of the kernel are important. Here, they are such that it is possible to obtain second order of accuracy, in line with the experimental finding. In particular, the following selectively optimal projector
\begin{align}
  \pi_{\vec k} \hat q := \left(q_{v^{(12)}} + a_0 \Delta x^3\right ) w 
\end{align}
in general gives $\widehat{\dof}(\hat Q) - \pi_{\vec k}\widehat{\dof}(\hat Q) \in \mathcal O(1)$, unless $a_0=0$, in which case{\footnotesize
\begin{align}
 \widehat{\dof}(\hat Q) - \pi_{\vec k}\widehat{\dof}(\hat Q) & = \\&\nonumber \!\!\!\!\!\!\!\!\!\!\!\!\! \left( *,0, \frac{\Delta x^2 k_y^3}{12\sqrt{5}},0,0,0,0,0,0,\Big| *,0,0,0,0,0, -\frac{ \Delta x^2 k_x^3}{12\sqrt{5}}, 0 ,0,  \boxed{9 \times 0} \right )^\text{T} + \mathcal O(\Delta x^3)
\end{align}}
Here, $*$ denotes other terms $\mathcal O(\Delta x^2)$ that are omitted for brevity. Thus, one can theoretically conclude that it is possible for DG with the Rusanov flux and $K=2$ to obtain second order of accuracy at stationary state. This seems to be the case experimentally.

Unfortunately, due to the length of the calculations it was not possible to obtain the kernel of $\mathcal E$ for the Rusanov flux in the case $K=3$. Experimentally, one observes that the accuracy at stationary state is order 3 instead of order 4.

\begin{figure}
 \centering
   \includegraphics[width=0.45\textwidth]{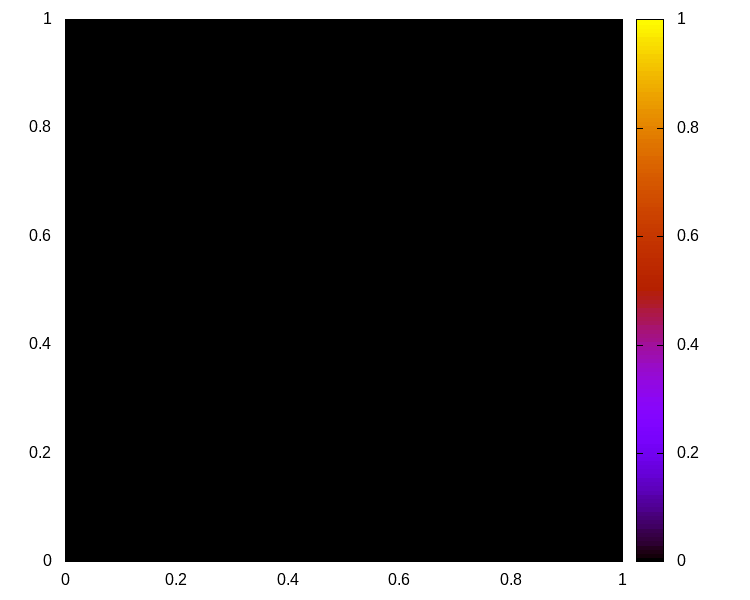} \includegraphics[width=0.45\textwidth]{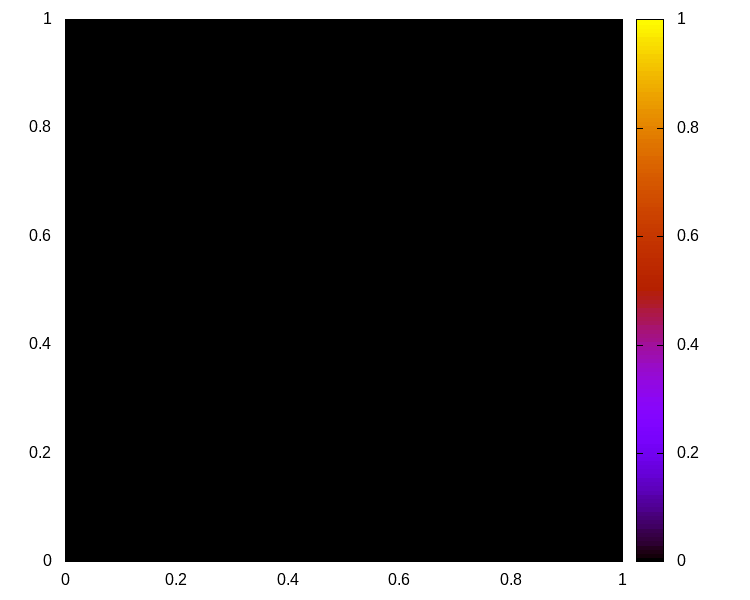} \\
 \includegraphics[width=0.45\textwidth]{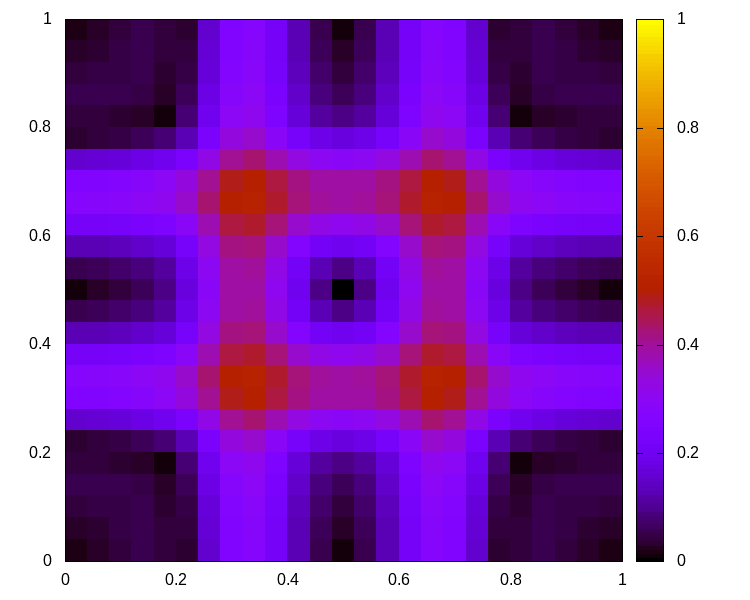} \includegraphics[width=0.45\textwidth]{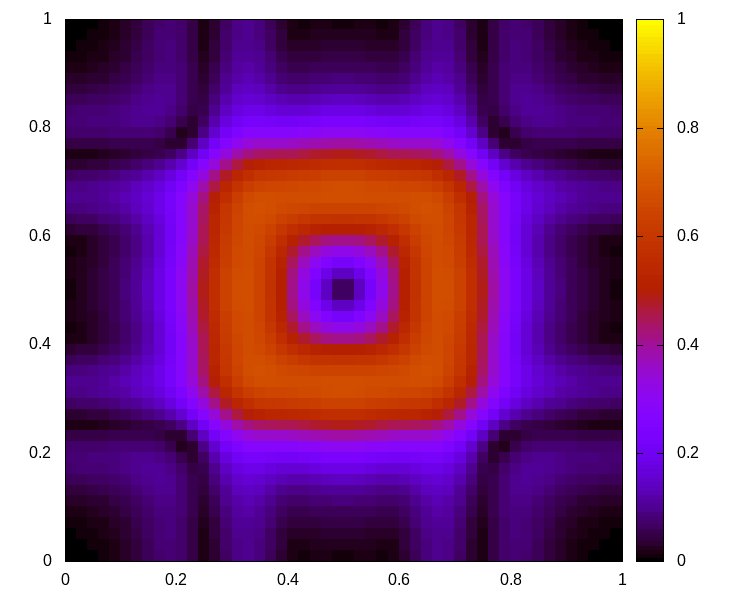} \\
 \includegraphics[width=0.45\textwidth]{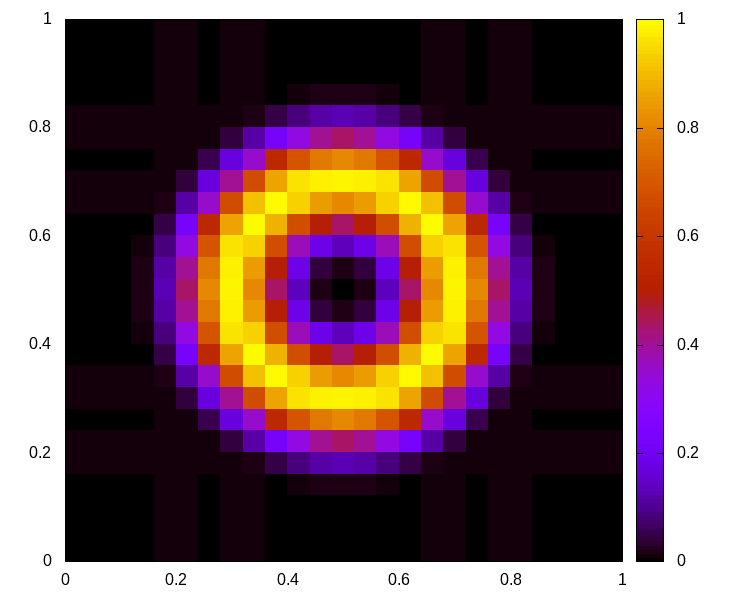} \includegraphics[width=0.45\textwidth]{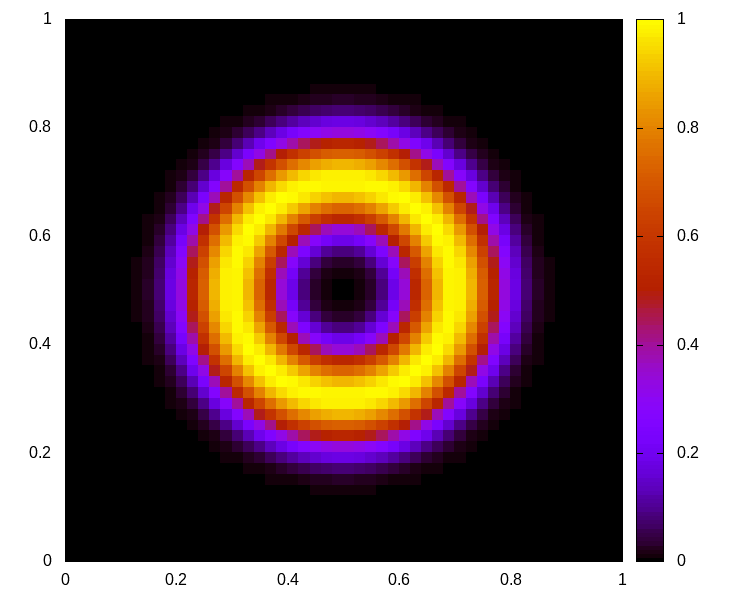} \\
\caption{Rusanov DG for linear acoustics; setup \eqref{eq:acousticvortex1}--\eqref{eq:acousticvortex2}. \emph{Top}: Numerical solution at time $t=300$ with $K=0$. It is diffused down to machine precision. \emph{Center}: Numerical solution at time $t=10^3$ with $K=1$. \emph{Bottom}: Same with $K=2$. The solution at cell center is shown. \emph{Left}: Grid of $25 \times 25$. \emph{Right}: Grid of $50 \times 50$.}
 \label{fig:rusanov}
\end{figure}

\begin{figure}
 \centering
\includegraphics[width=0.32\textwidth]{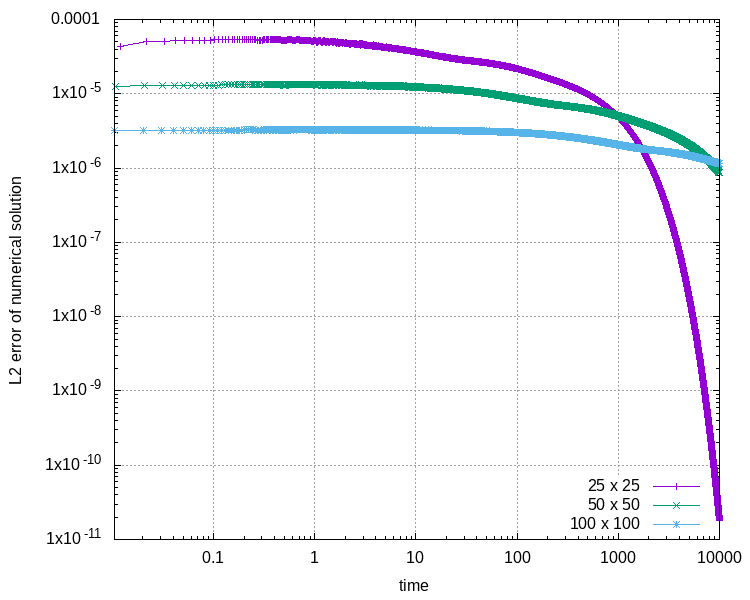}
 \hfill \includegraphics[width=0.32\textwidth]{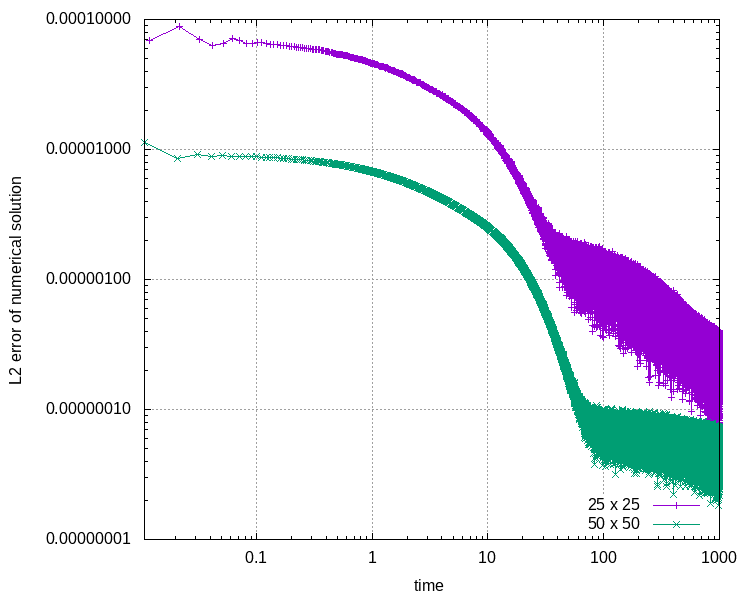}
 \hfill \includegraphics[width=0.32\textwidth]{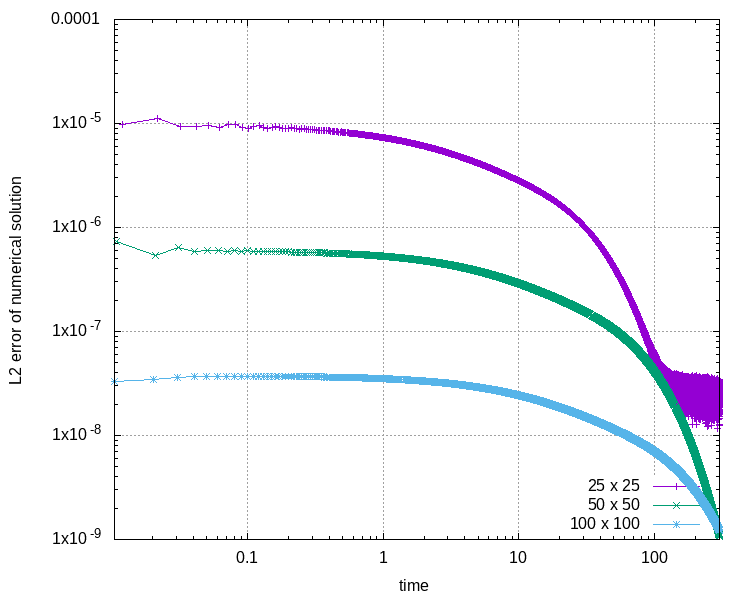}
 \caption{Error of the pressure for Rusanov DG as a function of time. \emph{Left}: $K=1$. \emph{Center}: $K=2$ \emph{Right}: $K=3$.}
 \label{fig:rusanovconvergencepressure}
\end{figure}

Figure \ref{fig:rusanov} shows the numerical solution for various values of $K$; it should be compared to Fig. \ref{fig:eulerM3rusanov} showing the low Mach number results for the Euler equations with the Rusanov flux. Figure \ref{fig:rusanovconvergencepressure} shows how the error of the pressure decays towards machine zero.

\subsubsection{Central flux} \label{sec:centralflux}

Central flux makes two changes necessary: first, the method is integrated with RK3 even for lower spatial accuracy in order to ensure stability. Second, in order to accelerate the transition to the discrete steady state, diffusion in the pressure is added, i.e. instead of $D_x = D_y = 0_{3 \times 3}$ the following matrices are used:
\begin{align}
 D_x &= D_y = \left( \begin{array}{ccc}  0 & 0 & 0 \\ 0 &0 &0 \\0 & 0 &1  \end{array} \right )
\end{align}
This change does not modify the discrete steady state, because the one of pure central flux is already characterized by a constant pressure. In \cite{gjonaj07} it has already been shown that the kernel of central DG is $\ndof$-dimensional.

\begin{figure}
 \centering
 \includegraphics[width=0.49\textwidth]{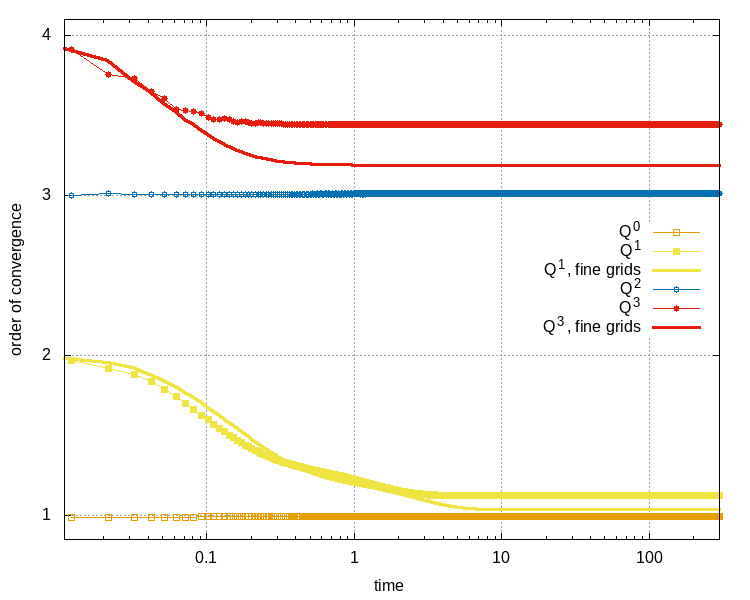} \hfill\includegraphics[width=0.5\textwidth]{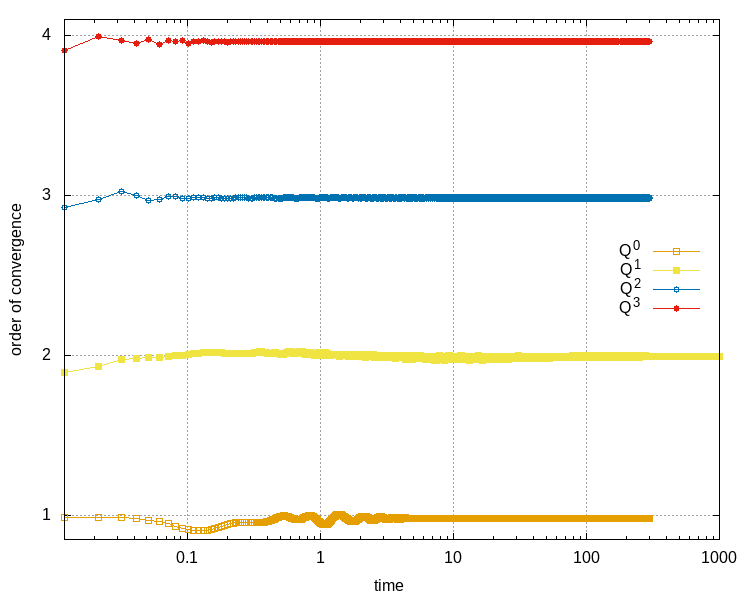}
 \caption{\emph{Left}: DG with a central flux in both velocity components: Order of convergence as a function of time, obtained from simulations on grids with $25 \times 25$ and $50\times 50$ cells (left; using $50\times 50$ and $100 \times 100$ on the right). One observes that with a $Q^1$ approximation the stationary state is only first-order accurate; third order as expected is observed for $Q^2$. For $Q^3$, again, the observed order is one less than expected. \emph{Right}: DG with a \fluxlowmach~flux: Order of convergence as a function of time, obtained from simulations on grids with $25 \times 25$ and $50\times 50$ cells. No loss of accuracy at steady state is observed.}
 \label{fig:central0201-order}
\end{figure}

Since no diffusion is added on the velocity, and since using centered discretizations is a popular cure for low Mach compliance / stationarity preservation for first-order finite volume / finite difference methods, one might expect that no loss of order would be observed with the central flux. Rather surprisingly, for $Q^1$ approximations central flux displays only first-order accuracy experimentally (Figure \ref{fig:central0201-order}, left). This is confirmed by the theoretical analysis, see Example \ref{example:centralker}. Thus, DG with central flux for $K=1$ is indeed at most first order accurate at stationary state.

For $K=2$ the kernel of $\mathcal E$ is given in Section \ref{app:ssec:central}. One finds this time, that as long as $a_6 = a_7 = 0$ in the selectively optimal projector given by
\begin{align}
  \pi_{\vec k} \hat q &:= \left(q_{v^{(22)}} + a_0 \Delta x^3\right ) w_1 + \left(q_{v^{(21)}} + a_1 \Delta x^3\right ) w_2 + \left(q_{v^{(20)}} + a_2 \Delta x^3\right ) w_3 \\ &\qquad \qquad \nonumber + \left(q_{v^{(12)}} + a_3 \Delta x^2\right ) w_4 + \ldots + \left(q_{v^{(00)}} + a_8 \Delta x^2\right ) w_9
\end{align}
then 
\begin{align}
 \widehat{\dof}(\hat Q) - \pi_{\vec k}\widehat{\dof}(\hat Q)  \in \mathcal O(\Delta x^3)
\end{align}

Thus, theoretically it is possible that DG with central flux achieves optimal (third-order) accuracy. This is confirmed experimentally. 

\begin{figure}
 \centering
 \includegraphics[width=0.45\textwidth]{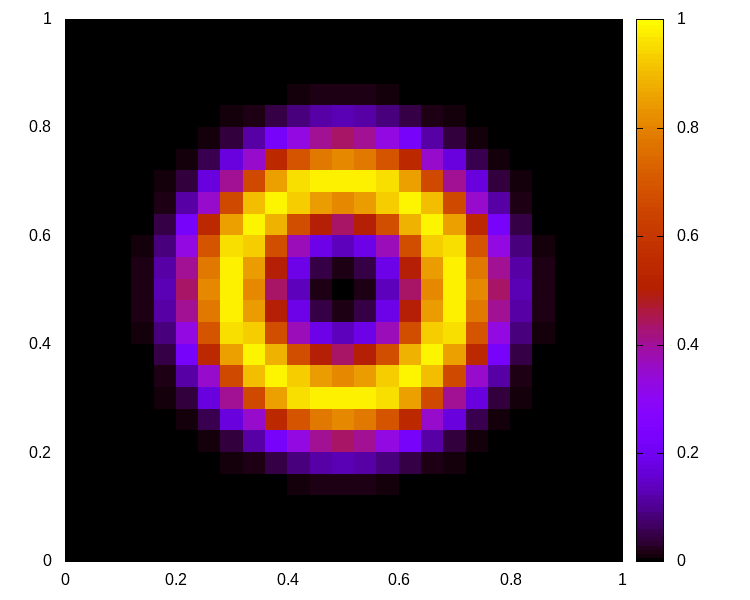} \includegraphics[width=0.45\textwidth]{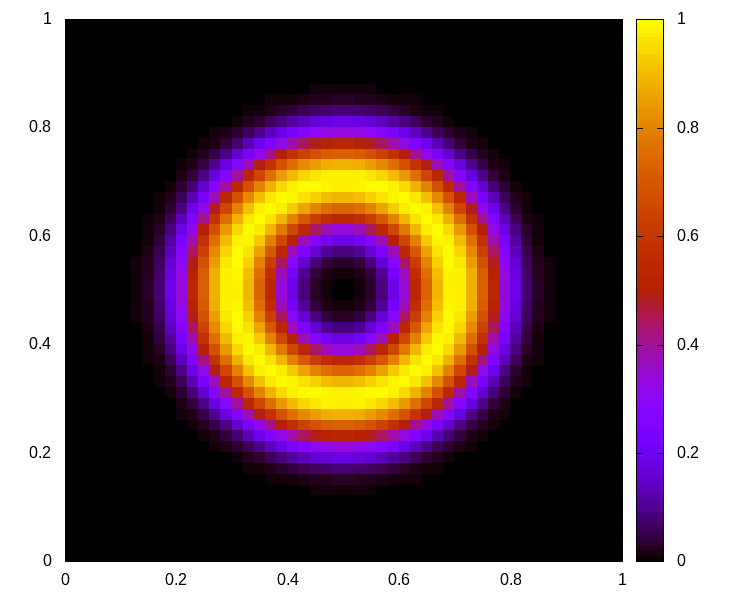} \\
 \includegraphics[width=0.45\textwidth]{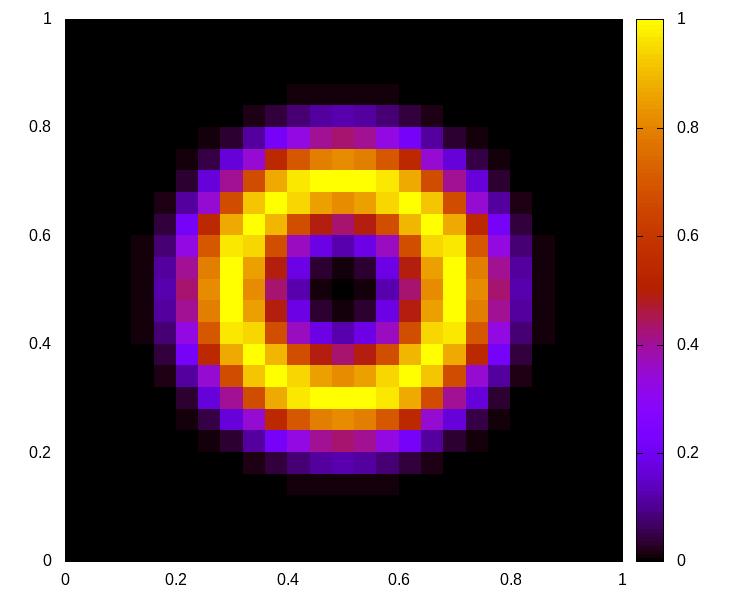} \includegraphics[width=0.45\textwidth]{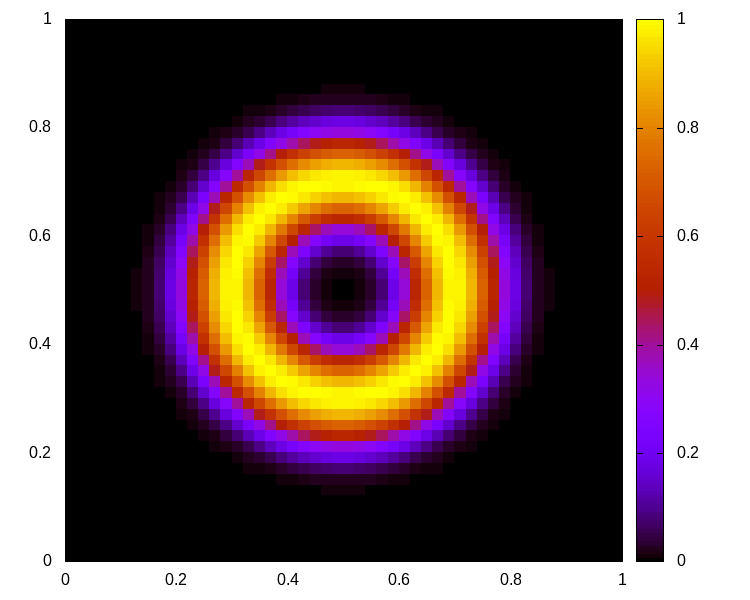} \\
 \includegraphics[width=0.45\textwidth]{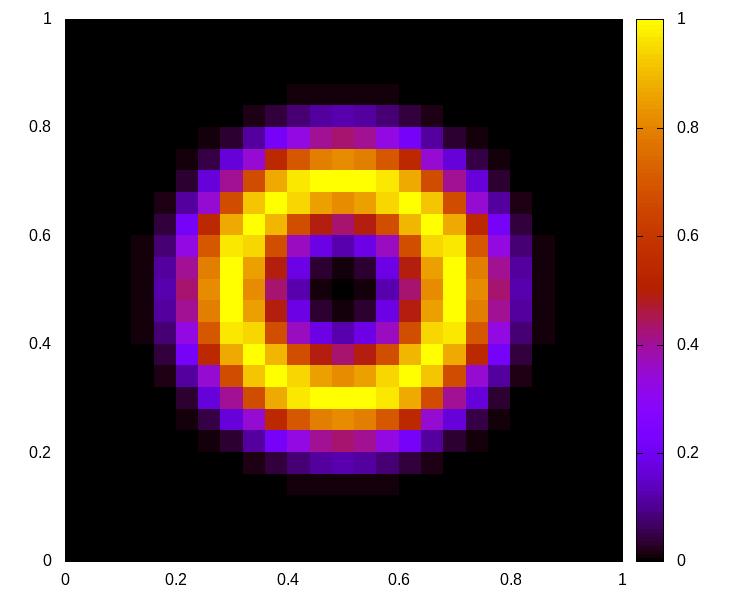} \includegraphics[width=0.45\textwidth]{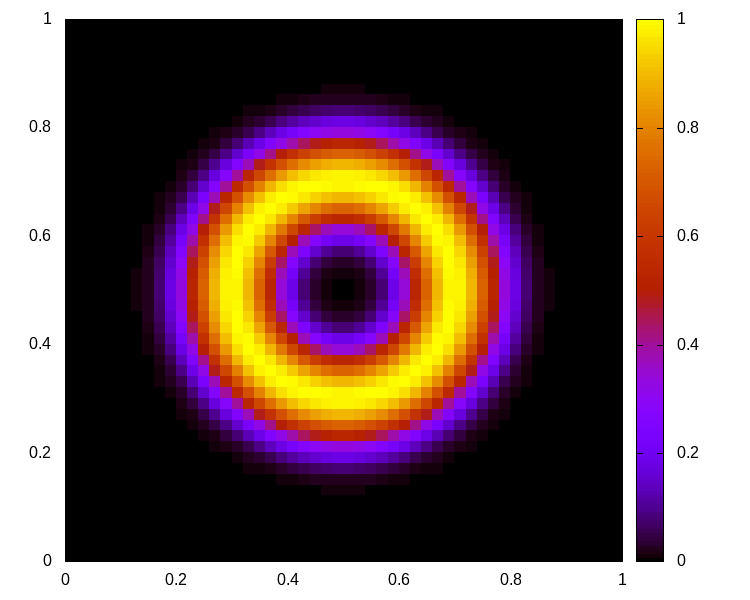} 
\caption{Central DG for linear acoustics; setup \eqref{eq:acousticvortex1}--\eqref{eq:acousticvortex2}. Numerical solution at cell center at time $t=300$ is shown. \emph{Top}: $K=1$. \emph{Center}: $K=2$. \emph{Bottom}: $K=3$. \emph{Left}: Grid of $25 \times 25$. \emph{Right}: Grid of $50 \times 50$.}
 \label{fig:central}
\end{figure}

\begin{figure}
 \centering
\includegraphics[width=0.32\textwidth]{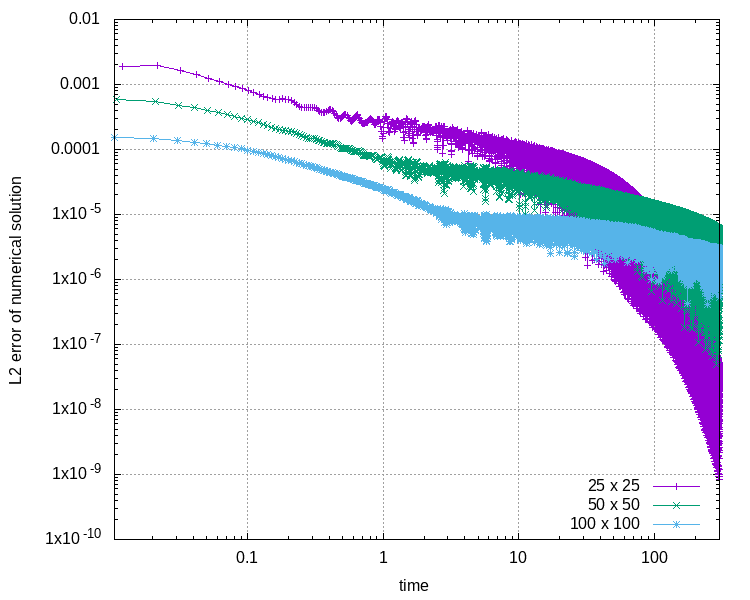}
 \hfill \includegraphics[width=0.32\textwidth]{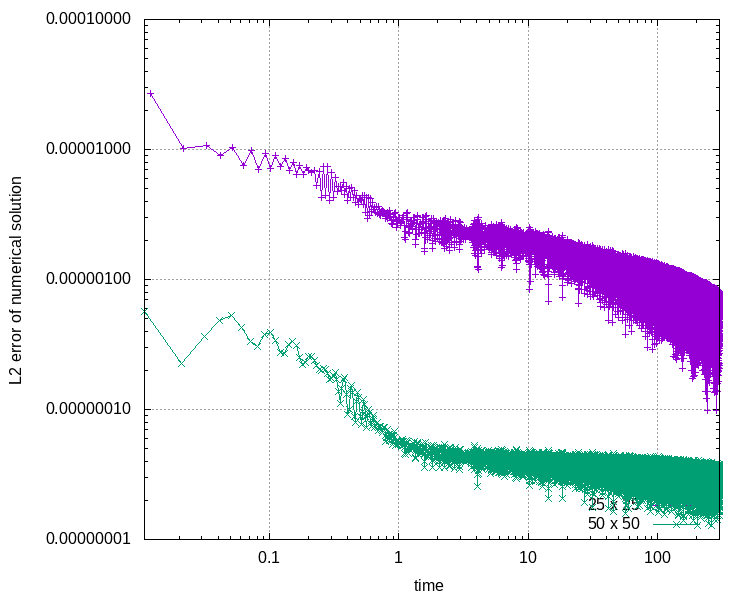}
 \hfill \includegraphics[width=0.32\textwidth]{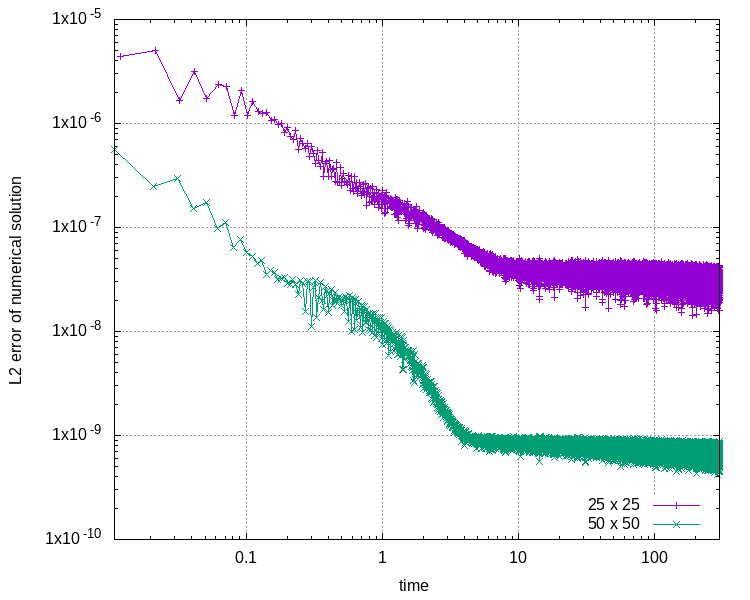}
 \caption{Error of the pressure for central DG as a function of time. \emph{Left}: $K=1$. \emph{Center}: $K=2$ \emph{Right}: $K=3$.}
 \label{fig:centralconvergencepressure}
\end{figure}

Finally, for $K=3$ (kernel of $\mathcal E$ is given in Section \ref{app:ssec:central}) one again finds that no accuracy better than third-order can be achieved, i.e. one order less than optimal. The calculation is omitted due to its length. These findings are reminiscent of \cite{liu20,hindenlang20}, where sub-optimal convergence rates are reported for DG with entropy conserving/central flux for polynomial approximations of odd degree, albeit without reference to steady state. Figure \ref{fig:central} shows the numerical solution for various values of $K$ and Figure \ref{fig:centralconvergencepressure} shows how the error of the pressure decays towards machine zero.

\subsubsection{\fluxlowmach~flux}

The last numerical flux to be investigated corresponds to the following choices of the stabilization matrices in \eqref{eq:numfluxx}--\eqref{eq:numfluxy}:
\begin{align}
 D_x &= \left( \begin{array}{ccc}  0 & 0 & 1 \\ 0 &0 &0 \\-1 & 0 &2  \end{array} \right ) &
 D_y &= \left( \begin{array}{ccc}  0 & 0 & 0 \\ 0 &0 &1 \\0 & -1 &2  \end{array} \right ) \label{eq:acoustic0201}
\end{align}
It is this flux from which the low Mach flux \eqref{eq:euler0201} for the Euler equations was obtained in \cite{barsukow18thesis}.
The modified equation associated to this method reads
\begin{align}
 \del_t u + \del_x p &= \Delta x \del_x^2 p \\
 \del_t v + \del_x p &= \Delta y \del_y^2 p\\
 \del_t p + \del_x u + \del_y v  &= -\Delta x \del_x^2 u - \Delta y \del_y^2 v + 2(\Delta x \del_x^2 p + \Delta y \del_y^2 p)
\end{align}
One can easily conjecture that therefore, the discrete stationary states are governed by $p=0$ and some non-standard discretization of the velocity divergence. As is shown in \cite{barsukow18thesis}, where this choice of flux has been introduced, it indeed leads to a low Mach compliant / stationarity preserving method when used inside a first-order finite volume / finite difference method, while being stable under explicit time discretization. Most importantly, it is stable under a CFL condition that involves $\Delta x$, as is customary for hyperbolic problems, and not $\Delta x^2$, as some other low Mach fixes (\cite{turkel87,birken05,barsukow16,bouchut20}). Figure \ref{fig:central0201-order} (right) shows that for DG with $K=0, \dots, 3$ it does not display a loss of order of accuracy. The corresponding evolution matrix $\mathcal E$ possesses a kernel of dimension $(K+1)^2$, just as the central method discussed above. For $K=1$ one finds, for example,

\begin{figure}
 \centering
 \includegraphics[width=0.45\textwidth]{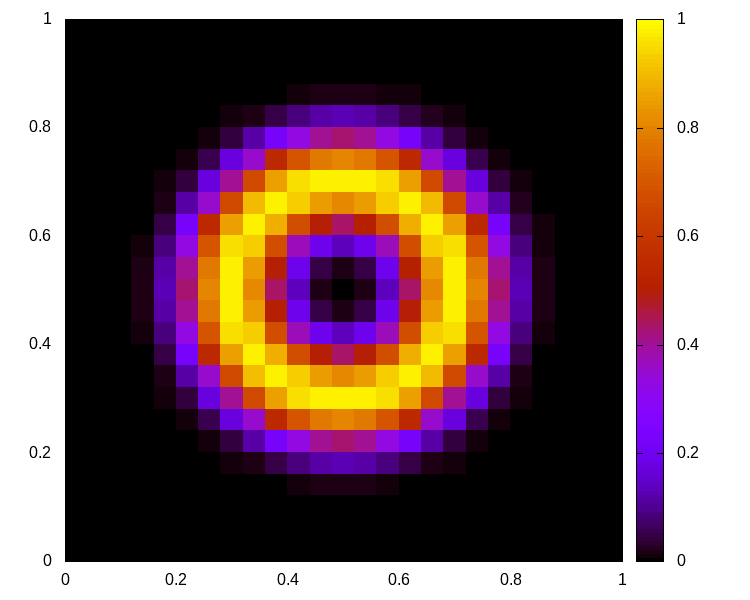} \includegraphics[width=0.45\textwidth]{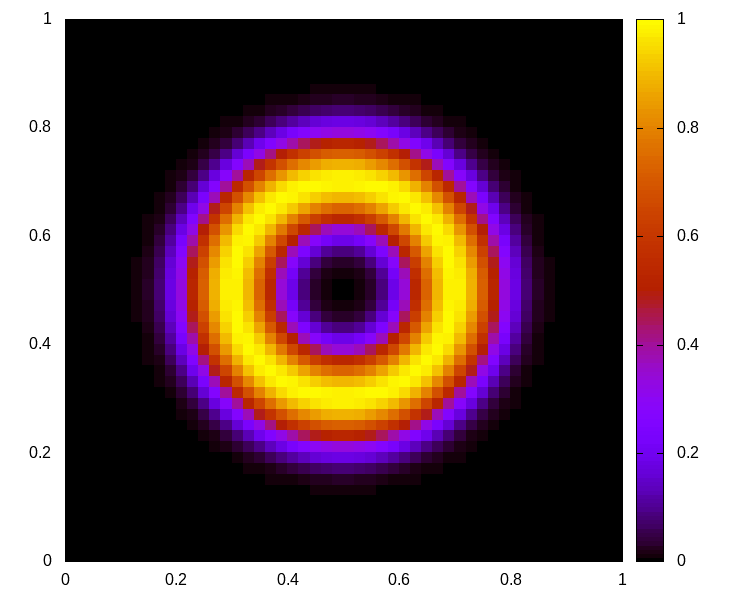} \\
 \includegraphics[width=0.45\textwidth]{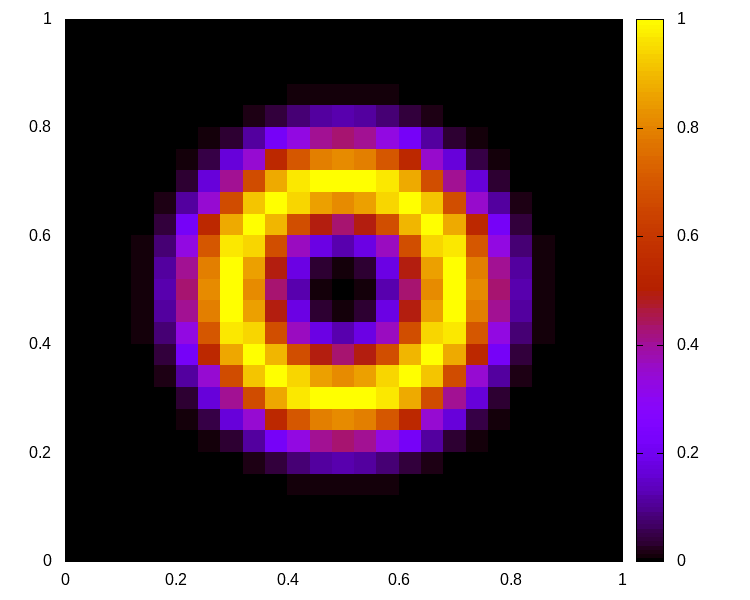} \includegraphics[width=0.45\textwidth]{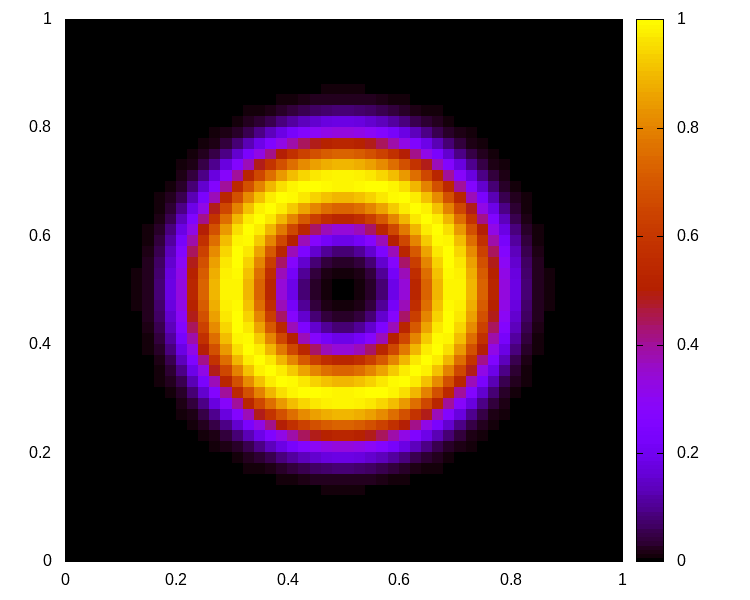} \\
  \includegraphics[width=0.45\textwidth]{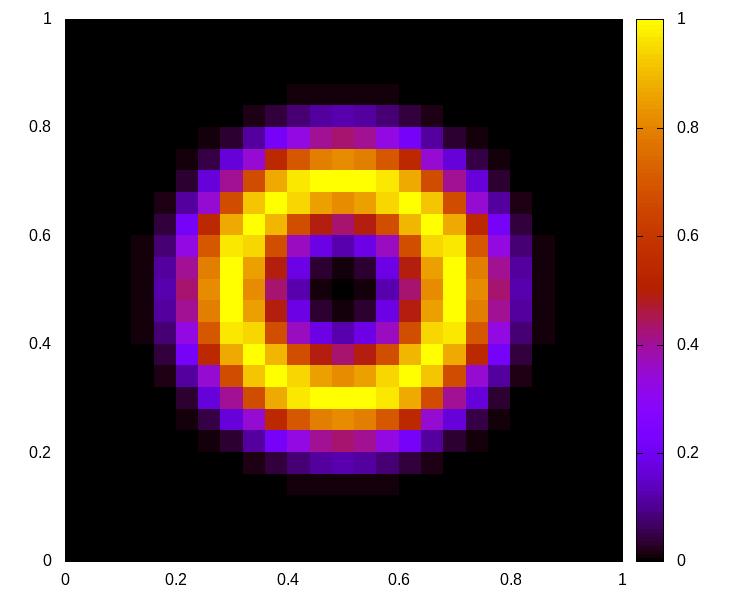} \includegraphics[width=0.45\textwidth]{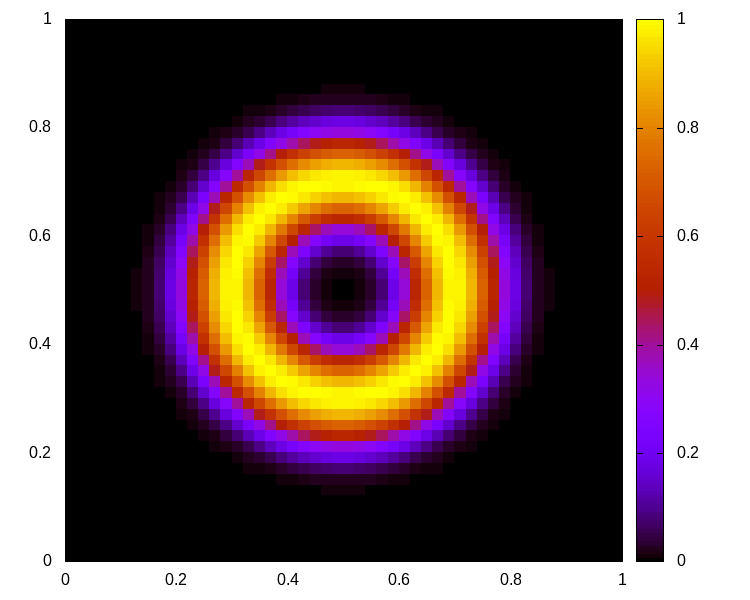}
\caption{DG with the \fluxlowmach~flux \eqref{eq:acoustic0201} for linear acoustics; setup \eqref{eq:acousticvortex1}--\eqref{eq:acousticvortex2}. Numerical solution at cell center at time $t=300$ is shown. \emph{Top}:  $K=1$. \emph{Center}: $K=2$. \emph{Bottom}: $K=3$. \emph{Left}: Grid of $25 \times 25$. \emph{Right}: Grid of $50 \times 50$.}
 \label{fig:0201}
\end{figure}

\begin{figure}
 \centering
 \includegraphics[width=0.32\textwidth]{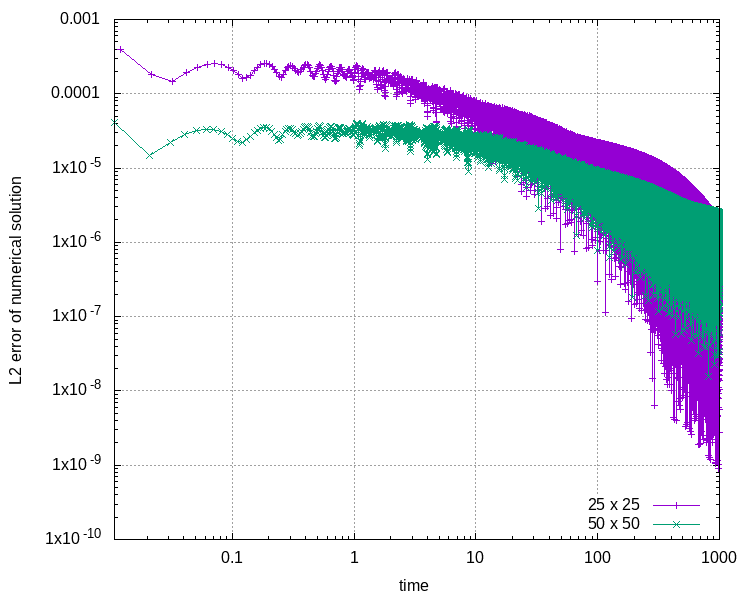} \hfill 
 \includegraphics[width=0.32\textwidth]{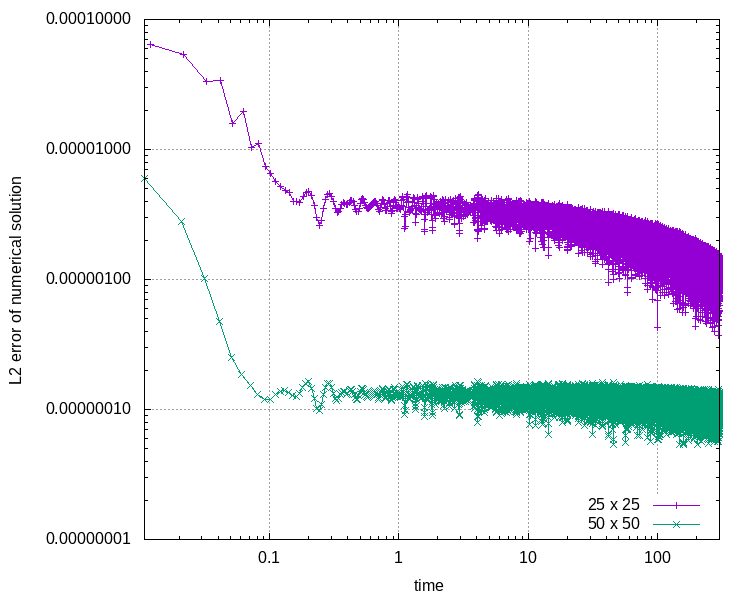} \hfill 
 \includegraphics[width=0.32\textwidth]{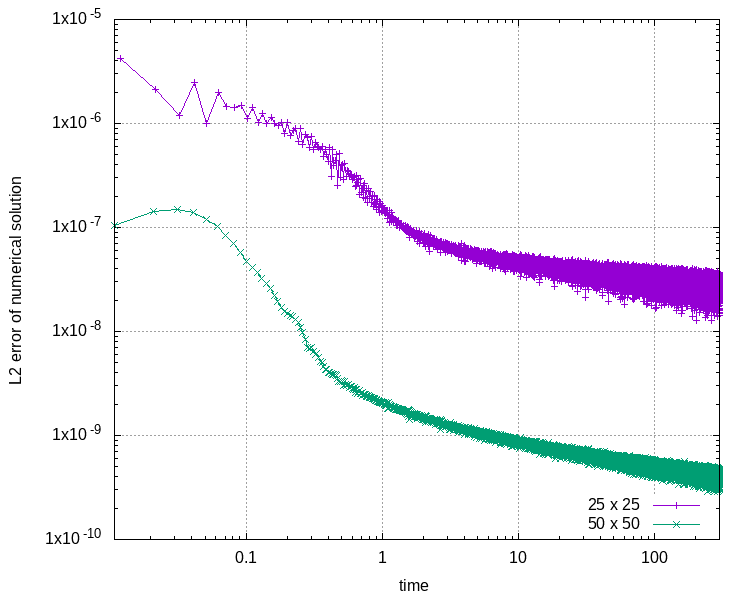} 
 \caption{DG with the \fluxlowmach~flux \eqref{eq:acoustic0201}: Error in the pressure as a function of time, obtained from simulations on grids with $25 \times 25$ and $50\times 50$ cells. The functions are highly oscillative and are partly on top of each other. The purpose of the plot is to show the general, qualitative behaviour, dominated by a near-exponential decay of the pressure, meaning that the solution is approaching the discrete steady state. \emph{Left}: $K=1$. \emph{Center}: $K=2$ \emph{Right}: $K=3$.}
 \label{fig:0201pressure}
\end{figure}

{\footnotesize
\begin{align}
w_1 &= \left (\frac{ \Delta x(1-  t_y)}{2  \Delta y},-\frac{\sqrt{3}  \Delta x (1+ t_y)}{2  \Delta y},-\frac{ \Delta x
(t_y-1)}{2 \sqrt{3}  \Delta y},-\frac{ \Delta x (1+ t_y)}{2  \Delta y},0,0,0,1,\boxed{4 \times 0}\right )^\text{T} \\ 
w_2 &= \left (\frac{ \Delta x (t_y-1)}{2
\sqrt{3}  \Delta y},\frac{ \Delta x (t_y-1)}{2  \Delta y},\frac{ \Delta x (t_y-1)}{6  \Delta y},\frac{ \Delta x (t_y-1)}{2 \sqrt{3}
 \Delta y},0,0,1,0,\boxed{4 \times 0}\right )^\text{T} \\ 
 w_3 &= \left (\frac{\sqrt{3}  \Delta x (1+ t_x) (t_y-1)}{2  \Delta y (t_x-1)},\frac{3  \Delta x (1+ t_x)
(1+ t_y)}{2  \Delta y (t_x-1)},\frac{ \Delta x (t_y-1)}{2  \Delta y},\frac{\sqrt{3}  \Delta x (1+ t_y)}{2  \Delta y},0,1,0,0,\boxed{4 \times 0}\right )^\text{T} \\ 
w_4 &= \left (-\frac{ \Delta x
(1+ t_x) (t_y-1)}{2  \Delta y (t_x-1)},-\frac{\sqrt{3}  \Delta x (1+ t_x) (t_y-1)}{2  \Delta y (t_x-1)},-\frac{ \Delta x
(t_y-1)}{2 \sqrt{3}  \Delta y},-\frac{ \Delta x (t_y-1)}{2  \Delta y},1,0,0,0,\boxed{4 \times 0}\right )^\text{T} 
\end{align}}
While having the same dimension, the elements in the kernel are different. This time they are such that the optimal order of accuracy
\begin{align}
 \widehat{\dof}(\hat Q) - \pi_{\vec k}\widehat{\dof}(\hat Q)  \in \mathcal O(\Delta x^2)
\end{align}
can be achieved with the following selectively optimal projector
\begin{align*}
  \pi_{\vec k} \hat q &:= \left(q_{v^{(11)}} + a_0 \Delta x^2\right ) w_1 + \left(q_{v^{(10)}} + a_1 \Delta x^2\right ) w_2 + \left(q_{v^{(01)}} + a_2 \Delta x^2\right ) w_3 \\&\qquad+ \left(q_{v^{(00)}} + a_3 \Delta x^2\right ) w_4
\end{align*}
as long as $a_2 = \frac{k_x k_y^2}{12 \sqrt{3}}$. One clearly observes the optimal order of accuracy in the experiment.

The kernels for other choices of $K$ are given in Appendix \ref{app:ssec:0201}. One similarly finds that it is possible to find selectively optimal projectors that lead to optimal order of accuracy $\mathcal O(\Delta x^{K+1})$. In particular, for $K=2$, such a projector is given by 
\begin{align}
  \pi_{\vec k} \hat q &:= \left(q_{v^{(22)}} + a_0 \Delta x^3\right ) w_1 + \left(q_{v^{(21)}} + a_1 \Delta x^3\right ) w_2 + \left(q_{v^{(20)}} + a_2 \Delta x^3\right ) w_3 \\ &\qquad \qquad \nonumber + \left(q_{v^{(12)}} + a_3 \Delta x^3\right ) w_4 + \ldots + \left(q_{v^{(00)}} + a_8 \Delta x^3\right ) w_9
\end{align}
with $a_6 = \frac{\ii k_x k_y^3}{120 \sqrt{5}}$ and for $K=3$ by
\begin{align}
  \pi_{\vec k} \hat q &:= \left(q_{v^{(33)}} + a_0 \Delta x^4\right ) w_1 + \left(q_{v^{(32)}} + a_1 \Delta x^4\right ) w_2 + \left(q_{v^{(31)}} + a_2 \Delta x^4\right ) w_3 \\ & \qquad \nonumber + \left(q_{v^{(30)}} + a_3 \Delta x^4\right ) w_4 +\left(q_{v^{(23)}} + a_4 \Delta x^4\right ) w_5 + \ldots + \left(q_{v^{(00)}} + a_{15} \Delta x^4\right ) w_{16}
\end{align}
with $a_{13} = a_{14} = 0$ and $a_{12} = -\frac{k_xk_y^4}{1680 \sqrt{7}}$.

Figure \ref{fig:0201} shows the numerical solution for various values of $K$ and Figure \ref{fig:0201pressure} shows how the error of the pressure decays towards machine zero.

\section{Conclusions and outlook}

The first striking observation is that the Discontinuous Galerkin methods of polynomial degree $K$ higher than first are structure preserving for linear acoustics on Cartesian grids upon usage of the upwind flux, and even the Rusanov flux for $K \geq 2$. An analogous statement holds for the low Mach number behaviour in the context of the Euler equations. This is good news, i.e. in principle no modifications (``low Mach fixes'') are necessary. The situation thus is significantly better than in the first-order case.

The discrete stationary states of the method are not always discretizations of those of the PDE. The statement that the situation becomes better with increasing polynomial degree is not just about them naturally becoming higher-order approximations, but it is a qualitative change that happens suddenly: Above a certain polynomial degree, DG becomes stationarity preserving. For lower polynomial degrees, the method simply does not have any non-trivial discrete stationary states. 

When they suddenly appear for sufficiently high $K$, they often turn out to have lower order of accuracy than expected for degree $K$, which is the second observation. Different numerical fluxes used in the DG method lead to vastly different orders of accuracy at stationary state. Upwind flux (for $K\geq 1$) and Rusanov flux (for $K\geq 2$) lose one order. A rather surprising result is that the central flux loses one order for odd $K$. The numerical flux \eqref{eq:acoustic0201} associated to a low Mach fix is the only numerical flux among those considered here, that retains its optimal order of accuracy at steady state. 

One can thus summarize that while low Mach fixes are not generally necessary to achieve a consistent approximation of steady states for high-order DG methods, low-Mach-fix-type numerical fluxes can help retain the design order of accuracy. To investigate whether this is the case for low Mach fixes in general is subject of future work. Stationarity preservation also implies that (for polynomial degrees mentioned above) DG is vorticity preserving. This interesting finding will also be exploited in future.

Stationarity preservation for linear acoustics is tightly linked to low Mach compliance for the Euler equations, and the theoretical and experimental results for linear acoustics presented here give a theoretical underpinning for the observed low Mach compliance of DG for the Euler equations. It would be of interest to study the loss of accuracy for the Euler equations in the low Mach number limit, that one thus would expect based on the results for linear acoustics. This, however, does not seem easy because of the high computational cost associated with such an undertaking. Also, it is not clear whether one can hope for as clear results as the ones for linear acoustics: If one thinks of the Euler equations as being composed of an advective and an acoustic operator, the loss of order would be restricted to the acoustic operator, while one will probably expect the advective operator to retain its design order. While some loss of order might be measurable, the value will possibly depend on the details of the setup.

\section*{Acknowledgement}
The author would like to express his gratitude to the Institute for Mathematics of the University of Zurich, whose computational facilities he was 
able to continue using during the extended time it took to bring this work to completion.

\newcommand{\etalchar}[1]{$^{#1}$}

\appendix

\begin{landscape}
\section{Bases of the kernel of $\mathcal E$}

Recall that the notation $\boxed{N \times 0}$ abbreviates $N$ vanishing entries of the vectors in components associated to the pressure.

\subsection{Upwind flux} \label{app:ssec:upwindK3}

The kernel of $\mathcal E$ for $K=3$ is spanned by the following 9 vectors: {\tiny
\begin{align*}
\left (0,0,0,0,\frac{\sqrt{\frac{7}{15}}  \Delta x}{ \Delta y},0,\frac{\sqrt{\frac{7}{3}}  \Delta x}{ \Delta y},0,0,0,0,0,-\frac{ \Delta x}{\sqrt{5}
 \Delta y},0,-\frac{ \Delta x}{ \Delta y},0,0,0,0,0,0,0,0,0,\frac{\sqrt{7} (1+ t_y)}{t_y-1},0,0,1,0,0,0,0,\boxed{16 \times 0}\right )^\text{T} \\ \left (0,0,0,0,0,\frac{ \Delta x}{ \Delta y},0,0,0,0,0,0,0,-\frac{\sqrt{\frac{3}{7}}
 \Delta x}{ \Delta y},0,0,0,0,0,0,0,0,0,0,-\sqrt{5},0,1,0,0,0,0,0,\boxed{16 \times 0}\right )^\text{T} \\ \left (0,0,0,0,\frac{ \Delta x}{\sqrt{5}
 \Delta y},0,0,0,0,0,0,0,-\frac{\sqrt{\frac{3}{35}}  \Delta x}{ \Delta y},0,0,0,0,0,0,0,0,0,0,0,\frac{\sqrt{3} (1+ t_y)}{t_y-1},1,0,0,0,0,0,0,\boxed{16 \times 0}\right )^\text{T} \\ \left (\frac{\sqrt{\frac{7}{3}}
 \Delta x}{ \Delta y},0,\frac{\sqrt{\frac{35}{3}}  \Delta x}{ \Delta y},0,0,0,0,0,-\frac{\sqrt{\frac{7}{15}}  \Delta x}{ \Delta y},0,-\frac{\sqrt{\frac{7}{3}}
 \Delta x}{ \Delta y},0,0,0,0,0,0,0,0,0,\frac{\sqrt{7} (1+ t_y)}{t_y-1},0,0,1,0,0,0,0,0,0,0,0,\boxed{16 \times 0}\right )^\text{T} \\ \left (0,\frac{\sqrt{5}
 \Delta x}{ \Delta y},0,0,0,0,0,0,0,-\frac{ \Delta x}{ \Delta y},0,0,0,0,0,0,0,0,0,0,-\sqrt{5},0,1,0,0,0,0,0,0,0,0,0,\boxed{16 \times 0}\right )^\text{T} \\ \left (\frac{ \Delta x}{ \Delta y},0,0,0,0,0,0,0,-\frac{ \Delta x}{\sqrt{5}
 \Delta y},0,0,0,0,0,0,0,0,0,0,0,\frac{\sqrt{3} (1+ t_y)}{t_y-1},1,0,0,0,0,0,0,0,0,0,0,\boxed{16 \times 0}\right )^\text{T} \\ \left (-\frac{\sqrt{7}
 \Delta x (1+ t_x)}{ \Delta y (t_x-1)},0,-\frac{\sqrt{35}  \Delta x (1+ t_x)}{ \Delta y (t_x-1)},0,-\frac{\sqrt{\frac{7}{3}}
 \Delta x}{ \Delta y},0,-\frac{\sqrt{\frac{35}{3}}  \Delta x}{ \Delta y},0,0,0,0,0,0,0,0,0,\frac{\sqrt{7} (1+ t_y)}{t_y-1},0,0,1,0,0,0,0,0,0,0,0,0,0,0,0,\boxed{16 \times 0}\right )^\text{T} \\ \left (0,-\frac{\sqrt{15}
 \Delta x (1+ t_x)}{ \Delta y (t_x-1)},0,0,0,-\frac{\sqrt{5}  \Delta x}{ \Delta y},0,0,0,0,0,0,0,0,0,0,-\sqrt{5},0,1,0,0,0,0,0,0,0,0,0,0,0,0,0,\boxed{16 \times 0}\right )^\text{T} \\ \left (-\frac{\sqrt{3}
 \Delta x (1+ t_x)}{ \Delta y (t_x-1)},0,0,0,-\frac{ \Delta x}{ \Delta y},0,0,0,0,0,0,0,0,0,0,0,\frac{\sqrt{3} (1+ t_y)}{t_y-1},1,0,0,0,0,0,0,0,0,0,0,0,0,0,0,\boxed{16 \times 0}\right )^\text{T}
\end{align*}}

\subsection{Central flux} \label{app:ssec:central}

The kernel of $\mathcal E$ for $K=2$ is spanned by the following 9 vectors: {\tiny
\begin{align*}
w_1 &= \Bigg (0,0,0,\frac{ \Delta x \left(t_y^2-1\right)}{4 \sqrt{3}  \Delta y  t_y},\frac{ \Delta x (1+ t_y)^2}{4  \Delta y  t_y},\frac{\sqrt{\frac{5}{3}}  \Delta x \left(t_y^2-1\right)}{4  \Delta y  t_y},\frac{ \Delta x (t_x-1) \left(t_y^2-1\right)}{4 \sqrt{5}  \Delta y (1+ t_x)  t_y},\frac{\sqrt{\frac{3}{5}}  \Delta x (t_x-1) (1+ t_y)^2}{4  \Delta y (1+ t_x)  t_y},\frac{ \Delta x (t_x-1) \left(t_y^2-1\right)}{4  \Delta y (1+ t_x)  t_y},   0,0,0,0,0,0,0,0,1,\boxed{9 \times 0}\Bigg )^\text{T} \\ 
w_2 &= \Bigg (0,0,0,-\frac{ \Delta x (t_y-1)^2}{4 \sqrt{5}  \Delta y  t_y},-\frac{\sqrt{\frac{3}{5}}  \Delta x \left(t_y^2-1\right)}{4  \Delta y  t_y},-\frac{ \Delta x (1+ t_y)^2}{4  \Delta y  t_y},-\frac{\sqrt{3}  \Delta x (t_x-1) (t_y-1)^2}{20  \Delta y (1+ t_x)  t_y},-\frac{3  \Delta x (t_x-1) \left(t_y^2-1\right)}{20  \Delta y (1+ t_x)  t_y},-\frac{\sqrt{\frac{3}{5}}  \Delta x (t_x-1) (1+ t_y)^2}{4  \Delta y (1+ t_x)  t_y}, \\\nonumber &\hspace{12cm}0,0,0,0,0,0,0,1,0,\boxed{9 \times 0}\Bigg )^\text{T} \\ 
w_3 &= \Bigg (0,0,0,\frac{ \Delta x \left(t_y^2-1\right)}{4 \sqrt{15}  \Delta y  t_y},\frac{ \Delta x (t_y-1)^2}{4 \sqrt{5}  \Delta y  t_y},\frac{ \Delta x \left(t_y^2-1\right)}{4 \sqrt{3}  \Delta y  t_y},\frac{ \Delta x (t_x-1) \left(t_y^2-1\right)}{20  \Delta y (1+ t_x)  t_y},\frac{\sqrt{3}  \Delta x (t_x-1) (t_y-1)^2}{20  \Delta y (1+ t_x)  t_y},\frac{ \Delta x (t_x-1) \left(t_y^2-1\right)}{4 \sqrt{5}  \Delta y (1+ t_x)  t_y},  0,0,0,0,0,0,1,0,0,\boxed{9 \times 0}\Bigg )^\text{T} \\ 
w_4 &= \Bigg (\frac{\sqrt{\frac{5}{3}}  \Delta x \left(t_y^2-1\right)}{4  \Delta y  t_y},\frac{\sqrt{5}  \Delta x (1+ t_y)^2}{4  \Delta y  t_y},\frac{5  \Delta x \left(t_y^2-1\right)}{4 \sqrt{3}  \Delta y  t_y},0,0,0,-\frac{ \Delta x \left(t_y^2-1\right)}{4 \sqrt{3}  \Delta y  t_y},-\frac{ \Delta x (1+ t_y)^2}{4  \Delta y  t_y},-\frac{\sqrt{\frac{5}{3}}  \Delta x \left(t_y^2-1\right)}{4  \Delta y  t_y},  0,0,0,0,0,1,0,0,0,\boxed{9 \times 0}\Bigg )^\text{T} \\ 
w_5 &= \Bigg (-\frac{ \Delta x (t_y-1)^2}{4  \Delta y  t_y},-\frac{\sqrt{3}  \Delta x \left(t_y^2-1\right)}{4  \Delta y  t_y},-\frac{\sqrt{5}  \Delta x (1+ t_y)^2}{4  \Delta y  t_y},0,0,0,\frac{ \Delta x (t_y-1)^2}{4 \sqrt{5}  \Delta y  t_y},\frac{\sqrt{\frac{3}{5}}  \Delta x \left(t_y^2-1\right)}{4  \Delta y  t_y},\frac{ \Delta x (1+ t_y)^2}{4  \Delta y  t_y},  0,0,0,0,1,0,0,0,0,\boxed{9 \times 0}\Bigg )^\text{T} \\ 
w_6 &= \Bigg (\frac{ \Delta x \left(t_y^2-1\right)}{4 \sqrt{3}  \Delta y  t_y},\frac{ \Delta x (t_y-1)^2}{4  \Delta y  t_y},\frac{\sqrt{\frac{5}{3}}  \Delta x \left(t_y^2-1\right)}{4  \Delta y  t_y},0,0,0,-\frac{ \Delta x \left(t_y^2-1\right)}{4 \sqrt{15}  \Delta y  t_y},-\frac{ \Delta x (t_y-1)^2}{4 \sqrt{5}  \Delta y  t_y},-\frac{ \Delta x \left(t_y^2-1\right)}{4 \sqrt{3}  \Delta y  t_y},  0,0,0,1,0,0,0,0,0,\boxed{9 \times 0}\Bigg )^\text{T} \\
w_7 &= \Bigg (-\frac{\sqrt{5}  \Delta x (1+ t_x) \left(t_y^2-1\right)}{4  \Delta y (t_x-1)  t_y},-\frac{\sqrt{15}  \Delta x (1+ t_x) (1+ t_y)^2}{4  \Delta y (t_x-1)  t_y},-\frac{5  \Delta x (1+ t_x) \left(t_y^2-1\right)}{4  \Delta y (t_x-1)  t_y},-\frac{\sqrt{\frac{5}{3}}  \Delta x \left(t_y^2-1\right)}{4  \Delta y  t_y},-\frac{\sqrt{5}  \Delta x (1+ t_y)^2}{4  \Delta y  t_y},-\frac{5  \Delta x \left(t_y^2-1\right)}{4 \sqrt{3}  \Delta y  t_y},0,0,0,  0,0,1,0,0,0,0,0,0,\boxed{9 \times 0}\Bigg )^\text{T} \\ 
w_8 &= \Bigg (\frac{\sqrt{3}  \Delta x (1+ t_x) (t_y-1)^2}{4  \Delta y (t_x-1)  t_y},\frac{3  \Delta x (1+ t_x) \left(t_y^2-1\right)}{4  \Delta y (t_x-1)  t_y},\frac{\sqrt{15}  \Delta x (1+ t_x) (1+ t_y)^2}{4  \Delta y (t_x-1)  t_y},\frac{ \Delta x (t_y-1)^2}{4  \Delta y  t_y},\frac{\sqrt{3}  \Delta x \left(t_y^2-1\right)}{4  \Delta y  t_y},\frac{\sqrt{5}  \Delta x (1+ t_y)^2}{4  \Delta y  t_y},0,0,0,  0,1,0,0,0,0,0,0,0,\boxed{9 \times 0}\Bigg )^\text{T} \\ 
w_9 &= \Bigg (-\frac{ \Delta x (1+ t_x) \left(t_y^2-1\right)}{4  \Delta y (t_x-1)  t_y},-\frac{\sqrt{3}  \Delta x (1+ t_x) (t_y-1)^2}{4  \Delta y (t_x-1)  t_y},-\frac{\sqrt{5}  \Delta x (1+ t_x) \left(t_y^2-1\right)}{4  \Delta y (t_x-1)  t_y},-\frac{ \Delta x \left(t_y^2-1\right)}{4 \sqrt{3}  \Delta y  t_y},-\frac{ \Delta x (t_y-1)^2}{4  \Delta y  t_y},-\frac{\sqrt{\frac{5}{3}}  \Delta x \left(t_y^2-1\right)}{4  \Delta y  t_y},0,0,0,  1,0,0,0,0,0,0,0,0,\boxed{9 \times 0} \Bigg)^\text{T} 
\end{align*}}

For $K=3$, the kernel of $\mathcal E$ is spanned by these 16 vectors: {\tiny
\begin{align*}
\left (0,0,0,0,0,0,0,0,-\frac{ \Delta x (t_y-1)^2}{4 \sqrt{5}  \Delta y  t_y},-\frac{\sqrt{\frac{3}{5}}  \Delta x \left(t_y^2-1\right)}{4
 \Delta y  t_y},-\frac{ \Delta x (t_y-1)^2}{4  \Delta y  t_y},-\frac{\sqrt{\frac{7}{5}}  \Delta x \left(t_y^2-1\right)}{4
 \Delta y  t_y},  -\frac{ \Delta x (1+ t_x) (t_y-1)^2}{4 \sqrt{7}  \Delta y (t_x-1)  t_y},-\frac{\sqrt{\frac{3}{7}}  \Delta x
(1+ t_x) \left(t_y^2-1\right)}{4  \Delta y (t_x-1)  t_y},-\frac{\sqrt{\frac{5}{7}}  \Delta x (1+ t_x) (t_y-1)^2}{4
 \Delta y (t_x-1)  t_y},-\frac{ \Delta x (1+ t_x) \left(t_y^2-1\right)}{4  \Delta y (t_x-1)  t_y},\right . \\ \nonumber \left.\phantom{\frac{\sqrt{\frac12}}2} 0,0,0,0,0,0,0,0,0,0,0,0,0,0,0,1,\boxed{16 \times 0}\right )^\text{T}
 \end{align*}
\begin{align*}
\left (0,0,0,0,0,0,0,0,\frac{ \Delta x
\left(t_y^2-1\right)}{4 \sqrt{7}  \Delta y  t_y},\frac{\sqrt{\frac{3}{7}}  \Delta x (1+ t_y)^2}{4  \Delta y  t_y},\frac{\sqrt{\frac{5}{7}}
 \Delta x \left(t_y^2-1\right)}{4  \Delta y  t_y},\frac{ \Delta x (t_y-1)^2}{4  \Delta y  t_y},  \frac{\sqrt{5}  \Delta x (1+ t_x)
\left(t_y^2-1\right)}{28  \Delta y (t_x-1)  t_y},\frac{\sqrt{15}  \Delta x (1+ t_x) (1+ t_y)^2}{28  \Delta y (t_x-1)
 t_y},\frac{5  \Delta x (1+ t_x) \left(t_y^2-1\right)}{28  \Delta y (t_x-1)  t_y},\frac{\sqrt{\frac{5}{7}}  \Delta x
(1+ t_x) (t_y-1)^2}{4  \Delta y (t_x-1)  t_y},\right . \\ \nonumber \left.\phantom{\frac{\sqrt{\frac12}}2} 0,0,0,0,0,0,0,0,0,0,0,0,0,0,1,0,\boxed{16 \times 0}\right )^\text{T}
\end{align*}
\begin{align*}
\left (0,0,0,0,0,0,0,0,-\frac{\sqrt{\frac{3}{35}}
 \Delta x (t_y-1)^2}{4  \Delta y  t_y},-\frac{3  \Delta x \left(t_y^2-1\right)}{4 \sqrt{35}  \Delta y  t_y},-\frac{\sqrt{\frac{3}{7}}
 \Delta x (1+ t_y)^2}{4  \Delta y  t_y},-\frac{\sqrt{\frac{3}{5}}  \Delta x \left(t_y^2-1\right)}{4  \Delta y  t_y},  -\frac{\sqrt{3}
 \Delta x (1+ t_x) (t_y-1)^2}{28  \Delta y (t_x-1)  t_y},-\frac{3  \Delta x (1+ t_x) \left(t_y^2-1\right)}{28
 \Delta y (t_x-1)  t_y},-\frac{\sqrt{15}  \Delta x (1+ t_x) (1+ t_y)^2}{28  \Delta y (t_x-1)  t_y},-\frac{\sqrt{\frac{3}{7}}
 \Delta x (1+ t_x) \left(t_y^2-1\right)}{4  \Delta y (t_x-1)  t_y},\right . \\ \nonumber \left.\phantom{\frac{\sqrt{\frac12}}2}  0,0,0,0,0,0,0,0,0,0,0,0,0,1,0,0,\boxed{16 \times 0}\right )^\text{T} 
 \end{align*}
\begin{align*}
\left (0,0,0,0,0,0,0,0,\frac{ \Delta x
\left(t_y^2-1\right)}{4 \sqrt{35}  \Delta y  t_y},\frac{\sqrt{\frac{3}{35}}  \Delta x (t_y-1)^2}{4  \Delta y  t_y},\frac{ \Delta x
\left(t_y^2-1\right)}{4 \sqrt{7}  \Delta y  t_y},\frac{ \Delta x (t_y-1)^2}{4 \sqrt{5}  \Delta y  t_y},  \frac{ \Delta x (1+ t_x)
\left(t_y^2-1\right)}{28  \Delta y (t_x-1)  t_y},\frac{\sqrt{3}  \Delta x (1+ t_x) (t_y-1)^2}{28  \Delta y (t_x-1)
 t_y},\frac{\sqrt{5}  \Delta x (1+ t_x) \left(t_y^2-1\right)}{28  \Delta y (t_x-1)  t_y},\frac{ \Delta x (1+ t_x)
(t_y-1)^2}{4 \sqrt{7}  \Delta y (t_x-1)  t_y},\right . \\ \nonumber \left.\phantom{\frac{\sqrt{\frac12}}2} 0,0,0,0,0,0,0,0,0,0,0,0,1,0,0,0,\boxed{16 \times 0}\right )^\text{T} 
\end{align*}
\begin{align*}
\left (0,0,0,0,-\frac{\sqrt{\frac{7}{15}}
 \Delta x (t_y-1)^2}{4  \Delta y  t_y},-\frac{\sqrt{\frac{7}{5}}  \Delta x \left(t_y^2-1\right)}{4  \Delta y  t_y},-\frac{\sqrt{\frac{7}{3}}
 \Delta x (t_y-1)^2}{4  \Delta y  t_y},-\frac{7  \Delta x \left(t_y^2-1\right)}{4 \sqrt{15}  \Delta y  t_y},  0,0,0,0,\frac{ \Delta x
(t_y-1)^2}{4 \sqrt{5}  \Delta y  t_y},\frac{\sqrt{\frac{3}{5}}  \Delta x \left(t_y^2-1\right)}{4  \Delta y  t_y},\frac{ \Delta x
(t_y-1)^2}{4  \Delta y  t_y},\frac{\sqrt{\frac{7}{5}}  \Delta x \left(t_y^2-1\right)}{4  \Delta y  t_y},\right . \\ \nonumber \left.\phantom{\frac{\sqrt{\frac12}}2} 0,0,0,0,0,0,0,0,0,0,0,1,0,0,0,0,\boxed{16 \times 0}\right )^\text{T} 
\end{align*}
\begin{align*}
\left (0,0,0,0,\frac{ \Delta x
\left(t_y^2-1\right)}{4 \sqrt{3}  \Delta y  t_y},\frac{ \Delta x (1+ t_y)^2}{4  \Delta y  t_y},\frac{\sqrt{\frac{5}{3}}  \Delta x
\left(t_y^2-1\right)}{4  \Delta y  t_y},\frac{\sqrt{\frac{7}{3}}  \Delta x (t_y-1)^2}{4  \Delta y  t_y},  0,0,0,0,-\frac{ \Delta x
\left(t_y^2-1\right)}{4 \sqrt{7}  \Delta y  t_y},-\frac{\sqrt{\frac{3}{7}}  \Delta x (1+ t_y)^2}{4  \Delta y  t_y},-\frac{\sqrt{\frac{5}{7}}
 \Delta x \left(t_y^2-1\right)}{4  \Delta y  t_y},-\frac{ \Delta x (t_y-1)^2}{4  \Delta y  t_y},\right . \\ \nonumber \left.\phantom{\frac{\sqrt{\frac12}}2} 0,0,0,0,0,0,0,0,0,0,1,0,0,0,0,0,\boxed{16 \times 0}\right )^\text{T}  
 \end{align*}
\begin{align*}
 \left (0,0,0,0,-\frac{ \Delta x
(t_y-1)^2}{4 \sqrt{5}  \Delta y  t_y},-\frac{\sqrt{\frac{3}{5}}  \Delta x \left(t_y^2-1\right)}{4  \Delta y  t_y},-\frac{ \Delta x
(1+ t_y)^2}{4  \Delta y  t_y},-\frac{\sqrt{\frac{7}{5}}  \Delta x \left(t_y^2-1\right)}{4  \Delta y  t_y},  0,0,0,0,\frac{\sqrt{\frac{3}{35}}
 \Delta x (t_y-1)^2}{4  \Delta y  t_y},\frac{3  \Delta x \left(t_y^2-1\right)}{4 \sqrt{35}  \Delta y  t_y},\frac{\sqrt{\frac{3}{7}}
 \Delta x (1+ t_y)^2}{4  \Delta y  t_y},\frac{\sqrt{\frac{3}{5}}  \Delta x \left(t_y^2-1\right)}{4  \Delta y  t_y},\right . \\ \nonumber \left.\phantom{\frac{\sqrt{\frac12}}2} 0,0,0,0,0,0,0,0,0,1,0,0,0,0,0,0,\boxed{16 \times 0}\right )^\text{T} 
 \end{align*}
\begin{align*}
\left (0,0,0,0,\frac{ \Delta x
\left(t_y^2-1\right)}{4 \sqrt{15}  \Delta y  t_y},\frac{ \Delta x (t_y-1)^2}{4 \sqrt{5}  \Delta y  t_y},\frac{ \Delta x \left(t_y^2-1\right)}{4
\sqrt{3}  \Delta y  t_y},\frac{\sqrt{\frac{7}{15}}  \Delta x (t_y-1)^2}{4  \Delta y  t_y},  0,0,0,0,-\frac{ \Delta x \left(t_y^2-1\right)}{4
\sqrt{35}  \Delta y  t_y},-\frac{\sqrt{\frac{3}{35}}  \Delta x (t_y-1)^2}{4  \Delta y  t_y},-\frac{ \Delta x \left(t_y^2-1\right)}{4
\sqrt{7}  \Delta y  t_y},-\frac{ \Delta x (t_y-1)^2}{4 \sqrt{5}  \Delta y  t_y},\right . \\ \nonumber \left.\phantom{\frac{\sqrt{\frac12}}2} 0,0,0,0,0,0,0,0,1,0,0,0,0,0,0,0,\boxed{16 \times 0}\right )^\text{T} 
\end{align*}
\begin{align*}
\left (-\frac{\sqrt{\frac{7}{3}}
 \Delta x (t_y-1)^2}{4  \Delta y  t_y},-\frac{\sqrt{7}  \Delta x \left(t_y^2-1\right)}{4  \Delta y  t_y},-\frac{\sqrt{\frac{35}{3}}
 \Delta x (t_y-1)^2}{4  \Delta y  t_y},-\frac{7  \Delta x \left(t_y^2-1\right)}{4 \sqrt{3}  \Delta y  t_y},  0,0,0,0,\frac{\sqrt{\frac{7}{15}}
 \Delta x (t_y-1)^2}{4  \Delta y  t_y},\frac{\sqrt{\frac{7}{5}}  \Delta x \left(t_y^2-1\right)}{4  \Delta y  t_y},\frac{\sqrt{\frac{7}{3}}
 \Delta x (t_y-1)^2}{4  \Delta y  t_y},\frac{7  \Delta x \left(t_y^2-1\right)}{4 \sqrt{15}  \Delta y  t_y}, 0,0,0,0,\right . \\ \nonumber \left.\phantom{\frac{\sqrt{\frac12}}2} 0,0,0,0,0,0,0,1,0,0,0,0,0,0,0,0,\boxed{16 \times 0}\right )^\text{T} 
 \end{align*}
\begin{align*}
\left (\frac{\sqrt{\frac{5}{3}}
 \Delta x \left(t_y^2-1\right)}{4  \Delta y  t_y},\frac{\sqrt{5}  \Delta x (1+ t_y)^2}{4  \Delta y  t_y},\frac{5  \Delta x
\left(t_y^2-1\right)}{4 \sqrt{3}  \Delta y  t_y},\frac{\sqrt{\frac{35}{3}}  \Delta x (t_y-1)^2}{4  \Delta y  t_y},  0,0,0,0,-\frac{ \Delta x
\left(t_y^2-1\right)}{4 \sqrt{3}  \Delta y  t_y},-\frac{ \Delta x (1+ t_y)^2}{4  \Delta y  t_y},-\frac{\sqrt{\frac{5}{3}}  \Delta x
\left(t_y^2-1\right)}{4  \Delta y  t_y},-\frac{\sqrt{\frac{7}{3}}  \Delta x (t_y-1)^2}{4  \Delta y  t_y}, 0,0,0,0,\right . \\ \nonumber \left.\phantom{\frac{\sqrt{\frac12}}2} 0,0,0,0,0,0,1,0,0,0,0,0,0,0,0,0,\boxed{16 \times 0}\right )^\text{T} 
\end{align*}
\begin{align*}
\left (-\frac{ \Delta x
(t_y-1)^2}{4  \Delta y  t_y},-\frac{\sqrt{3}  \Delta x \left(t_y^2-1\right)}{4  \Delta y  t_y},-\frac{\sqrt{5}  \Delta x
(1+ t_y)^2}{4  \Delta y  t_y},-\frac{\sqrt{7}  \Delta x \left(t_y^2-1\right)}{4  \Delta y  t_y},  0,0,0,0,\frac{ \Delta x (t_y-1)^2}{4
\sqrt{5}  \Delta y  t_y},\frac{\sqrt{\frac{3}{5}}  \Delta x \left(t_y^2-1\right)}{4  \Delta y  t_y},\frac{ \Delta x (1+ t_y)^2}{4
 \Delta y  t_y},\frac{\sqrt{\frac{7}{5}}  \Delta x \left(t_y^2-1\right)}{4  \Delta y  t_y},0,0,0,0,\right . \\ \nonumber \left.\phantom{\frac{\sqrt{\frac12}}2}  0,0,0,0,0,1,0,0,0,0,0,0,0,0,0,0,\boxed{16 \times 0}\right )^\text{T} 
 \end{align*}
\begin{align*}
 \left (\frac{ \Delta x
\left(t_y^2-1\right)}{4 \sqrt{3}  \Delta y  t_y},\frac{ \Delta x (t_y-1)^2}{4  \Delta y  t_y},\frac{\sqrt{\frac{5}{3}}  \Delta x
\left(t_y^2-1\right)}{4  \Delta y  t_y},\frac{\sqrt{\frac{7}{3}}  \Delta x (t_y-1)^2}{4  \Delta y  t_y},  0,0,0,0,-\frac{ \Delta x
\left(t_y^2-1\right)}{4 \sqrt{15}  \Delta y  t_y},-\frac{ \Delta x (t_y-1)^2}{4 \sqrt{5}  \Delta y  t_y},-\frac{ \Delta x
\left(t_y^2-1\right)}{4 \sqrt{3}  \Delta y  t_y},-\frac{\sqrt{\frac{7}{15}}  \Delta x (t_y-1)^2}{4  \Delta y  t_y}, 0,0,0,0,\right . \\ \nonumber \left.\phantom{\frac{\sqrt{\frac12}}2} 0,0,0,0,1,0,0,0,0,0,0,0,0,0,0,0,\boxed{16 \times 0}\right )^\text{T} 
\end{align*}
\begin{align*}
\left (\frac{\sqrt{7}
 \Delta x (1+ t_x) (t_y-1)^2}{4  \Delta y (t_x-1)  t_y},\frac{\sqrt{21}  \Delta x (1+ t_x) \left(t_y^2-1\right)}{4
 \Delta y (t_x-1)  t_y},\frac{\sqrt{35}  \Delta x (1+ t_x) (t_y-1)^2}{4  \Delta y (t_x-1)  t_y},\frac{7  \Delta x
(1+ t_x) \left(t_y^2-1\right)}{4  \Delta y (t_x-1)  t_y},  \frac{\sqrt{\frac{7}{3}}  \Delta x (t_y-1)^2}{4  \Delta y
 t_y},\frac{\sqrt{7}  \Delta x \left(t_y^2-1\right)}{4  \Delta y  t_y},\frac{\sqrt{\frac{35}{3}}  \Delta x (t_y-1)^2}{4  \Delta y
 t_y},\frac{7  \Delta x \left(t_y^2-1\right)}{4 \sqrt{3}  \Delta y  t_y}, 0,0,0,0,0,0,0,0,\right . \\ \nonumber \left.\phantom{\frac{\sqrt{\frac12}}2} 0,0,0,1,0,0,0,0,0,0,0,0,0,0,0,0,\boxed{16 \times 0}\right )^\text{T} 
 \end{align*}
\begin{align*}
\left (-\frac{\sqrt{5}
 \Delta x (1+ t_x) \left(t_y^2-1\right)}{4  \Delta y (t_x-1)  t_y},-\frac{\sqrt{15}  \Delta x (1+ t_x) (1+ t_y)^2}{4
 \Delta y (t_x-1)  t_y},-\frac{5  \Delta x (1+ t_x) \left(t_y^2-1\right)}{4  \Delta y (t_x-1)  t_y},-\frac{\sqrt{35}
 \Delta x (1+ t_x) (t_y-1)^2}{4  \Delta y (t_x-1)  t_y},  -\frac{\sqrt{\frac{5}{3}}  \Delta x \left(t_y^2-1\right)}{4
 \Delta y  t_y},-\frac{\sqrt{5}  \Delta x (1+ t_y)^2}{4  \Delta y  t_y},-\frac{5  \Delta x \left(t_y^2-1\right)}{4 \sqrt{3}
 \Delta y  t_y},-\frac{\sqrt{\frac{35}{3}}  \Delta x (t_y-1)^2}{4  \Delta y  t_y}, 0,0,0,0,0,0,0,0,\right . \\ \nonumber \left.\phantom{\frac{\sqrt{\frac12}}2} 0,0,1,0,0,0,0,0,0,0,0,0,0,0,0,0,\boxed{16 \times 0}\right )^\text{T} 
 \end{align*}
\begin{align*}
\left (\frac{\sqrt{3}
 \Delta x (1+ t_x) (t_y-1)^2}{4  \Delta y (t_x-1)  t_y},\frac{3  \Delta x (1+ t_x) \left(t_y^2-1\right)}{4  \Delta y
(t_x-1)  t_y},\frac{\sqrt{15}  \Delta x (1+ t_x) (1+ t_y)^2}{4  \Delta y (t_x-1)  t_y},\frac{\sqrt{21}  \Delta x
(1+ t_x) \left(t_y^2-1\right)}{4  \Delta y (t_x-1)  t_y},  \frac{ \Delta x (t_y-1)^2}{4  \Delta y  t_y},\frac{\sqrt{3}
 \Delta x \left(t_y^2-1\right)}{4  \Delta y  t_y},\frac{\sqrt{5}  \Delta x (1+ t_y)^2}{4  \Delta y  t_y},\frac{\sqrt{7}  \Delta x
\left(t_y^2-1\right)}{4  \Delta y  t_y}, 0,0,0,0,0,0,0,0,\right . \\ \nonumber \left.\phantom{\frac{\sqrt{\frac12}}2} 0,1,0,0,0,0,0,0,0,0,0,0,0,0,0,0,\boxed{16 \times 0}\right )^\text{T} 
\end{align*}
\begin{align*}
\left (-\frac{ \Delta x
(1+ t_x) \left(t_y^2-1\right)}{4  \Delta y (t_x-1)  t_y},-\frac{\sqrt{3}  \Delta x (1+ t_x) (t_y-1)^2}{4  \Delta y
(t_x-1)  t_y},-\frac{\sqrt{5}  \Delta x (1+ t_x) \left(t_y^2-1\right)}{4  \Delta y (t_x-1)  t_y},-\frac{\sqrt{7}
 \Delta x (1+ t_x) (t_y-1)^2}{4  \Delta y (t_x-1)  t_y},  -\frac{ \Delta x \left(t_y^2-1\right)}{4 \sqrt{3}  \Delta y
 t_y},-\frac{ \Delta x (t_y-1)^2}{4  \Delta y  t_y},-\frac{\sqrt{\frac{5}{3}}  \Delta x \left(t_y^2-1\right)}{4  \Delta y
 t_y},-\frac{\sqrt{\frac{7}{3}}  \Delta x (t_y-1)^2}{4  \Delta y  t_y},  0,0,0,0,0,0,0,0,\right . \\ \nonumber \left.\phantom{\frac{\sqrt{\frac12}}2} 1,0,0,0,0,0,0,0,0,0,0,0,0,0,0,0,\boxed{16 \times 0}\right )^\text{T} 
\end{align*}}

\subsection{\fluxlowmach~flux} \label{app:ssec:0201}

The kernel of $\mathcal E$ for $K=2$ is spanned by the following 9 vectors: {\tiny
\begin{align*}
w_1 = \left (0,0,0,\frac{ \Delta x (t_y-1)}{2 \sqrt{3}  \Delta y},\frac{ \Delta x (1+ t_y)}{2  \Delta y},\frac{\sqrt{\frac{5}{3}}  \Delta x (t_y-1)}{2  \Delta y},\frac{ \Delta x (t_y-1)}{2 \sqrt{5}  \Delta y},\frac{\sqrt{\frac{3}{5}}  \Delta x (1+ t_y)}{2  \Delta y},\frac{ \Delta x (t_y-1)}{2  \Delta y},0,0,0,0,0,0,0,0,1,\boxed{9 \times 0}\right )^\text{T} 
 \end{align*}
\begin{align*}
w_2 = \left (0,0,0,-\frac{ \Delta x (t_y-1)}{2 \sqrt{5}  \Delta y},-\frac{\sqrt{\frac{3}{5}}  \Delta x (1+ t_y)}{2  \Delta y},-\frac{ \Delta x (1+ t_y)}{2  \Delta y},-\frac{\sqrt{3}  \Delta x (t_y-1)}{10  \Delta y},-\frac{3  \Delta x (1+ t_y)}{10  \Delta y},-\frac{\sqrt{\frac{3}{5}}  \Delta x (1+ t_y)}{2  \Delta y},0,0,0,0,0,0,0,1,0,\boxed{9 \times 0}\right )^\text{T} 
 \end{align*}
\begin{align*}
w_3 = \left (0,0,0,\frac{ \Delta x (t_y-1)}{2 \sqrt{15}  \Delta y},\frac{ \Delta x (t_y-1)}{2 \sqrt{5}  \Delta y},\frac{ \Delta x (t_y-1)}{2 \sqrt{3}  \Delta y},\frac{ \Delta x (t_y-1)}{10  \Delta y},\frac{\sqrt{3}  \Delta x (t_y-1)}{10  \Delta y},\frac{ \Delta x (t_y-1)}{2 \sqrt{5}  \Delta y},0,0,0,0,0,0,1,0,0,\boxed{9 \times 0}\right )^\text{T} 
 \end{align*}
\begin{align*}
w_4 = \left (\frac{\sqrt{\frac{5}{3}}  \Delta x (t_y-1)}{2  \Delta y},\frac{\sqrt{5}  \Delta x (1+ t_y)}{2  \Delta y},\frac{5  \Delta x (t_y-1)}{2 \sqrt{3}  \Delta y},0,0,0,-\frac{ \Delta x (t_y-1)}{2 \sqrt{3}  \Delta y},-\frac{ \Delta x (1+ t_y)}{2  \Delta y},-\frac{\sqrt{\frac{5}{3}}  \Delta x (t_y-1)}{2  \Delta y},0,0,0,0,0,1,0,0,0,\boxed{9 \times 0}\right )^\text{T} 
 \end{align*}
\begin{align*}
w_5 = \left (-\frac{ \Delta x (t_y-1)}{2  \Delta y},-\frac{\sqrt{3}  \Delta x (1+ t_y)}{2  \Delta y},-\frac{\sqrt{5}  \Delta x (1+ t_y)}{2  \Delta y},0,0,0,\frac{ \Delta x (t_y-1)}{2 \sqrt{5}  \Delta y},\frac{\sqrt{\frac{3}{5}}  \Delta x (1+ t_y)}{2  \Delta y},\frac{ \Delta x (1+ t_y)}{2  \Delta y},0,0,0,0,1,0,0,0,0,\boxed{9 \times 0}\right )^\text{T} 
 \end{align*}
\begin{align*}
w_6 = \left (\frac{ \Delta x (t_y-1)}{2 \sqrt{3}  \Delta y},\frac{ \Delta x (t_y-1)}{2  \Delta y},\frac{\sqrt{\frac{5}{3}}  \Delta x (t_y-1)}{2  \Delta y},0,0,0,-\frac{ \Delta x (t_y-1)}{2 \sqrt{15}  \Delta y},-\frac{ \Delta x (t_y-1)}{2 \sqrt{5}  \Delta y},-\frac{ \Delta x (t_y-1)}{2 \sqrt{3}  \Delta y},0,0,0,1,0,0,0,0,0,\boxed{9 \times 0}\right )^\text{T} 
 \end{align*}
\begin{align*}
w_7 = \left (-\frac{\sqrt{5}  \Delta x (1+ t_x) (t_y-1)}{2  \Delta y (t_x-1)},-\frac{\sqrt{15}  \Delta x (1+ t_x) (1+ t_y)}{2  \Delta y (t_x-1)},-\frac{5  \Delta x (1+ t_x) (t_y-1)}{2  \Delta y (t_x-1)},-\frac{\sqrt{\frac{5}{3}}  \Delta x (t_y-1)}{2  \Delta y},-\frac{\sqrt{5}  \Delta x (1+ t_y)}{2  \Delta y},-\frac{5  \Delta x (t_y-1)}{2 \sqrt{3}  \Delta y},0,0,0,0,0,1,0,0,0,0,0,0,\boxed{9 \times 0}\right )^\text{T} 
 \end{align*}
\begin{align*}
w_8 = \left (\frac{\sqrt{3}  \Delta x (1+ t_x) (t_y-1)}{2  \Delta y (t_x-1)},\frac{3  \Delta x (1+ t_x) (1+ t_y)}{2  \Delta y (t_x-1)},\frac{\sqrt{15}  \Delta x (1+ t_x) (1+ t_y)}{2  \Delta y (t_x-1)},\frac{ \Delta x (t_y-1)}{2  \Delta y},\frac{\sqrt{3}  \Delta x (1+ t_y)}{2  \Delta y},\frac{\sqrt{5}  \Delta x (1+ t_y)}{2  \Delta y},0,0,0,0,1,0,0,0,0,0,0,0,\boxed{9 \times 0}\right )^\text{T} 
 \end{align*}
\begin{align*}
w_9 = \left (-\frac{ \Delta x (1+ t_x) (t_y-1)}{2  \Delta y (t_x-1)},-\frac{\sqrt{3}  \Delta x (1+ t_x) (t_y-1)}{2  \Delta y (t_x-1)},-\frac{\sqrt{5}  \Delta x (1+ t_x) (t_y-1)}{2  \Delta y (t_x-1)},-\frac{ \Delta x (t_y-1)}{2 \sqrt{3}  \Delta y},-\frac{ \Delta x (t_y-1)}{2  \Delta y},-\frac{\sqrt{\frac{5}{3}}  \Delta x (t_y-1)}{2  \Delta y},0,0,0,1,0,0,0,0,0,0,0,0,\boxed{9 \times 0}\right )^\text{T} 
\end{align*}}

For $K=3$, the kernel of $\mathcal E$ is spanned by these 16 vectors: {\tiny
\begin{align*}
w_1 = \Bigg (0,0,0,0,0,0,0,0,-\frac{ \Delta x (t_y-1)}{2 \sqrt{5}  \Delta y},-\frac{\sqrt{\frac{3}{5}}  \Delta x (1+ t_y)}{2
 \Delta y},-\frac{ \Delta x (t_y-1)}{2  \Delta y},-\frac{\sqrt{\frac{7}{5}}  \Delta x (1+ t_y)}{2  \Delta y},-\frac{ \Delta x (t_y-1)}{2
\sqrt{7}  \Delta y},-\frac{\sqrt{\frac{3}{7}}  \Delta x (1+ t_y)}{2  \Delta y},-\frac{\sqrt{\frac{5}{7}}  \Delta x (t_y-1)}{2  \Delta y},-\frac{ \Delta x
(1+ t_y)}{2  \Delta y}, \\ \nonumber 0,0,0,0,0,0,0,0,0,0,0,0,0,0,0,1,\boxed{16 \times 0}\Bigg )^\text{T} 
 \end{align*}
\begin{align*}
w_2 = \Bigg (0,0,0,0,0,0,0,0,\frac{ \Delta x (t_y-1)}{2
\sqrt{7}  \Delta y},\frac{\sqrt{\frac{3}{7}}  \Delta x (1+ t_y)}{2  \Delta y},\frac{\sqrt{\frac{5}{7}}  \Delta x (t_y-1)}{2  \Delta y},\frac{ \Delta x
(t_y-1)}{2  \Delta y},\frac{\sqrt{5}  \Delta x (t_y-1)}{14  \Delta y},\frac{\sqrt{15}  \Delta x (1+ t_y)}{14  \Delta y},\frac{5
 \Delta x (t_y-1)}{14  \Delta y},\frac{\sqrt{\frac{5}{7}}  \Delta x (t_y-1)}{2  \Delta y},\\ \nonumber 0,0,0,0,0,0,0,0,0,0,0,0,0,0,1,0,\boxed{16 \times 0}\Bigg )^\text{T} 
  \end{align*}
\begin{align*}
w_3 = \Bigg (0,0,0,0,0,0,0,0,-\frac{\sqrt{\frac{3}{35}}
 \Delta x (t_y-1)}{2  \Delta y},-\frac{3  \Delta x (1+ t_y)}{2 \sqrt{35}  \Delta y},-\frac{\sqrt{\frac{3}{7}}  \Delta x (1+ t_y)}{2
 \Delta y},-\frac{\sqrt{\frac{3}{5}}  \Delta x (1+ t_y)}{2  \Delta y},-\frac{\sqrt{3}  \Delta x (t_y-1)}{14  \Delta y},-\frac{3  \Delta x
(1+ t_y)}{14  \Delta y},-\frac{\sqrt{15}  \Delta x (1+ t_y)}{14  \Delta y},-\frac{\sqrt{\frac{3}{7}}  \Delta x (1+ t_y)}{2  \Delta y},\\ \nonumber 0,0,0,0,0,0,0,0,0,0,0,0,0,1,0,0,\boxed{16 \times 0}\Bigg )^\text{T} 
 \end{align*}
\begin{align*}
w_4 = \Bigg (0,0,0,0,0,0,0,0,\frac{ \Delta x
(t_y-1)}{2 \sqrt{35}  \Delta y},\frac{\sqrt{\frac{3}{35}}  \Delta x (t_y-1)}{2  \Delta y},\frac{ \Delta x (t_y-1)}{2 \sqrt{7}
 \Delta y},\frac{ \Delta x (t_y-1)}{2 \sqrt{5}  \Delta y},\frac{ \Delta x (t_y-1)}{14  \Delta y},\frac{\sqrt{3}  \Delta x (t_y-1)}{14
 \Delta y},\frac{\sqrt{5}  \Delta x (t_y-1)}{14  \Delta y},\frac{ \Delta x (t_y-1)}{2 \sqrt{7}  \Delta y},\\ \nonumber 0,0,0,0,0,0,0,0,0,0,0,0,1,0,0,0,\boxed{16 \times 0}\Bigg )^\text{T} 
  \end{align*}
\begin{align*}
w_5 = \Bigg (0,0,0,0,-\frac{\sqrt{\frac{7}{15}}
 \Delta x (t_y-1)}{2  \Delta y},-\frac{\sqrt{\frac{7}{5}}  \Delta x (1+ t_y)}{2  \Delta y},-\frac{\sqrt{\frac{7}{3}}  \Delta x (t_y-1)}{2
 \Delta y},-\frac{7  \Delta x (1+ t_y)}{2 \sqrt{15}  \Delta y},0,0,0,0,\frac{ \Delta x (t_y-1)}{2 \sqrt{5}  \Delta y},\frac{\sqrt{\frac{3}{5}}
 \Delta x (1+ t_y)}{2  \Delta y},\frac{ \Delta x (t_y-1)}{2  \Delta y},\frac{\sqrt{\frac{7}{5}}  \Delta x (1+ t_y)}{2  \Delta y},\\ \nonumber 0,0,0,0,0,0,0,0,0,0,0,1,0,0,0,0,\boxed{16 \times 0}\Bigg )^\text{T} 
  \end{align*}
\begin{align*}
w_6 = \Bigg (0,0,0,0,\frac{ \Delta x
(t_y-1)}{2 \sqrt{3}  \Delta y},\frac{ \Delta x (1+ t_y)}{2  \Delta y},\frac{\sqrt{\frac{5}{3}}  \Delta x (t_y-1)}{2  \Delta y},\frac{\sqrt{\frac{7}{3}}
 \Delta x (t_y-1)}{2  \Delta y},0,0,0,0,-\frac{ \Delta x (t_y-1)}{2 \sqrt{7}  \Delta y},-\frac{\sqrt{\frac{3}{7}}  \Delta x (1+ t_y)}{2
 \Delta y},-\frac{\sqrt{\frac{5}{7}}  \Delta x (t_y-1)}{2  \Delta y},-\frac{ \Delta x (t_y-1)}{2  \Delta y},\\ \nonumber 0,0,0,0,0,0,0,0,0,0,1,0,0,0,0,0,\boxed{16 \times 0}\Bigg )^\text{T} 
  \end{align*}
\begin{align*}
w_7 = \Bigg (0,0,0,0,-\frac{ \Delta x
(t_y-1)}{2 \sqrt{5}  \Delta y},-\frac{\sqrt{\frac{3}{5}}  \Delta x (1+ t_y)}{2  \Delta y},-\frac{ \Delta x (1+ t_y)}{2  \Delta y},-\frac{\sqrt{\frac{7}{5}}
 \Delta x (1+ t_y)}{2  \Delta y},0,0,0,0,\frac{\sqrt{\frac{3}{35}}  \Delta x (t_y-1)}{2  \Delta y},\frac{3  \Delta x (1+ t_y)}{2
\sqrt{35}  \Delta y},\frac{\sqrt{\frac{3}{7}}  \Delta x (1+ t_y)}{2  \Delta y},\frac{\sqrt{\frac{3}{5}}  \Delta x (1+ t_y)}{2  \Delta y},\\ \nonumber 0,0,0,0,0,0,0,0,0,1,0,0,0,0,0,0,\boxed{16 \times 0}\Bigg )^\text{T} 
 \end{align*}
\begin{align*}
w_8 = \Bigg (0,0,0,0,\frac{ \Delta x
(t_y-1)}{2 \sqrt{15}  \Delta y},\frac{ \Delta x (t_y-1)}{2 \sqrt{5}  \Delta y},\frac{ \Delta x (t_y-1)}{2 \sqrt{3}  \Delta y},\frac{\sqrt{\frac{7}{15}}
 \Delta x (t_y-1)}{2  \Delta y},0,0,0,0,-\frac{ \Delta x (t_y-1)}{2 \sqrt{35}  \Delta y},-\frac{\sqrt{\frac{3}{35}}  \Delta x (t_y-1)}{2
 \Delta y},-\frac{ \Delta x (t_y-1)}{2 \sqrt{7}  \Delta y},-\frac{ \Delta x (t_y-1)}{2 \sqrt{5}  \Delta y},\\ \nonumber 0,0,0,0,0,0,0,0,1,0,0,0,0,0,0,0,\boxed{16 \times 0}\Bigg )^\text{T}
  \end{align*}
\begin{align*}
w_9 = \Bigg (-\frac{\sqrt{\frac{7}{3}}
 \Delta x (t_y-1)}{2  \Delta y},-\frac{\sqrt{7}  \Delta x (1+ t_y)}{2  \Delta y},-\frac{\sqrt{\frac{35}{3}}  \Delta x (t_y-1)}{2
 \Delta y},-\frac{7  \Delta x (1+ t_y)}{2 \sqrt{3}  \Delta y},0,0,0,0,\frac{\sqrt{\frac{7}{15}}  \Delta x (t_y-1)}{2  \Delta y},\frac{\sqrt{\frac{7}{5}}
 \Delta x (1+ t_y)}{2  \Delta y},\frac{\sqrt{\frac{7}{3}}  \Delta x (t_y-1)}{2  \Delta y},\frac{7  \Delta x (1+ t_y)}{2 \sqrt{15}
 \Delta y},0,0,0,0,\\ \nonumber 0,0,0,0,0,0,0,1,0,0,0,0,0,0,0,0,\boxed{16 \times 0}\Bigg )^\text{T} 
  \end{align*}
\begin{align*}
w_{10} = \Bigg (\frac{\sqrt{\frac{5}{3}}  \Delta x (t_y-1)}{2
 \Delta y},\frac{\sqrt{5}  \Delta x (1+ t_y)}{2  \Delta y},\frac{5  \Delta x (t_y-1)}{2 \sqrt{3}  \Delta y},\frac{\sqrt{\frac{35}{3}}
 \Delta x (t_y-1)}{2  \Delta y},0,0,0,0,-\frac{ \Delta x (t_y-1)}{2 \sqrt{3}  \Delta y},-\frac{ \Delta x (1+ t_y)}{2  \Delta y},-\frac{\sqrt{\frac{5}{3}}
 \Delta x (t_y-1)}{2  \Delta y},-\frac{\sqrt{\frac{7}{3}}  \Delta x (t_y-1)}{2  \Delta y},0,0,0,0,\\ \nonumber 0,0,0,0,0,0,1,0,0,0,0,0,0,0,0,0,\boxed{16 \times 0}\Bigg )^\text{T} 
  \end{align*}
\begin{align*}
w_{11} = \Bigg (-\frac{ \Delta x
(t_y-1)}{2  \Delta y},-\frac{\sqrt{3}  \Delta x (1+ t_y)}{2  \Delta y},-\frac{\sqrt{5}  \Delta x (1+ t_y)}{2  \Delta y},-\frac{\sqrt{7}
 \Delta x (1+ t_y)}{2  \Delta y},0,0,0,0,\frac{ \Delta x (t_y-1)}{2 \sqrt{5}  \Delta y},\frac{\sqrt{\frac{3}{5}}  \Delta x (1+ t_y)}{2
 \Delta y},\frac{ \Delta x (1+ t_y)}{2  \Delta y},\frac{\sqrt{\frac{7}{5}}  \Delta x (1+ t_y)}{2  \Delta y},0,0,0,0,\\ \nonumber 0,0,0,0,0,1,0,0,0,0,0,0,0,0,0,0,\boxed{16 \times 0}\Bigg )^\text{T} 
  \end{align*}
\begin{align*}
w_{12} = \Bigg (\frac{ \Delta x
(t_y-1)}{2 \sqrt{3}  \Delta y},\frac{ \Delta x (t_y-1)}{2  \Delta y},\frac{\sqrt{\frac{5}{3}}  \Delta x (t_y-1)}{2  \Delta y},\frac{\sqrt{\frac{7}{3}}
 \Delta x (t_y-1)}{2  \Delta y},0,0,0,0,-\frac{ \Delta x (t_y-1)}{2 \sqrt{15}  \Delta y},-\frac{ \Delta x (t_y-1)}{2 \sqrt{5}
 \Delta y},-\frac{ \Delta x (t_y-1)}{2 \sqrt{3}  \Delta y},-\frac{\sqrt{\frac{7}{15}}  \Delta x (t_y-1)}{2  \Delta y},0,0,0,0,\\ \nonumber 0,0,0,0,1,0,0,0,0,0,0,0,0,0,0,0,\boxed{16 \times 0}\Bigg )^\text{T} 
  \end{align*}
\begin{align*}
w_{13} = \Bigg (\frac{\sqrt{7}
 \Delta x (1+ t_x) (t_y-1)}{2  \Delta y (t_x-1)},\frac{\sqrt{21}  \Delta x (1+ t_x) (1+ t_y)}{2  \Delta y (t_x-1)},\frac{\sqrt{35}
 \Delta x (1+ t_x) (t_y-1)}{2  \Delta y (t_x-1)},\frac{7  \Delta x (1+ t_x) (1+ t_y)}{2  \Delta y (t_x-1)},\frac{\sqrt{\frac{7}{3}}
 \Delta x (t_y-1)}{2  \Delta y},\frac{\sqrt{7}  \Delta x (1+ t_y)}{2  \Delta y},\frac{\sqrt{\frac{35}{3}}  \Delta x (t_y-1)}{2
 \Delta y},\frac{7  \Delta x (1+ t_y)}{2 \sqrt{3}  \Delta y},0,0,0,0,0,0,0,0,\\ \nonumber 0,0,0,1,0,0,0,0,0,0,0,0,0,0,0,0,\boxed{16 \times 0}\Bigg )^\text{T} 
  \end{align*}
\begin{align*}
w_{14} = \Bigg (-\frac{\sqrt{5}
 \Delta x (1+ t_x) (t_y-1)}{2  \Delta y (t_x-1)},-\frac{\sqrt{15}  \Delta x (1+ t_x) (1+ t_y)}{2  \Delta y (t_x-1)},-\frac{5
 \Delta x (1+ t_x) (t_y-1)}{2  \Delta y (t_x-1)},-\frac{\sqrt{35}  \Delta x (1+ t_x) (t_y-1)}{2  \Delta y (t_x-1)},-\frac{\sqrt{\frac{5}{3}}
 \Delta x (t_y-1)}{2  \Delta y},-\frac{\sqrt{5}  \Delta x (1+ t_y)}{2  \Delta y},-\frac{5  \Delta x (t_y-1)}{2 \sqrt{3}  \Delta y},-\frac{\sqrt{\frac{35}{3}}
 \Delta x (t_y-1)}{2  \Delta y},0,0,0,0,0,0,0,0,\\ \nonumber 0,0,1,0,0,0,0,0,0,0,0,0,0,0,0,0,\boxed{16 \times 0}\Bigg )^\text{T} 
  \end{align*}
\begin{align*}
w_{15} = \Bigg (\frac{\sqrt{3}
 \Delta x (1+ t_x) (t_y-1)}{2  \Delta y (t_x-1)},\frac{3  \Delta x (1+ t_x) (1+ t_y)}{2  \Delta y (t_x-1)},\frac{\sqrt{15}
 \Delta x (1+ t_x) (1+ t_y)}{2  \Delta y (t_x-1)},\frac{\sqrt{21}  \Delta x (1+ t_x) (1+ t_y)}{2  \Delta y (t_x-1)},\frac{ \Delta x
(t_y-1)}{2  \Delta y},\frac{\sqrt{3}  \Delta x (1+ t_y)}{2  \Delta y},\frac{\sqrt{5}  \Delta x (1+ t_y)}{2  \Delta y},\frac{\sqrt{7}
 \Delta x (1+ t_y)}{2  \Delta y},0,0,0,0,0,0,0,0,\\ \nonumber 0,1,0,0,0,0,0,0,0,0,0,0,0,0,0,0,\boxed{16 \times 0}\Bigg )^\text{T} 
  \end{align*}
\begin{align*}
w_{16} = \Bigg (-\frac{ \Delta x
(1+ t_x) (t_y-1)}{2  \Delta y (t_x-1)},-\frac{\sqrt{3}  \Delta x (1+ t_x) (t_y-1)}{2  \Delta y (t_x-1)},-\frac{\sqrt{5}
 \Delta x (1+ t_x) (t_y-1)}{2  \Delta y (t_x-1)},-\frac{\sqrt{7}  \Delta x (1+ t_x) (t_y-1)}{2  \Delta y (t_x-1)}, -\frac{ \Delta x
(t_y-1)}{2 \sqrt{3}  \Delta y},-\frac{ \Delta x (t_y-1)}{2  \Delta y},-\frac{\sqrt{\frac{5}{3}}  \Delta x (t_y-1)}{2  \Delta y},-\frac{\sqrt{\frac{7}{3}}
 \Delta x (t_y-1)}{2  \Delta y},0,0,0,0,0,0,0,0,\\ \nonumber 1,0,0,0,0,0,0,0,0,0,0,0,0,0,0,0,\boxed{16 \times 0}\Bigg )^\text{T} 
\end{align*}}
\end{landscape}


\begin{thebibliography}{NRDB{\etalchar{+}}13}

\bibitem[AIP19]{abbate19}
Emanuela Abbate, Angelo Iollo, and Gabriella Puppo.
\newblock An asymptotic-preserving all-speed scheme for fluid dynamics and
  nonlinear elasticity.
\newblock {\em SIAM Journal on Scientific Computing}, 41(5):A2850--A2879, 2019.

\bibitem[Bar18]{barsukow18thesis}
Wasilij Barsukow.
\newblock {\em Low {M}ach number finite volume methods for the acoustic and
  {E}uler equations}.
\newblock Doctoral thesis, University of Wuerzburg, 2018.

\bibitem[Bar19]{barsukow17a}
Wasilij Barsukow.
\newblock Stationarity preserving schemes for multi-dimensional linear systems.
\newblock {\em Mathematics of Computation}, 88(318):1621--1645, 2019.

\bibitem[Bar21]{barsukow20cgk}
Wasilij Barsukow.
\newblock Truly multi-dimensional all-speed schemes for the euler equations on
  cartesian grids.
\newblock {\em Journal of Computational Physics}, 435:110216, 2021.

\bibitem[Bar23]{barsukow21yee}
Wasilij Barsukow.
\newblock All-speed numerical methods for the {E}uler equations via a
  sequential explicit time integration.
\newblock {\em Journal of Scientific Computing}, 95(2):53, 2023.

\bibitem[BCG20]{bouchut20}
Fran{\c{c}}ois Bouchut, Christophe Chalons, and S{\'e}bastien Guisset.
\newblock An entropy satisfying two-speed relaxation system for the barotropic
  euler equations: application to the numerical approximation of low mach
  number flows.
\newblock {\em Numerische Mathematik}, 145(1):35--76, 2020.

\bibitem[BDBHN09]{bassi09}
Francesco Bassi, Carmine De~Bartolo, Ralf Hartmann, and Alessandra Nigro.
\newblock A discontinuous galerkin method for inviscid low mach number flows.
\newblock {\em Journal of Computational Physics}, 228(11):3996--4011, 2009.

\bibitem[BDL{\etalchar{+}}20]{boscheri20}
Walter Boscheri, Giacomo Dimarco, Rapha{\"e}l Loub{\`e}re, Maurizio Tavelli,
  and Marie-H{\'e}l{\`e}ne Vignal.
\newblock A second order all mach number imex finite volume solver for the
  three dimensional euler equations.
\newblock {\em Journal of Computational Physics}, 415:109486, 2020.

\bibitem[BDT21]{boscheri21}
Walter Boscheri, Giacomo Dimarco, and Maurizio Tavelli.
\newblock An efficient second order all {M}ach finite volume solver for the
  compressible {N}avier-{S}tokes equations.
\newblock {\em Computer Methods in Applied Mechanics and Engineering},
  374:113602, 2021.

\bibitem[BEK{\etalchar{+}}17]{barsukow16}
Wasilij Barsukow, Philipp~VF Edelmann, Christian Klingenberg, Fabian Miczek,
  and Friedrich~K R{\"o}pke.
\newblock A numerical scheme for the compressible low-{M}ach number regime of
  ideal fluid dynamics.
\newblock {\em Journal of Scientific Computing}, 72(2):623--646, 2017.

\bibitem[BHKR19]{barsukow18activeflux}
Wasilij Barsukow, Jonathan Hohm, Christian Klingenberg, and Philip~L Roe.
\newblock The active flux scheme on {C}artesian grids and its low {M}ach number
  limit.
\newblock {\em Journal of Scientific Computing}, 81(1):594--622, 2019.

\bibitem[BKKL25]{barsukow24affourier}
Wasilij Barsukow, Janina Kern, Christian Klingenberg, and Lisa Lechner.
\newblock Analysis of the multi-dimensional semi-discrete {A}ctive {F}lux
  method using the fourier transform.
\newblock {\em Communications on Applied Mathematics and Computation}, pages
  1--49, 2025.

\bibitem[BKL{\etalchar{+}}25]{barsukow2025sbp}
Wasilij Barsukow, Christian Klingenberg, Lisa Lechner, Jan Nordstr{\"o}m,
  Sigrun Ortleb, and Hendrik Ranocha.
\newblock Proving stability of the active flux method in the framework of
  summation-by-parts operators.
\newblock {\em submitted}, 2025.

\bibitem[BLM23]{barsukow23nodal}
Wasilij Barsukow, Raphael Loubere, and Pierre-Henri Maire.
\newblock A node-conservative vorticity-preserving finite volume method for
  linear acoustics on unstructured grids.
\newblock {\em accepted in Math.Comp.}, 2023.

\bibitem[BM05]{birken05}
Philipp Birken and Andreas Meister.
\newblock Stability of preconditioned finite volume schemes at low {M}ach
  numbers.
\newblock {\em BIT Numerical Mathematics}, 45(3):463--480, 2005.

\bibitem[BP21]{boscheri21a}
Walter Boscheri and Lorenzo Pareschi.
\newblock High order pressure-based semi-implicit {IMEX} schemes for the 3{D}
  {N}avier-{S}tokes equations at all {M}ach numbers.
\newblock {\em Journal of Computational Physics}, 434:110206, 2021.

\bibitem[BQRX19]{boscarino19}
Sebastiano Boscarino, Jing-Mei Qiu, Giovanni Russo, and Tao Xiong.
\newblock A high order semi-implicit {IMEX WENO} scheme for the all-{M}ach
  isentropic {E}uler system.
\newblock {\em Journal of Computational Physics}, 392:594--618, 2019.

\bibitem[BR97]{bassi97}
Francesco Bassi and Stefano Rebay.
\newblock A high-order accurate discontinuous finite element method for the
  numerical solution of the compressible navier--stokes equations.
\newblock {\em Journal of computational physics}, 131(2):267--279, 1997.

\bibitem[CDK12]{cordier12}
Floraine Cordier, Pierre Degond, and Anela Kumbaro.
\newblock An asymptotic-preserving all-speed scheme for the {E}uler and
  {N}avier-{S}tokes equations.
\newblock {\em Journal of Computational Physics}, 231(17):5685--5704, 2012.

\bibitem[CGK13]{chalons13}
Christophe Chalons, Mathieu Girardin, and Samuel Kokh.
\newblock Large time step and asymptotic preserving numerical schemes for the
  gas dynamics equations with source terms.
\newblock {\em SIAM Journal on Scientific Computing}, 35(6):A2874--A2902, 2013.

\bibitem[CKS12]{cockburn12}
Bernardo Cockburn, George~E Karniadakis, and Chi-Wang Shu.
\newblock {\em Discontinuous Galerkin methods: theory, computation and
  applications}, volume~11.
\newblock Springer Science \& Business Media, 2012.

\bibitem[CS89]{cockburn1989a}
Bernardo Cockburn and Chi-Wang Shu.
\newblock {TVB} {R}unge-{K}utta local projection discontinuous {G}alerkin
  finite element method for conservation laws {II}: General framework.
\newblock {\em Mathematics of computation}, 52(186):411--435, 1989.

\bibitem[Del10]{dellacherie10}
St{\'e}phane Dellacherie.
\newblock Analysis of {G}odunov type schemes applied to the compressible
  {E}uler system at low {M}ach number.
\newblock {\em Journal of Computational Physics}, 229(4):978--1016, 2010.

\bibitem[DJOR16]{dellacherie16}
St{\'e}phane Dellacherie, Jonathan Jung, Pascal Omnes, and P-A Raviart.
\newblock Construction of modified {G}odunov-type schemes accurate at any
  {M}ach number for the compressible {E}uler system.
\newblock {\em Mathematical Models and Methods in Applied Sciences},
  26(13):2525--2615, 2016.

\bibitem[DJY07]{degond07}
Pierre Degond, S~Jin, and J~Yuming.
\newblock {M}ach-number uniform asymptotic-preserving gauge schemes for
  compressible flows.
\newblock {\em Bulletin-Institute of Mathematics Academia Sinica}, 2(4):851,
  2007.

\bibitem[DLV17]{dimarco17}
Giacomo Dimarco, Rapha{\"e}l Loub{\`e}re, and Marie-H{\'e}l{\`e}ne Vignal.
\newblock Study of a new asymptotic preserving scheme for the {E}uler system in
  the low {M}ach number limit.
\newblock {\em SIAM Journal on Scientific Computing}, 39(5):A2099--A2128, 2017.

\bibitem[DOR10]{dellacherierieper10}
St{\'e}phane Dellacherie, Pascal Omnes, and Felix Rieper.
\newblock The influence of cell geometry on the {G}odunov scheme applied to the
  linear wave equation.
\newblock {\em Journal of Computational Physics}, 229(14):5315--5338, 2010.

\bibitem[EG23]{ern23}
Alexandre Ern and Jean-Luc Guermond.
\newblock The discontinuous galerkin approximation of the grad-div and
  curl-curl operators in first-order form is involution-preserving and
  spectrally correct.
\newblock {\em SIAM Journal on Numerical Analysis}, 61(6):2940--2966, 2023.

\bibitem[FK07]{feistauer07}
Miloslav Feistauer and V~Ku{\v{c}}era.
\newblock On a robust discontinuous galerkin technique for the solution of
  compressible flow.
\newblock {\em Journal of Computational Physics}, 224(1):208--221, 2007.

\bibitem[GC90]{gresho90}
Philip~M Gresho and Stevens~T Chan.
\newblock On the theory of semi-implicit projection methods for viscous
  incompressible flow and its implementation via a finite element method that
  also introduces a nearly consistent mass matrix. {P}art 2: {I}mplementation.
\newblock {\em International Journal for Numerical Methods in Fluids},
  11(5):621--659, 1990.

\bibitem[GLW07]{gjonaj07}
E~Gjonaj, T~Lau, and T~Weiland.
\newblock Conservation properties of the discontinuous {G}alerkin method for
  {M}axwell equations.
\newblock In {\em 2007 International Conference on Electromagnetics in Advanced
  Applications}, pages 356--359. IEEE, 2007.

\bibitem[GM04]{guillard04}
Herv{\'e} Guillard and Angelo Murrone.
\newblock On the behavior of upwind schemes in the low {M}ach number limit:
  {II}. {G}odunov type schemes.
\newblock {\em Computers \& fluids}, 33(4):655--675, 2004.

\bibitem[Gui09]{guillard09}
Herv{\'e} Guillard.
\newblock On the behavior of upwind schemes in the low {M}ach number limit.
  {IV}: {P}0 approximation on triangular and tetrahedral cells.
\newblock {\em Computers \& Fluids}, 38(10):1969--1972, 2009.

\bibitem[HG20]{hindenlang20}
Florian~J Hindenlang and Gregor~J Gassner.
\newblock On the order reduction of entropy stable {DGSEM} for the compressible
  {E}uler equations.
\newblock In {\em Spectral and High Order Methods for Partial Differential
  Equations ICOSAHOM 2018: Selected Papers from the ICOSAHOM Conference,
  London, UK, July 9-13, 2018}, pages 21--44. Springer International
  Publishing, 2020.

\bibitem[HJL12]{haack12}
Jeffrey Haack, Shi Jin, and Jian-Guo Liu.
\newblock An all-speed asymptotic-preserving method for the isentropic {E}uler
  and {N}avier-{S}tokes equations.
\newblock {\em Communications in Computational Physics}, 12(4):955--980, 2012.

\bibitem[HTL21]{hennink21}
Aldo Hennink, Marco Tiberga, and Danny Lathouwers.
\newblock A pressure-based solver for low-mach number flow using a
  discontinuous galerkin method.
\newblock {\em Journal of Computational Physics}, 425:109877, 2021.

\bibitem[JP22]{jung22}
Jonathan Jung and Vincent Perrier.
\newblock Steady low mach number flows: identification of the spurious mode and
  filtering method.
\newblock {\em Journal of Computational Physics}, 468:111462, 2022.

\bibitem[JT06]{jeltsch06}
Rolf Jeltsch and Manuel Torrilhon.
\newblock On curl-preserving finite volume discretizations for shallow water
  equations.
\newblock {\em BIT Numerical Mathematics}, 46(1):35--53, 2006.

\bibitem[Kle95]{klein95}
Rupert Klein.
\newblock Semi-implicit extension of a {G}odunov-type scheme based on low
  {M}ach number asymptotics {I}: {O}ne-dimensional flow.
\newblock {\em Journal of Computational Physics}, 121(2):213--237, 1995.

\bibitem[Lan24]{lannabi24}
Ibtissem Lannabi.
\newblock {\em Analysis of spurious oscillations problem of Finite Volume
  Methods for low Mach number flows in fluid mechanics}.
\newblock PhD thesis, Universit{\'e} de Pau et des Pays de l'Adour, 2024.

\bibitem[LG08]{li08}
Xue-song Li and Chun-wei Gu.
\newblock An all-speed {R}oe-type scheme and its asymptotic analysis of low
  {M}ach number behaviour.
\newblock {\em Journal of Computational Physics}, 227(10):5144--5159, 2008.

\bibitem[LG13]{li13}
Xue-song Li and Chun-wei Gu.
\newblock Mechanism of {R}oe-type schemes for all-speed flows and its
  application.
\newblock {\em Computers \& Fluids}, 86:56--70, 2013.

\bibitem[LR14]{lung14}
TB~Lung and PL~Roe.
\newblock Toward a reduction of mesh imprinting.
\newblock {\em International Journal for Numerical Methods in Fluids},
  76(7):450--470, 2014.

\bibitem[LSZ20]{liu20}
Yong Liu, Chi-Wang Shu, and Mengping Zhang.
\newblock Sub-optimal convergence of discontinuous galerkin methods with
  central fluxes for linear hyperbolic equations with even degree polynomial
  approximations.
\newblock {\em arXiv preprint arXiv:2001.03825}, 2020.

\bibitem[MR01]{morton01}
Keith~William Morton and Philip~L Roe.
\newblock Vorticity-preserving {L}ax-{W}endroff-type schemes for the system
  wave equation.
\newblock {\em SIAM Journal on Scientific Computing}, 23(1):170--192, 2001.

\bibitem[MRE15]{miczek15}
F~Miczek, FK~R{\"o}pke, and PVF Edelmann.
\newblock New numerical solver for flows at various {M}ach numbers.
\newblock {\em Astronomy \& Astrophysics}, 576:A50, 2015.

\bibitem[MT09]{mishra09preprint}
Siddhartha Mishra and Eitan Tadmor.
\newblock Constraint preserving schemes using potential-based fluxes {II}.
  genuinely multi-dimensional central schemes for systems of conservation laws.
\newblock {\em ETH preprint}, (2009-32), 2009.

\bibitem[NRDB{\etalchar{+}}13]{bassi13}
Alessandra Nigro, Salvatore Renda, Carmine De~Bartolo, Ralf Hartmann, and
  Francesco Bassi.
\newblock A high-order accurate discontinuous galerkin finite element method
  for laminar low mach number flows.
\newblock {\em International Journal for Numerical Methods in Fluids},
  72(1):43--68, 2013.

\bibitem[OSB{\etalchar{+}}16]{birken16}
Kai O{\ss}wald, Alexander Siegmund, Philipp Birken, Volker Hannemann, and
  Andreas Meister.
\newblock L2roe: a low dissipation version of {R}oe's approximate {R}iemann
  solver for low {M}ach numbers.
\newblock {\em International Journal for Numerical Methods in Fluids},
  81(2):71--86, 2016.

\bibitem[Per24a]{perrier24dg}
Vincent Perrier.
\newblock Development of discontinuous galerkin methods for hyperbolic systems
  that preserve a curl or a divergence constraint.
\newblock {\em arXiv preprint arXiv:2405.04347}, 2024.

\bibitem[Per24b]{perrier24}
Vincent Perrier.
\newblock Discrete de-{R}ham complex involving a discontinuous finite element
  space for velocities: the case of periodic straight triangular and
  {C}artesian meshes.
\newblock {\em arXiv preprint arXiv:2404.19545}, 2024.

\bibitem[PP08]{persson08}
P-O Persson and Jaime Peraire.
\newblock Newton-gmres preconditioning for discontinuous galerkin
  discretizations of the navier--stokes equations.
\newblock {\em SIAM Journal on Scientific Computing}, 30(6):2709--2733, 2008.

\bibitem[Rie11]{rieper11}
Felix Rieper.
\newblock A low-{M}ach number fix for {R}oe’s approximate {R}iemann solver.
\newblock {\em Journal of Computational Physics}, 230(13):5263--5287, 2011.

\bibitem[Sid02]{sidilkover02}
David Sidilkover.
\newblock Factorizable schemes for the equations of fluid flow.
\newblock {\em Applied numerical mathematics}, 41(3):423--436, 2002.

\bibitem[TD08]{thornber08}
BJR Thornber and D~Drikakis.
\newblock Numerical dissipation of upwind schemes in low {M}ach flow.
\newblock {\em International journal for numerical methods in fluids},
  56(8):1535--1541, 2008.

\bibitem[TPK20]{thomann20}
Andrea Thomann, Gabriella Puppo, and Christian Klingenberg.
\newblock An all speed second order well-balanced imex relaxation scheme for
  the euler equations with gravity.
\newblock {\em Journal of Computational Physics}, 420:109723, 2020.

\bibitem[Tur87]{turkel87}
Eli Turkel.
\newblock Preconditioned methods for solving the incompressible and low speed
  compressible equations.
\newblock {\em Journal of computational physics}, 72(2):277--298, 1987.

\bibitem[VBW11]{viallet11}
M~Viallet, I~Baraffe, and R~Walder.
\newblock Towards a new generation of multi-dimensional stellar evolution
  models: development of an implicit hydrodynamic code.
\newblock {\em Astronomy \& Astrophysics}, 531:A86, 2011.

\bibitem[WS95]{weiss95}
Jonathan~M Weiss and Wayne~A Smith.
\newblock Preconditioning applied to variable and constant density flows.
\newblock {\em AIAA journal}, 33(11):2050--2057, 1995.

\end{thebibliography}
\end{document}